\documentclass[12pt]{article}

\usepackage{amsthm,amsmath,stmaryrd,bbm,hyperref,geometry,color}
\usepackage{amssymb}
\usepackage[english]{babel}
 \usepackage[utf8]{inputenc}
\usepackage{graphicx}
\usepackage{verbatim}
\usepackage{enumitem}

\newcommand{\po}{\left(}
\newcommand{\pf}{\right)}
\newcommand{\co}{\left[}
\newcommand{\cf}{\right]}
\newcommand{\cco}{\llbracket}
\newcommand{\ccf}{\rrbracket}
\newcommand{\R}{\mathbb R} 
 
\newcommand{\Z}{\mathbb Z} 
 
\newcommand{\N}{\mathbb N} 
\newcommand{\dd}{\text{d}}
\newcommand{\na}{\nabla}
\newcommand{\1}{\mathbbm{1}}
\newcommand{\V}{\mathcal V}
\newcommand{\W}{\mathcal W}
\newcommand{\dist}{\bar d_1}

\newcommand{\PV}{\mathcal P_{\V,1}(E)}
\newcommand{\MRP}{\mathcal M(\R_+,\mathcal P(E))}

\newtheorem{thrm}{Theorem}
\newtheorem{prpstn}{Proposition}

\newtheorem{lmm}{Lemma}

\newtheorem{dfntn}{Definition}

\newtheorem{assu}{Assumption}
\newtheorem*{example}{Example}

\title{Elementary coupling approach for non-linear perturbation of Markov processes with mean-field jump mechanisms and related problems}
\author{Pierre Monmarché}

\begin{document}
%
\maketitle

\begin{abstract} Mean-field integro-differential equations are studied in an abstract framework, through couplings of the corresponding stochastic processes. In the perturbative regime, the equation is proven to admit a unique equilibrium, toward which the process converges exponentially fast. Similarly, in this case, the associated particle system is proven to converge toward its equilibrium at a rate independent from the number of particles. \end{abstract}
%
%

\maketitle

\section{Introduction}

The initial motivation of this work is the study of the long-time behaviour of mean-field semi-linear integro-differential equations of the form
\begin{eqnarray}\label{Eq-EquaNonLin}
\partial_t m_t(x) & = & L'm_t(x) + Q'_{m_t}\po \lambda_{m_t} m_t\pf(x) - \lambda_{m_t}(x) m_t(x)\,,
\end{eqnarray}
where $m_t(x)$ is the density of particles at time $t$ at point $x$, $L$ is the generator of a Markov process, $A'$ denotes the dual on measures of an operator $A$ on functions, and for all probability distribution $\nu$, $\lambda_{\nu}>0$ is a jump rate and $Q_{\nu}$ is a jump kernel. In other words, the stochastic process whose law $m_t$ solves \eqref{Eq-EquaNonLin} follows the dynamics of $L$ and, at a rate $\lambda_{m_t}$, jumps to a new position drawn according to the kernel $Q_{m_t}$.

In particular, this encompasses the case of mean-field continuous-time Markov chains or of mean field Piecewise Deterministic Markov Processes (PDMP). Mixed jump/diffusion procceses are also concerned, but only when the non-linearity lies in the jump mechanism. For instance, the kinetic neuron model studied in \cite{PerthameSalort2013}, for which the mean-field interaction lies in the drift, doesn't enter this scope.

Equation \eqref{Eq-EquaNonLin} can model various phenomena like chemical reactions \cite{Schlogl1972}, neuron activity \cite{PerthameSalortPakdaman,CanizoYoldas}, economical games \cite{Dumitrescu2017,DuncanTembine2018}, communication networks \cite{Borkar2012}. In the field of stochastic algorithms, where PDMP samplers have recently gained interest \cite{MonmarcheRTP,DoucetPDMCMC,BierkensFearnheadRoberts}, mean-field equations describe the limiting behaviour, as the number of particles go to infinity, of algorithms like parallel tempering \cite{Swendsen1986}, adaptive biasing force methods \cite{LelievreABF}, selection/mutation techniques like Kalman-Bucy filter \cite{DelMoral2013,DelMoralMiclo}, Fleming-Viot algorithms \cite{FlemingViot79,CloezThai}, adaptive multi-level splitting \cite{CerouGuyader} and so on. 

We are interested in the long-time behaviour of the solutions of \eqref{Eq-EquaNonLin} in the perturbative regime, namely when the effects of the interaction is small with respect to the mixing effect of $L$.
The general idea is simple and based on coupling methods. A classical way to obtain some contraction in (weighted) total variation distance for Markov processes is to define simultaneously two processes starting at different positions in such a way that they have merged at some positive time with some positive probability. Once they have met, in the classical setting, then they remain equal forever. However, in the non-linear case, even if we manage to merge two processes that started at different positions,  they still don't have exactly the same law, hence the same dynamics, and thus they may have to split after some time. Nevertheless, as they stay  together, their laws get  closer one to the other, which makes a splitting less and less likely. Provided the non-linearity is small enough, the balance between this two merge/split mechanisms is  on the side of contraction.

This strategy of coupling two processes that don't have the same dynamics in order to keep them as close  as possible is classical in the study of mean-field diffusions like the McKean-Vlasov equation \cite{Malrieu2001,Sznitman}. In that case, a parallel coupling can typically be used, in other words the processes are defined as solutions of different SDE but driven by the same Brownian motion. Even if they start at the same position, both processes then drift away one from the other, but the evolution of the expectation of their squared distance is bounded by the difference in their drifts which, for mean-field processes, is typically controlled by the Wasserstein-2 distance between their law, which is itself controlled by the expectation of their squared distance. The expansion of the latter is thus controlled through the Gronwall Lemma{{, see e.g. \cite{Malrieu2001,Sznitman,MonmarcheKineticMF} for details.}}

When the interaction lies in the jump mechanism rather than in the deterministic drift one,  although jump processes may lack the regularization properties of diffusions,  the picture is somehow a bit simpler. Indeed, the analogous of the parallel coupling in that case means that, as much as possible, the two processes are prescribed to jump at the same times and to the same locations, and in particular, starting from the same position, they will stay equal for a (random) positive time. As a consequence, total variation estimates are obtained, which are then compatible with classical contraction estimates obtained via Foster-Lyapunov techniques. This synchronous coupling is used in \cite{MonmarcheSIVJP} to study self-interacting PDMP and in \cite{DGM1} to control the difference between   a smooth approximation of a PDMP and the latter. 

This coupling method works both for the non-linear process and for its associated interacting particle system. Actually, the general idea is not restricted to mean-field interactions, and can be applied for instance to self-interacting processes with memory (see e.g. the discussion in Section \ref{Sec-Delay}), sampling processes when the target measure is close to be the product of its marginal (see Section \ref{Sec-ZZ}), chain of oscillators, diffusions with colored noise, killed or branching processes, SPDEs etc. It works for discrete time interacting chains. Mattingly in \cite{Mattingly2002,Mattingly2003,Mattingly2005} and Cloez and Thai in \cite{CloezThai} use the same kind of arguments. For diffusion processes, rather than the simple parallel coupling, the analoguous of our merge/split coupling should be the sticky coupling of Eberle and his co-authors \cite{DurmusEberleGuillin,Eberle,EberleZimmer}. In a word, it is a very natural and general strategy, and we expect the arguments developed in the particular framework of this paper to be of interest in many other settings\footnote{{In fact, between the submission of this work and its publication, we followed this approach to study Fleming-Viot processes in \cite{MonmarcheJournel1,MonmarcheJournel2}}}.

\bigskip

The paper is organized as follows. The definitions and main results are presented in Section \ref{Section-def-non-lineaire} for the non-linear problem and in Section \ref{Section-Def-Particule} for the interacting particle system. These results are proven in Section \ref{Sec-Preuves}. More precisely, Section \ref{Sec-Preuves-Couple} is devoted to general considerations on couplings, while Section \ref{Sec-Preuves-NL} for the non-linear problem (and similarly \ref{Sec-Preuves-Particules} for the particle system) is decomposed in three parts: a coupling is constructed in the first one, which is then used in a second part to prove the main result in a simple case (uniform Doeblin condition), while the third part treats the more general case (local Doeblin condition with a Lyapunov function). Finally, Section \ref{Sec-Exemples} gathers various examples of applications of the main results.

\section{Definitions and main results}

\subsection{The non-linear process}\label{Section-def-non-lineaire}

Let $E$ be a Polish space,  $(\Omega,\mathcal F,\mathbb P)$ be a probability space  on which all the random variables in this work will be implicitly defined, and let $L$ be the infinitesimal generator of a conservative Feller   semi-group $(P_t)_{t\geqslant0}$  on $E$. Denote $\mathcal P(E)$ the set of probability measures on $E$ and $\mathcal D(E)$ the set of c\`adl\`ag functions for $\R_+$ to $E$, endowed with the Skorokhod topology.  From \cite[Chapter 17]{Kallenberg}, for all initial condition $x\in E$, there exists an associated  strong Markov process $(Z_t^x)_{t\geqslant 0}$ on $\mathcal D(E)$   such that $P_t f(x) = \mathbb E \po f(Z_t^x)\pf$ for all $t\geqslant 0$ and $f\in \mathcal C_0(E)$. In order to avoid regularity or well-posedness  considerations, we will rather work at the level of the stochastic process than of the generator.

For all $\nu \in \mathcal P(E)$, let $\lambda_\nu$ be a measurable function from $E$ to $\R_+$, and let $Q_\nu$ be a Markov kernel on $E$, namely a measurable function from $E$ to $\mathcal P(E)$. Throughout all this work, if $Q:E\rightarrow \mathcal P(E)$ is a Markov kernel, we also denote $Q$ the Markov operator defined by $Qf(x) = (Q(x))(f)$ for all $x\in E$ and all bounded measurable functions $f$. We call $\lambda_\nu$ the non-linear jump rate, and $Q_\nu$ the non-linear jump kernel.

Let $G_\nu:E\times[0,1] \rightarrow E$ be a representation of $Q_\nu$, namely be such that, if $U$ is a uniformly distributed random variable (r.v.) on $[0,1]$ then, for all $x\in E$, $G_\nu(x,U)$ is distributed according to $Q_\nu(x)$. From \cite[Corollary 7.16.1]{BertsekasShreve} such a representation always exists. Similarly, let $H:(x,u)\in E\times[0,1] \mapsto (H_t(x,u))_{t\geqslant0}\in  \mathcal D(E)$ be a representation of the kernel $x \mapsto  Law\po (Z_t^x)_{t\geqslant0}\pf$.

\begin{assu}\label{Hyp-taux_borne}
There exists $\lambda_*$ such that for all $\nu \in \mathcal P(E)$ and $x\in E$, $\lambda_{\nu}(x) \leqslant \lambda_*$.
\end{assu}
This assumption that the jump rate is bounded uniformly, both in $x\in E$ and $\nu\in\mathcal P(E)$, is made for simplicity. Although it already holds in many interesting applications, it may sometimes be too restrictive, but in many such cases it can be circumvented via proper a priori bounds (see e.g. Section \ref{Sec-TCP}), specific to the problem at hand.

 Let $\MRP$ be the set of  measurable functions from $\R_+$ to $\mathcal P(E)$. Under Assumption \ref{Hyp-taux_borne}, for  $x\in E$ and  $\mu: t\mapsto \mu_t$ in $\MRP$, we want to define an inhomogeneous Markov process $(X_t^{\mu,x})_{t\geqslant 0}$ starting from $x$ which, loosely speaking, follows the Markov dynamics of the semi-group $(P_t)_{t\geqslant0}$ between some random jump times  drawn at rate $t\mapsto \lambda_{\mu_t}(X_t^{\mu,x})$ at which it jumps to  new values drawn according to the distribution $Q_{\mu_t}(X_t^{\mu,x})$.   More precisely, consider  an i.i.d. sequence $(S_k,U_k,V_k,W_k)_{k\in\mathbb N}$ where, for all $k\in\mathbb N$, $S_k$, $U_k$, $V_k$ and $W_k$ are independent one from the other, $S_k$ follows an  exponential law with parameter 1  and $U_k$, $V_k$ and $W_k$ a uniform distribution over $[0,1]$. Set $T_0 = 0$, $X_0^{\mu,x} = x$ and suppose that $T_n\geqslant 0$ and $(X_t^{\mu,x})_{t\in [0,T_n]}$ have been defined for some $n\in \mathbb N$ and are independent from $(S_k,U_k,V_k,W_k)_{k\geqslant n}$. Set $T_{n+1} = T_n +   S_n / \lambda_*$.
\begin{itemize}
\item For all $t\in(T_n,T_{n+1})$, set $X_t^{\mu,x} = H_{t-T_n}(X_{T_n}^{\mu,x} ,U_n)$.
\item If $V_n \leqslant \lambda_{\mu_{T_{n+1}}}(H_{{S_n} / \lambda_*}(X_{T_n}^{\mu,x} ,U_n))/\lambda_*$ then set 
\[X_{T_{n+1}}^{\mu,x} = G_{\mu_{T_{n+1}}}(H_{S_n / \lambda_*}(X_{T_n}^{\mu,x} ,U_n),W_n)\,.\]
\item Else, set $X_{T_{n+1}}^{\mu,x} = H_{S_n / \lambda_*}(X_{T_n}^{\mu,x},U_n)$.
\end{itemize}
By induction, $(X_t^{\mu,x})_{t\in [0,T_n]}$ is thus defined for all $n\in\N$ and, since $T_n$ is the $n^{th}$ jump time of a Poisson process, it goes almost surely to infinity as $n$ goes to infinity, so that $X_t^{\mu,x}$ is almost surely defined for all $t\geqslant 0$. Denote by $(R_{s,t}^\mu)_{t\geqslant s\geqslant 0}$ the inhomogeneous Markov semi-group defined by
\begin{eqnarray*}
R_{s,t}^\mu f(x) & = & \mathbb E f\po X^{\sigma_s \mu,x}_{t-s}\pf
\end{eqnarray*}
for all bounded measurable functions $f$, $t\geqslant s\geqslant 0$ and $x\in E$, where $\sigma_s$ is the shift of time $s$ defined by $(\sigma_s \mu)_t = \mu_{s+t}$ for all $(\mu_t)_{t\geqslant0}$ and all $t\geqslant 0$. Note that, at least formally for suitable functions $f_0$, $f_t = R_{0,t}^\mu f_0$ is a solution of
\begin{eqnarray}\label{Eq-Pit}
\partial_t f_t &  = & L_t f_t\ := \   L f_t   +  \lambda_{\mu_t} \po Q_{\mu_t} f_t- f_t\pf\,.
\end{eqnarray} 
\begin{dfntn}\label{Def-solution}
We say that  $m\in\MRP$ is a solution of \eqref{Eq-EquaNonLin} with initial distribution $m_0$ if  $m_0 R_{0,t}^m  = m_t$ for all $t\geqslant 0$.
\end{dfntn}
 Note that Definition \ref{Def-solution} gives solutions in a quite weak sense, but in many cases additional informations on $L$ allows to   precise the meaning of \eqref{Eq-Pit} and thus obtain, by uniqueness of the solution, that a solution in the sense of Definition \ref{Def-solution} is in fact a strong solution of \eqref{Eq-EquaNonLin}.
 
 For $\nu_1,\nu_2 \in\mathcal P(E)$, we denote their total variation distance
\[\|\nu_1 - \nu_2\|_{TV} \ = \ \underset{\|f\|_\infty \leqslant 1}\sup |(\nu_1  - \nu_2)(f)|\,.\]
\begin{assu}\label{Hyp-taux_Lipschitz}
{The function $\nu \mapsto \lambda_\nu (Q_\nu - \mathrm{Id})$ is Lipschitz for the total variation norm, in the sense that there} exists $\theta  >0$ such that for all $\nu_1,\nu_2 \in \mathcal P(E)$,
\[\underset{\|f\|_\infty \leqslant 1}\sup \|\lambda_{\nu_1} (Q_{\nu_1} f - f)  - \lambda_{\nu_2} (Q_{\nu_2} f -f) \|_\infty \ \leqslant \ \theta \|\nu_1-\nu_2\|_{TV}\,.\]
\end{assu}

\begin{example}
{Let us consider a simple process in order to illustrate the assumptions in all this section. We consider a birth and death Markov chain on $E=\N$. The linear Markov part is given by
\begin{equation}\label{eq:Lexample}
Lf(x) =  b(x) \po f(x+1)-f(x)\pf + d(x) \po f(x-1)-f(x)\pf\,, 
\end{equation}
for some intrinsic birth and death rates $b,d\geqslant 0$, with $d(0)=0$. Then, we add non-linear rates given by
\[b_\nu(x) = \sum_{k=0}^\infty B(x,k) \nu(k)\,,\qquad d_\nu(x) = \sum_{k=0}^\infty D(x,k) \nu(k)\]
for some functions $B,D:\N^2\rightarrow \R_+$ for all probability measure $\nu$ on $\N$, with $D(0,k)=0$ for all $k\in\N$. Then we set
\[\lambda_\nu(x) = b_\nu(x) + d_\nu(x)\,,\qquad Q_\nu (x) = \frac{1}{\lambda_\nu(x)} \po b_\nu(x)\delta_{x+1} + d_\nu(x) \delta_{x-1}\pf\,,\]
with the convention $Q_\nu (x) = \delta_x$ if $\lambda_\nu(x)=0$. In this example, Assumption~\ref{Hyp-taux_borne} holds if $B$ and $D$ are bounded, and for a function $f$ with $\|f\|_\infty\leqslant 1$ and $\nu_1,\nu_2$ two probability measures on $\N$,
\begin{eqnarray*}
 \|\lambda_{\nu_1} (Q_{\nu_1} f - f)  - \lambda_{\nu_2} (Q_{\nu_2} f -f) \|_\infty & = &  \|(b_{\nu_1}-b_{\nu_2})(f(\cdot +1)-f) + (d_{\nu_1}-d_{\nu_2})(f(\cdot -1)-f) \|_\infty \\
 & \leqslant & 2\|b_{\nu_1}-b_{\nu_2}\|_\infty + 2 \|d_{\nu_1}-d_{\nu_2}\|_\infty \\
 & \leqslant  & 2 \po \|B\|_\infty + \|D\|_\infty\pf \|\nu_1-\nu_2\|_{TV}\,,
\end{eqnarray*}
which means Assumption~\ref{Hyp-taux_Lipschitz} also holds in this case.
}
\end{example}

\begin{prpstn}\label{Prop-ExistUniqueNonLin}
Under Assumptions \ref{Hyp-taux_borne} and \ref{Hyp-taux_Lipschitz}, for all $m_0\in\mathcal P(E)$, there exists a unique solution of \eqref{Eq-EquaNonLin}. Moreover, if $t\mapsto m_t,h_t$ are solutions with respective initial conditions $m_0$ and $h_0$, then for all $t\geqslant 0$
\[\| m_t - h_t\|_{TV} \ \leqslant \ e^{\theta t} \| m_0-h_0\|_{TV}\,.\]
\end{prpstn}
This will be proven in Section \ref{Section-NLTV} through a fixed point procedure.  
 As far as the long time behavior of the process is concerned, let us first state a simple result:
\begin{assu}\label{Hyp-DoeblinUnif}
There exist $t_0,\alpha\geqslant 0$ such that for all $x,y\in E$,
\begin{eqnarray*}
\| \delta_x P_{t_0} - \delta_y P_{t_0}\|_{TV} & \leqslant &  2(1-\alpha)\,.
\end{eqnarray*}
\end{assu}

\begin{example}
{
In the previous example, this is true for instance if there exists $x_0\in\N$ such that $b(x)+b_{\nu}(x)=0$ for all $x \geqslant x_0$ for all $\nu$,  in which case we can restrict the chain to the finite state $\cco 0,x_0\ccf$. If $b(x)>0$ for all $x<x_0$ and $d(x)>0$ for all $x>0 $ then the Markov chain on $\cco 0,x_0\ccf$ with generator \eqref{eq:Lexample} is irreducible on a finite state, hence satisfies Assumption~\ref{Hyp-DoeblinUnif}.
}
\end{example}

\begin{thrm}\label{Thm-TVNonLin}
Under Assumptions \ref{Hyp-taux_borne}, \ref{Hyp-taux_Lipschitz} and \ref{Hyp-DoeblinUnif}, if $t\mapsto m_t,h_t$ are solutions of \eqref{Eq-EquaNonLin} with respective initial distributions $m_0,h_0\in\mathcal P(E)$ then, for all $t\geqslant 0$,
\[\| m_t - h_t\|_{ TV} \ \leqslant \ e^{\theta t_0} \po e^{\theta t_0 } - \alpha e^{-\lambda_* t_0}\pf^{\lfloor t/t_0 \rfloor} \| m_{0} - h_{0}\|_{TV}\,.\]
\end{thrm}
An equivalent result is proved in \cite{CanizoYoldas}, for a particular model, via the Duhamel formula. We will rather establish it (in Section \ref{Section-NLTV}) with a similar but probabilistic point of view. This proof will introduce, in a simpler context, most of the ideas then used in the more general proof of Theorem \ref{Thm_non-lin} below.

The uniform Doeblin condition of Assumption \ref{Hyp-DoeblinUnif} is quite demanding, in particular if $E$ is not compact. The classical Foster-Lyapunov approach \cite{HairerMattingly2008,MeynTweedie} for Markov processes usually gathers a local Doeblin condition on compact sets and the existence of a Lyapunov function that tends to decrease on average along a trajecteroy, away from some compact. The following  counterpart of this strategy for non-linear perturbations seems to be new.

For a measurable function $\V$ from $E$ to $[1,\infty)$, denote $P_\V(E):=\{\nu\in \mathcal P(E),\ \nu(\V)<\infty\}$ and for $\mu,\nu\in\mathcal P_\V(E)$,  $\|\mu-\nu\|_\V\ = \ |\mu-\nu|(\V)$ the so-called $\V$-norm of $\mu-\nu$. In particular, if $\V=1$, then $\|\cdot\|_\V = \|\cdot\|_{TV}$. 
\begin{assu}\label{Hyp-MeanField-Exist}
There exist a measurable $\V:E\rightarrow [1,\infty)$, $\theta\geqslant 0$, $\rho, t_0>0$, $\rho_*\in\R$, $\alpha\in[0,1]$, $\eta\in[0,1)$ and $M ,\gamma_*\geqslant 1$ such that  the following holds.
\begin{itemize}
\item  {Lyapunov conditions.} For all  $\mu\in\MRP$, $x\in E$, $t\geqslant 0$, and almost all $u\in[0,1]$,
\begin{eqnarray}
R_{0,t}^\mu \V(x) & \leqslant & e^{- \rho t }\V(x) + \rho \int_0^t e^{ \rho(s-t) } \po \eta \mu_s(\V) + M\pf \dd s\label{Eq-Hyp-NonLin-Lyap}\\ 
P_t \V(x) & \leqslant & e^{\rho_* t} \V(x)   \label{Eq-Hyp-NonLin-Gene-1}\\
\V\po G_{\mu_t}(x,u)\pf & \leqslant & \gamma_* \V(x)\,.\label{Eq-Hyp-NonLin-Gene-4}
\end{eqnarray} 
\item {Lipschitz condition in $\V$-norm.} For all $\nu_1,\nu_2 \in \mathcal P_\V(E)$,
\begin{eqnarray}\label{Eq-Hyp-NonLin-Contract} 
\underset{\| f\|_\infty \leqslant 1 }\sup \|\lambda_{\nu_1} (Q_{\nu_1} f - f)  - \lambda_{\nu_2} (Q_{\nu_2} f -f) \|_\infty & \leqslant & \theta \|\nu_1-\nu_2\|_{\V}\,.
\end{eqnarray}
%
%
%
%
\item {Local Doeblin condition.} For all $x,y\in E$ with $\V(x)+\V(y) \leqslant 8(1+M/(1-\eta))$,
\begin{eqnarray}\label{Eq-NonLin-Gene-Doeblin}
\| \delta_{x} P_{t_0} - \delta_{y} P_{t_0}\|_{TV} & \leqslant & 2(1-\alpha)\,.
\end{eqnarray}
\end{itemize}
\end{assu}

Note that \eqref{Eq-Hyp-NonLin-Lyap} and \eqref{Eq-Hyp-NonLin-Gene-1} can usually be established if, for all $\nu\in\mathcal P_\V$,
\begin{eqnarray}\label{eq:LyapunovLV}
L\V + \lambda_\nu \po Q_\nu \V - \V\pf & \leqslant & -\rho \po  \V - \eta \nu(\V) - M\pf \,, \qquad L\V \leqslant \rho_* \V\,,
\end{eqnarray}
provided  additional regularity and/or truncation arguments, see e.g. \cite{MeynTweedie} or Section \ref{Sec-Exemples}. The Doeblin condition \eqref{Eq-NonLin-Gene-Doeblin} is always satisfied with $\alpha=0$.

\begin{example}
{Consider the previous running example. Assume that $B$ and $D$ are bounded, so that Assumptions~\ref{Hyp-taux_borne} and \ref{Hyp-taux_Lipschitz} (hence \eqref{Eq-Hyp-NonLin-Contract}) hold. Assume that $b(x)>0$ for all $x\in\N$ and $d(x)>0$ for all $x>0$, so that the Markov chain with generator~\eqref{eq:Lexample} is irreducible, hence satisfies \eqref{Eq-NonLin-Gene-Doeblin} on all finite subsets of $\N$. In order for Assumption~\ref{Hyp-MeanField-Exist} to holds, it only remains to establish the Lyapunov conditions~\eqref{Eq-Hyp-NonLin-Lyap}, \eqref{Eq-Hyp-NonLin-Gene-1} and \eqref{Eq-Hyp-NonLin-Gene-4} with a function $\V$ which goes to infinity at infinity. Let $\V(x) = e^{ax}$ for some $a>0$ to be chosen later on. Remark that, for all $\nu,x$, if $X\sim Q_\nu(x)$, then $|X-x|=1$ almost surely, hence $\V(X) \leqslant e^a\V(x)$, so that  \eqref{Eq-Hyp-NonLin-Gene-1} holds with $\gamma_* = e^a$. To establish \eqref{Eq-Hyp-NonLin-Lyap} and \eqref{Eq-Hyp-NonLin-Gene-1}, let us prove \eqref{eq:LyapunovLV}, which is sufficient by standard arguments (again, we refer to Section \ref{Sec-Exemples} for details). For  $\nu\in\mathcal P_\V$,
\[
\frac1{\V(x)}\po L\V(x) + \lambda_\nu(x) \po Q_\nu \V(x) - \V(x)\pf\pf   =  \po b(x)+b_{\nu}(x)\pf \po e^{a} -1\pf + \po d(x)+d_{\nu}(x)\pf \po e^{-a}-1\pf\,.
\]
Assume that 
\[c:=\liminf_{x\rightarrow \infty} \po d(x) - b(x) - \sup_{k\in\N} B(x,k) \pf >0\,, \]
i.e. for  large populations the death rate is predominant, uniformly in the non-linear birth rate. 
Let $a,\rho>0$ be small enough and $x_0\in\N$ be large enough so that $(1-e^{-a})d(x) -(e^a-1)(b(x) +\sup_{k\in\N} B(x,k))\geqslant\rho$ for all $x\geqslant x_0$. Then for all $\nu\in\mathcal P_\V$ and $x\geqslant x_0$,
\[
 L\V(x) + \lambda_\nu(x) \po Q_\nu \V(x) - \V(x)\pf   \leqslant -\rho \V(x)\,.
\]
from which the first part of \eqref{eq:LyapunovLV} holds with $\eta=0$ and
\[M = \sup_{x\in\N} \V(x) + \frac1\rho \po L\V(x) + \lambda_\nu(x) \po Q_\nu \V(x) - \V(x)\pf \pf \,.\]
The second part of \eqref{eq:LyapunovLV}  corresponds to the same computation in the case $B=D=0$, hence it holds with $\rho_*=-\rho$. As a conclusion, in this example, Assumption~\ref{Hyp-MeanField-Exist} holds as soon as the rates $B$ and $D$ are bounded and $c>0$. Besides, these conditions do not imply the stronger Assumption~\ref{Hyp-DoeblinUnif}, which would require that the process with generator $L$ comes from infinity in a finite time, which is for instance not the case if $b$ and $d$ are bounded (in which case it is easily checked that, at any fixed time $t_0$, the probability that a process starting at $x$ has reached a position lower than $x/2$ goes to $0$ as $x\rightarrow \infty$, which prevents  Assumption~\ref{Hyp-DoeblinUnif} to hold with $\alpha>0$ uniformly in $x,y$).
 }
\end{example}

\begin{thrm}\label{Thm-NonLin-Exist}
Under Assumptions \ref{Hyp-taux_borne} and \ref{Hyp-MeanField-Exist}, for all $m_0\in\mathcal P_\V$, there exists a unique solution $t\mapsto m_t$ of \eqref{Eq-EquaNonLin} with initial condition $m_0$. Moreover, for all $t\geqslant 0$,
\begin{eqnarray}\label{Eq-NL-Lyap}
m_t(\V) & \leqslant & e^{-\rho(1-\eta) t} m_0(\V) + (1-e^{-\rho (1-\eta)t}) \frac{M}{1-\eta}\,.
\end{eqnarray}
\end{thrm}

\begin{thrm}\label{Thm_non-lin}
Under Assumptions \ref{Hyp-taux_borne} and \ref{Hyp-MeanField-Exist}, if $\alpha>0$, denote
\begin{eqnarray*}
\beta & =&    \frac{2 \po 1-e^{-\rho t_0} \pf e^{\lambda_* t_0} }{\alpha} \po \frac{M}{1-\eta} +1\pf   \\
C_* & =& \theta \frac{  2(1+\beta) \gamma_*  \po \frac{M}{1-\eta}+1\pf}{\po\rho + \theta  \gamma_* \pf \po\frac{M}{1-\eta} +1\pf}  e^{\po \rho_* + \lambda_*(\gamma_* -1)\pf t_0}  \po e^{\po\rho + \theta  \gamma_* \pf \po\frac{M}{1-\eta} +1\pf t_0} - 1\pf\\
\kappa & = & \frac{1+e^{-\rho t_0}}{2} \vee \po  1-\frac12 \alpha e^{-\lambda_* t_0}\pf\\
\tilde\kappa & = & e^{-\rho(1-\eta) t_0} \vee \po \kappa  + C_*\pf\,.
\end{eqnarray*}
Suppose that $\tilde  \kappa < 1$. 
 Then Equation \eqref{Eq-EquaNonLin} admits an equilibrium $\mu_\infty\in \mathcal P_\V(E)$, i.e. a  constant solution $t\mapsto m_t = \mu_\infty$ for all $t\geqslant 0$, which is unique in $\mathcal P_\V(E)$. Moreover, $\mu_\infty(\V) \leqslant M/(1-\eta )$, and  for all $m_0\in \mathcal P_\V(E)$ and all $t\geqslant 0$,
\begin{eqnarray} \label{Eq-Contract_ThmNonLin}
 \|m_t-\mu_\infty\|_\V &\leqslant & 2 \tilde  \kappa^{ t/t_0 -1}\po \beta + 1 + \frac{M}{1-\eta}\pf   m_0(\V)\,.
\end{eqnarray}
\end{thrm}
Remark that, for fixed values of the other parameters of Assumptions \ref{Hyp-taux_borne} and \ref{Hyp-MeanField-Exist},  $\tilde \kappa<1$ if $\theta$ is small enough.

\begin{example}
{In the previous example, let $\theta = 2(\|B\|_\infty + \|D\|_\infty)$. It is readily checked that all the estimates obtained in this example to establish Assumption~\ref{Hyp-MeanField-Exist} are uniform for small $\theta$. As a conclusion, for fixed rates $b,d$ satisfying $b(x)>0$ for all $x\in\N$, $d(x)>0$ for all $x>0$ and $\limsup(b-d)<0$, there is an explicit $\theta_0\geqslant 0$ so that  Theorem~\ref{Thm_non-lin} holds with $\tilde\kappa<1$ as soon as $B$ and $D$ are such that $\theta \leqslant \theta_0$. }
\end{example}



\subsection{Weakly interacting Markov particles}\label{Section-Def-Particule}

Let $N\in\N$ and, for $i\in\cco 1,N\ccf$, let $E_i$ be a Polish space and $(P_t^i)_{t\geqslant0}$  be a conservative Feller semi-group  on $E_i$. Let $E = \prod_{i=1}^N E_i$ and $(P_t)_{t\geqslant 0}$ be the semi-group on $E$ defined by
\[P_t \po \prod_{i=1}^N f_i\pf(x)\ = \ \prod_{i=1}^N P_t ^if_i(x_i)\]
for all $x\in E$, $t\geqslant 0$ and all measurable bounded functions $f_i$ on $E_i$, $i\in\cco 1,N\ccf$. For $i\in\cco 1,N\ccf$ let $H^i$ be a representation of $z\in E_i \mapsto Law((Z^{z,i}_t)_{t\geqslant 0})$ on $\mathcal D(E_i)$ where, for all $z\in E_i$, $(Z^{z,i}_t)_{t\geqslant 0}$ is a Markov process associated to  $(P_t^i)_{t\geqslant0}$ with initial value $z$. In particular, if $(U_1,\dots,U_N)$ are i.i.d. r.v. uniformly distributed over $[0,1]$ and $x\in E$ then $\po H_t^i(x_i,U_i)\pf_{i\in\cco 1,N\ccf,t\geqslant 0}$ is a Markov process on $E$ associated to  $(P_t)_{t\geqslant0}$ with initial value $x$. For all $i\in\cco 1, N\ccf$ let $\lambda_i$ be a jump rate on $E$, i.e. be a measurable function  from $E$ to $\R_+$, and $Q_i$ be a measurable function from $E$ to $\mathcal P(E_i)$ with representation $G_i$. We extend $Q_i$ as a Markov kernel $Q_i'$ on $E$ defined for all $x\in E$ by 
\[Q_i'(x)\ = \ \delta_{x_1}\otimes \dots \otimes \delta_{x_{i-1}} \otimes Q_i(x) \otimes \delta_{x_{i+1}} \otimes \dots \otimes \delta_{x_N}\,.\]

\begin{assu}\label{Hyp-TauxBorne-Particules}
There exists $\lambda_*>0$  such that, for all $x\in E$ and all $i\in\cco 1,N\ccf$, $\lambda_i(x) \leqslant \lambda_*$.
\end{assu}

Under Assumption \ref{Hyp-TauxBorne-Particules}, let $Q$ be the Markov kernel on $E$ defined by
\[Q(x) \ = \ \frac1N \sum_{i=1}^N \frac{\lambda_i(x)}{\lambda_*} Q_i'(x) + \po 1 - \frac{\lambda_i(x)}{\lambda_*}\pf \delta_x\]
for all $x\in E$, and let $G$ be a representation of $Q$. To jump from a position $x$ to a r.v. drawn according to $Q(x)$ is equivalent to chose a coordinate $i$ uniformly in $\cco 1,N\ccf$ and, with probability $\lambda_i(x)/\lambda_*$, make it jump to a new position drawn according to $Q_i(x)$, else leave it at its current position, and in either case leave all the other coordinates unchanged.

Let $x\in E$. We define a Markov process $(X_t^x)_{t\geqslant 0} = (X_{i,t}^x)_{i\in\cco 1,N\ccf, t\geqslant 0}$ on $E$ as follows. Consider  an i.i.d. sequence $(S_k,(U_{i,k})_{i\in\cco1,N\ccf},V_k)_{k\in\mathbb N}$ where, for all $k\in\N$, $S_k,V_k$ and $(U_{i,k})_{i\in\cco1,N\ccf}$ are independent one from the other, $S_k$ follows an exponential law with parameter 1, $V_k$   a uniform distribution over $[0,1]$   and $(U_{i,k})_{i\in\cco1,N\ccf}$ are i.i.d. r.v. uniformly distributed over $[0,1]$. Set $T_0 = 0$,  $X_{0}^x = x$ and, for all $n\in\N$, $T_{n+1} = T_n + (N\lambda_*)^{-1} S_n$. Suppose that  $(X_{t}^x)_{ t\in [0,T_n]}$ have been defined for some $n\in \mathbb N$ and is independent from $(S_k,(U_{i,k})_{i\in\cco1,N\ccf},V_k)_{k\geqslant n}$. For all $t\in(T_n,T_{n+1})$ and all $i\in\cco1,N\ccf$, set $X_{i,t}^x = H_{t-T_n}^i(X_{i,T_n}^x ,U_{i,n})$. Finally, set 
\[X_{T_{n+1}}^x \ = \ G\po \po H_{T_{n+1}-T_n}^i(X_{i,T_n}^x ,U_{i,n})\pf_{i\in\cco 1,N\ccf},V_n\pf\,.\]
By induction, $(X_{t}^x)_{t\in [0,T_n]}$ is thus defined for all $n\in\N$ and, since $T_n$ is the $n^{th}$ jump time of a Poisson process, it goes almost surely to $+\infty$ as $n$ goes to infinity, so that $X_t$ is almost surely defined for all $t\geqslant 0$. Let $(R_t)_{t\geqslant 0}$ be the associated Markov semi-group, i.e. 
\[R_t f(x) \ = \ \mathbb E\po f\po X_t^x\pf\pf\]
for all bounded measurable functions $f$ on $E$.

\begin{example}
{Let us briefly introduce the particle system which is the analogous of the running example of the previous section, when the non-linearity is replaced by a mean-field interaction. The Markov part is independent from the particle, i.e. $P_t^i = P_t$ for all $i$, and it is the semi-group associated to the generator given by \eqref{eq:Lexample}. Given $B,D$ as in the previous section, set
\[b_i(x) = \frac1N\sum_{j=1}^N  B(x_i,x_j)\qquad d_i(x) = \frac1N\sum_{j=1}^N D (x_i,x_j)\,.\]
Then the interacting jump mechanisms are given by
\[\lambda_i(x) = b_i(x)+d_i(x)\,,\qquad Q_i(x) = \frac{1}{\lambda_i(x)} \po b_i(x) \delta_{x_i+1} + d_i(x)\delta_{x_i-1} \pf\,.  \]
Then, all the assumptions stated in the rest of this section can be proven under suitable conditions on $b,d,B,D$ similarly to the non-linear case of the previous section. The computations being exactly the same, we will not repeat this discussion.
  }
\end{example}



Consider on $E$ the   $l_0$  distance defined by
\[\bar d_1(x,y)\ = \  2 \sharp \{ i\in\cco 1,N\ccf,\ : \ x_i\neq y_i\}\,, \]
for all  $x,y\in E$, where $\sharp A$ denote the cardinality of a finite set $A$. 
 Remark that the topology induced by $\dist$ is equivalent to the trivial one, with more precisely
\[2\1_{x\neq y}\ \leqslant\ \dist(x,y) \ \leqslant\ 2N \1_{x\neq y}\,.\]


\begin{assu}\label{Hyp-TauxLipschitz-Particules}
There exists $\theta,t_0,\alpha\geqslant0$ such that, for all $i\in\cco 1,N\ccf$, the following holds.
\begin{itemize}
\item {Lipschitz condition.} For all $x,y\in E$ with $y_i=x_i$,
\begin{eqnarray}\label{Eq-TauxLipschitz-Particules}
  \underset{\|f\|_\infty \leqslant 1}\sup \|\lambda_i(x) \po Q_i f(x)-f(x_i)\pf  - \lambda_i(y) \po Q_i  f(y)-f(y_i)\pf \|_\infty & \leqslant &  \theta \frac{\dist(x,y)}N\,.
\end{eqnarray}
\item {(Uniform) Doeblin condition.} For all $x,y\in E_i$,
\begin{eqnarray}\label{Eq_Prop-Particule_HypDoeblin_Unif} 
\| \delta_x P_{t_0}^i - \delta_y P_{t_0}^i\|_{TV} & \leqslant &  2(1-\alpha)\,.
\end{eqnarray}
\end{itemize}
\end{assu}

Here is the analogous of Theorem \ref{Thm-TVNonLin} for interacting particles.

\begin{thrm}\label{Thm_particules_TV}
Under Assumptions \ref{Hyp-TauxBorne-Particules} and \ref{Hyp-TauxLipschitz-Particules}, for all $m_0,h_0\in\mathcal P(E)$ and $t\geqslant 0$,
\begin{eqnarray*}
 \|m_0 R_t - h_0 R_t\|_{TV} & \leqslant &  N e^{ \theta t_0  } \po e^{  \theta t_0 } - \alpha e^{-\lambda_* t_0}\pf^{\lfloor t/t_0\rfloor} \|m_0   - h_0  \|_{TV}\,.
\end{eqnarray*}
Moreover, denoting $m'_t$ and $h_t'$ the respective laws of $X_{I,t}$ and $Y_{I,t}$ where $X_t \sim m_0 R_t$, $Y_t\sim h_0 R_t$ and $I$ is uniformly distributed over $\cco 1,N\ccf$ and independent from $X_t$ and $Y_t$, then
\begin{eqnarray*}
 \|m_t' - h_t'\|_{TV} & \leqslant &  e^{  \theta t_0  } \po e^{ \theta t_0 } - \alpha e^{-\lambda_* t_0}\pf^{\lfloor t/t_0\rfloor} \|m_0' - h_0'\|_{TV}\,.
\end{eqnarray*}
\end{thrm}

Note that, using a naive approach (trying to merge at time $t_0$ two whole systems of particles), it is easy to prove that the condition \eqref{Eq_Prop-Particule_HypDoeblin_Unif} together with Assumption \ref{Hyp-TauxBorne-Particules}  imply the uniform Doeblin condition
\begin{eqnarray*}
\frac12 \| \delta_x R_{t_0} - \delta_y R_{t_0}\|_{TV} & \leqslant &   1-\po e^{-\lambda_* t_0} \alpha\pf^N  \,,
\end{eqnarray*}
but then the contraction rate is geometrically poor with respect to $N$. This is why the Foster-Lyapunov approach is sometimes thought to scale badly with dimension.

\bigskip

Let us show how Theorems \ref{Thm_particules_TV} may apply in the case of a mean-field particle system associated to a non-linear process such as defined in Section \ref{Section-def-non-lineaire}. Suppose that $E_i = E_1$, $P_t^i = P_t^1$ and
\[\lambda_i(x)\ = \ \ \lambda_{\mu_N(x)}(x_i)\,,\qquad Q_i(x)\ = \ Q_{\mu_N(x)}(x_i)\]
 for  all $i\in\cco 1,N\ccf$ and $x\in E = E_1^N$, where $\lambda_\nu$ and $Q_\nu$ are given in Section \ref{Section-def-non-lineaire} and
 \[\mu_N(x) \ = \ \frac1N \sum_{i=1}^N \delta_{x_i}\]
is the empirical distribution of $x$. Remark that, for $x,y\in E$,
\[\|\mu_N(x) - \mu_N(y) \|_{TV} \ \leqslant \ \dist(x,y)/N\,.\]
Indeed, $(x_I,y_I)$ with a r.v. $I$ uniformly distributed over $\cco 1,N\ccf$ is a coupling of $\mu_N(x)$ and $ \mu_N(y)$ with $\mathbb P (x_I=y_I) =\dist(x,y)/N$. Hence, Assumptions \ref{Hyp-taux_borne}, \ref{Hyp-taux_Lipschitz} and \ref{Hyp-DoeblinUnif} imply Assumptions \ref{Hyp-TauxBorne-Particules} and \ref{Hyp-TauxLipschitz-Particules}, with the same $\lambda_*,\theta,t_0,\alpha$.

If the initial positions $(X_{i,0})_{i\in\cco 1,N\ccf }$ are i.i.d. r.v. with a given law $m_0\in \mathcal P(E_1)$, then the particles are 
  exchangeables, in the sense that for all $t\geqslant 0$, $(X_{\sigma(i),t})_{i\in\cco 1,N\ccf }$ and $(X_{i,t})_{i\in\cco 1,N\ccf }$ have the same distribution for all permutation $\sigma$ of $\cco 1,N\ccf$. In particular, if $I$ is uniformly distributed over $\cco 1,N\ccf$ and independent from $(X_{t})_{t\geqslant 0}$ then $X_{I,t}$ and $X_{1,t}$ have the same law. This means that, under Assumptions \ref{Hyp-taux_borne}, \ref{Hyp-taux_Lipschitz}  and \ref{Hyp-DoeblinUnif}, then 
\begin{eqnarray*}
 \|m_t - h_t\|_{TV} & \leqslant &  e^{\lambda_* \theta t_0  } \po e^{\lambda_* \theta t_0 } - \alpha e^{-\lambda_* t_0}\pf^{\lfloor t/t_0\rfloor} \|m_0 - h_0\|_{TV}
\end{eqnarray*}
holds  with either $m_t$ and $h_t$ solutions of the non-linear equation \eqref{Eq-EquaNonLin} or the law of the first particle of the associated system.

\bigskip

We go back to the general case. In order to relax the strong uniform Doeblin condition \eqref{Eq_Prop-Particule_HypDoeblin_Unif}  under the assumption that the process admits a suitable Lyapunov function, we now state an analogous of Theorem \ref{Thm_non-lin} for the particle system. 

\begin{assu}\label{Hyp-Particules-Gene}
There exist $\theta\geqslant 0$, $\rho, t_0>0$, $\rho_*\in\R$, $\alpha\in(0,1]$, $\eta\in[1/2,1)$, $M ,\gamma_*\geqslant 1$ and, for all $i\in\cco 1,N\ccf$,  a measurable $\V_i:E_i\rightarrow [1,\infty)$,  such that, denoting $\V(x) = \sum_{i=1}^N \V_i(x_i)$,  the following holds for all $i\in\cco 1,N\ccf$.
\begin{itemize}
\item {Lyapunov conditions.} For all $x\in E$, $t\geqslant 0$, and almost all $u\in[0,1]$,
\begin{eqnarray}
R_{t} \V_i(x) & \leqslant & e^{- \rho t }\V_i(x_i) + \rho \int_0^t e^{ \rho(s-t) } \po \frac{\eta}N R_s\V(x) + M\pf \dd s\label{Eq-Hyp-Particule-Lyap}\\ 
P_t^i \V_i(x_i) & \leqslant & e^{\rho_* t} \V_i(x_i)   \label{Eq-Hyp-Particule-Gene-1}\\
\V_i\po G_{i}(x,u)\pf & \leqslant & \gamma_* \V_i(x_i)\,.\label{Eq-Hyp-Particule-Gene-4}
\end{eqnarray}
\item {Lipschitz condition.} For all $x,y\in E$ such that $x_i=y_i$,
\begin{eqnarray}\label{Eq-TauxLipschitz-Particules-V}
\underset{\| f\|_\infty \leqslant 1 }\sup |\lambda_i(x) (Q_i f(x) - f(x_i))  - \lambda_i(y)(Q_i  f(y) -f(y_i)) \|_\infty & \leqslant & \theta \frac{\bar d_1(x,y)}N\,.
\end{eqnarray}
\item {Local Doeblin condition.} For all $x_i,y_i\in E_i$ with $\V_i(x_i)+\V_i(y_i) \leqslant 16(1+\eta )M/(1-\eta)^2$,
\begin{eqnarray}\label{Eq-Particule-Gene-Doeblin}
\| \delta_{x_i} P_{t_0}^i - \delta_{y_i} P_{t_0}^i\|_{TV} & \leqslant & 2(1-\alpha)\,.
\end{eqnarray}
\end{itemize}
\end{assu}

Under Assumption \ref{Hyp-Particules-Gene}, denote $\bar d_\V$ and $\bar \rho$ the distances, respectively  on $E$ and $\mathcal P_\V(E)$, defined for all $x,y\in E$ by
\[\bar d_\V(x,y) =  \sum_{i=1}^N \mathbbm{1}_{x_i\neq y_i} \po \V_i(x_i)+  \V_i(y_i) + \frac{\V(x)+\V(y)}{N}\pf\,,\]
and for all $\mu,\nu\in\mathcal P_\V(E)$ by
\[\bar \rho (\mu,\nu) = \inf\left \{ \mathbb E \po \bar d_\V(X,Y)\pf,\ X\sim \mu,\ Y\sim \nu\right\}\,.\]
(See Section \ref{Sec-Preuves-Couple} for a more detailed definition of couplings.) Remark that $\|\mu-\nu\|_{TV} \leqslant  \bar \rho(\mu,\nu)$ and that $\|\mu-\nu\|_{\V} \leqslant  N\bar \rho(\mu,\nu)$. Note also that, in the case of a particle system associated to a mean-field equation, i.e. if the particles are exchangeables and $\V_i = \V_j$ for all $i,j\in\cco 1,N\ccf$, denoting $\mu'$ and $\nu'$ the laws of $X_1$ and $Y_1$  when $X\sim\mu$ and $Y\sim \nu$ (hence of $X_I$ and $Y_I$ with an independent $I$ uniform over $\cco 1,N\ccf$) then  $\|\mu'-\nu'\|_{\V_1} \leqslant\bar  \rho(\mu,\nu)/N$.

\begin{thrm}\label{Thm_particules_gene}
Under Assumption \ref{Hyp-TauxBorne-Particules} and \ref{Hyp-Particules-Gene}, denote
\begin{eqnarray*}
\beta & =&    \frac{4(1+\eta) M e^{\lambda_* t_0} }{\alpha (1-\eta)^2} \po 1 - e^{- \frac12\rho(1-\eta) t_0 } \pf \\
 \kappa & = & \frac{1+e^{- \frac12\rho(1-\eta) t_0 }}2  \vee \po  1-\frac12 \alpha e^{-\lambda_* t_0}\pf\\
\tilde \kappa & =& \kappa +  \kappa + \po e^{\theta t_0} - 1\pf (\gamma_*+1) e^{ \rho_* +\lambda_*(\gamma_*-1)  t_0} (1+\eta)\,.
\end{eqnarray*}
Suppose that $\tilde \kappa <1$. Then, $(R_t)_{t\geqslant 0}$ admits a unique invariant measure $\mu_\infty$ and for all $m_0\in\mathcal P_\V(E)$ and all $t\geqslant 0$,
\begin{eqnarray*}
\bar  \rho(m_0 R_t,\mu_\infty) & \leqslant &  \tilde \kappa^{\lfloor t/t_0\rfloor} \po 2N \po \beta  + M \frac{1+\eta}{(1-\eta)^2 }\pf  + \frac{1+\eta}{1-\eta }  \nu (\V) \pf\,.
\end{eqnarray*}
\end{thrm}

Note that, in the case of a system of interacting particles associated to a mean-field equation, the assumptions of Theorem \ref{Thm_particules_gene} are slightly stronger than the assumptions of Theorem \ref{Thm_non-lin}. Indeed, the right-hand side of \eqref{Eq-TauxLipschitz-Particules-V} involves $\bar d_1$ instead of  $\bar d_\V$, 
which would be the analogous of \eqref{Eq-Hyp-NonLin-Contract}. In other words, in this mean-field case, the assumptions of Theorem  \ref{Thm_particules_gene} are equivalent to the assumptions of Theorem \ref{Thm_non-lin} with the addition of Assumption \ref{Hyp-taux_Lipschitz}. The interesting question of weakening the condition \eqref{Eq-TauxLipschitz-Particules-V} in Theorem \ref{Thm_particules_gene} raises non-trivial considerations about the non-independence of $\V_i(X_{i,t})$ for different $i\in\cco 1,N\ccf$, and is postponed to a future work.

\section{Proofs}\label{Sec-Preuves}

\subsection{General considerations on couplings}\label{Sec-Preuves-Couple}



Let us first give an alternative representation of the $\V$-norm, for a given measurable function $\V$ from $E$ to $[1,\infty[$. 
Consider the distance $d_\V$ on $E$ defined for all $x,y\in E$ by
\[d_\V(x,y) =  \mathbbm{1}_{x\neq y} \po \V(x)+  \V(y)\pf\,.\]
Remark that it induces the trivial topology on $E$, since $d_\V({x,y}) \geqslant 2 \1_{x\neq y} $ for all $x,y\in E$. 
Hairer and Mattingly proved in \cite{HairerMattingly2008} that
\begin{eqnarray}\label{Eq-HairerMattingly}
\|\mu_1-\mu_2\|_\V & = & \sup_{\| \varphi\|_{\V,1}\leqslant1}\po \mu_1(\varphi) - \mu_2(\varphi)\pf \ = \ \sup_{\| \varphi\|_{\V,2}\leqslant1}\po \mu_1(\varphi) - \mu_2(\varphi)\pf 
\end{eqnarray}
where, for all measurable $\varphi$ from $E$ to $\R$,
\[\| \varphi\|_{\V,1} = \sup_{x\in E} \frac{\varphi{(x)}}{\V(x)}\,,\qquad \|\varphi\|_{\V,2} = \sup_{x\neq y} \frac{|\varphi(x) - \varphi(y)|}{d_\V(x,y)}\,.\]
Let us now give a coupling interpretation of $\|\cdot\|_\V$. For $\mu_1,\mu_2\in\mathcal{P}(E)$, we denote $\xi(\mu_1,\mu_2)$ the set of transference plane between $\mu_1$ and $\mu_2$, i.e. the set of probability measures $\nu$ on $E^2$ such that $(X,Y)\sim\nu$ implies that $X\sim \mu_1$ and $Y\sim\mu_2$. Such a random variable $(X,Y)$ is called a coupling of  $\mu_1$ and $\mu_2$. Alternatively, in that case, if $Z_i \sim \mu_i$ for $i=1,2$, then we also say that $(X,Y)$ is a coupling of $Z_1$ and $Z_2$.

\begin{lmm}\label{Lem-couplage}
For all $\mu_1$ and $\mu_2$, 
\[\|\mu_1-\mu_2\|_\V\ =\ \inf\left\{\mathbb E\po d_\V(X,Y)\pf\,:\, (X,Y)\sim\nu,\,\nu\in \xi(\mu_1,\mu_2)\right \}\,,\]
and the infimum is attained for some $\nu_{opt}\in \xi(\mu_1,\mu_2)$ that does not depend on $\V$.
\end{lmm}
\begin{proof}
Note that, according to \eqref{Eq-HairerMattingly}, $\|\cdot\|_\V$ is the Wasserstein-1 distance on $\mathcal P(E)$ associated to the distance $d_\V$ on $E$, so that the first statement of the lemma stems from the general result of duality for Wasserstein distances. Nevertheless, in the general case, the coupling would depend on the metric. Let us then give an elementary proof in the present specific case, straightforwardly adapted from the classical case of the total variation distance. First, note that $|\mu_1-\mu_2|=\mu_1+\mu_2 - 2 \mu_1\wedge\mu_2$. If $\nu \in \xi(\mu_1,\mu_2)$ and $(X,Y)\sim\nu$, then for all event $A$ of $E$,
\[\mathbb P(X\in A,\, X=Y) \ \leqslant\ \mu_1(A)\wedge\mu_2(A)\,.\]
In other word, the measure $\tilde \nu$ on $E$ defined by $\tilde \nu(\dd x) = \int_E \1_{x=y} \nu(\dd x,\dd y)$ satisfies $\tilde \nu \leqslant \mu_1\wedge\mu_2$.
Therefore, writing $\1_{x\neq y}=1-\1_{x=y}$,
\begin{eqnarray*}
 \mathbb E\po d_\V(X,Y)\pf &= &\mu_1(  \V) + \mu_2( \V) - \int_{E^2} \po   \V(x)+ \V(y)\pf\1_{x=y} \nu(\dd x,\dd y)\\
&=& \mu_1( \V) + \mu_2(  \V) - 2\int_E    \V(x)  \tilde \nu(\dd x)\\
&\geqslant & \mu_1( \V) + \mu_2(  \V) - 2\mu_1\wedge \mu_2(  \V)  \ =\ \|\mu_1-\mu_2\|_\V\,.
\end{eqnarray*}
It now remains to construct an optimal coupling, i.e. a measure $\nu\in\xi(\mu_1,\mu_2)$ such that equality holds. Let $p=\mu_1\wedge \mu_2(E)$. If $p=1$ then $\mu_1=\mu_2$ and the optimal coupling is to draw $X$ according to $\mu_1$ and then set $Y=X$. If $p=0$ then $\mu_1$ and $\mu_2$ are singular one to the other, and all couplings are optimal since
\[\mathbb E\po d_\V(X,Y)\pf \ \leqslant \ \mathbb E\po  \V(X) + \V(Y)\pf \ =\ \|\mu_1-\mu_2\|_\V\,.\]
If $p\in(0,1)$, we consider the probability measures
\[\nu_0 = \frac{\mu_1\wedge\mu_2}p \,,\qquad \nu_1 = \frac{\mu_1-\mu_1\wedge\mu_2}{1-p}\,,\qquad \nu_2 = \frac{\mu_2-\mu_1\wedge\mu_2}{1-p}\,.\]
Remark that $\nu_1$ and $\nu_2$ are singular one to the other, and that $\mu_i = p\nu_0+(1-p)\nu_i$ for $i=1,2$. Let $E_1,E_2,E_3$ and $B$ be independent random variables with $E_i\sim\nu_i$ for $i=1,2,3$ and $B\sim\mathcal B(p)$ the Bernoulli law with parameter $p$. Set $X=Y=E_0$ if $B=1$ and $X=E_1$, $Y=E_2$ if $B=0$. Then $X\sim \mu_1$ and $Y\sim\mu_2$, and $\{X=Y\} = \{B=1\}$. Conditionning on the values of  $B$,
\begin{eqnarray*}
\mathbb E\po d_\V(X,Y)\pf & = & (1-p) \mathbb E\po d_\V(X,Y)\ |\ B= 0\pf \\
& = & (1-p) \mathbb E\po  \V(E_1) +  \V(E_2)\pf \\
&= & \int  \V(x) \po \mu_1 + \mu_2 - 2 \mu_1\wedge \mu_2\pf(\dd x)\ = \ \|\mu_1-\mu_2\|_\V\,,
\end{eqnarray*}
which concludes. 
\end{proof}
This representation of $\|\cdot \|_\V$ is convenient since, as soon are we are able to construct a coupling of two probability measures, we get an upper bound on their distance.

\begin{lmm}\label{Lem-CouplageMarkov}
Let $(P_t)_{t\geqslant 0}$ be a conservative Feller semi-group on $E$, $t_0\geqslant 0$, $x,y\in E$ and $\V$ be a measurable function from $E$ to $[1,\infty)$. Then there exists a random variable $(Z^x_t,Z^y_t)_{t\geqslant 0}$ on $\mathcal D(E^2)$ such that, for $z=x,y$,  $(Z^z_t)_{t\geqslant 0}$ is a Markov process associated to $(P_t)_{t\geqslant 0}$ with $Z^z_0=z$, 
\begin{eqnarray}\label{Eq-LemCoupl1}
\mathbb E \po d_\V\po  Z_{t_0}^x,  Z^y_{t_0} \pf \pf & = & \|\delta_x P_{t_0} - \delta_y P_{t_0}\|_\V\,,
\end{eqnarray}
and moreover,
\begin{eqnarray}\label{Eq-LemCoupl2}
\forall s\geqslant t\geqslant 0\,,\qquad Z_t^x  =  Z_t^y & \Rightarrow &  Z_s^x = Z_s^x\,.
\end{eqnarray}
 We call such a process an optimal coupling of $(\delta_x P_t)_{t\geqslant 0}$ and of $(\delta_y P_t)_{t\geqslant 0}$ at time $t_0$.
\end{lmm}

\begin{proof}
Let $(X^x,X^y)$ be an optimal coupling of $\delta_x P_{t_0} $ and $\delta_y P_{t_0} $ such as constructed in Lemma \ref{Lem-couplage}. 
For $z\in E$, let $(V_t^z)_{t\geqslant 0}$ be a Markov process associated with $(P_t)_{t\geqslant0} $ with $V_0^z=z$, and let $W_t=V^z_{t-t_0}$ for $t\in[0,t_0]$. Then by usual conditionning arguments, $(W_t)_{t\in[0,t_0]}$ is an inhomogeneous Markov process with transitions
\[\mathbb P \po W_t \in A \ | \ W_s\in B \pf \ = \ \ \mathbb P \po V_{t_0-s} \in B \ | \ V_{t_0-t}\in A \pf \frac{\mathbb P \po V_{t_0-t}\in A \pf }{\mathbb P \po V_{t_0-s}\in B \pf }\]  
for all $0\leqslant s\leqslant t \leqslant t_0$ and all events $A,B$ of $E$. Let $(R_{s,t}^z)_{0\leqslant s\leqslant t\leqslant t_0}$ be the inhomogeneous  semi-group associated to $(W_t)_{t\in[0,t_0]}$. For $z=x,y$, let $(U^z_t)_{t\in[0,t_0]}$ be a Markov process associated to  $(R_{s,t}^z)_{0\leqslant s\leqslant t\leqslant t_0}$  with $U^z_0 = X^z$. For all $t\in [0,t_0]$ and $z=x,y$, set $Z_t^z = U^z_{t_0-t}$. In other words, we have defined $(Z_t^z)_{t\in(0,t_0)}$ to be a Markov bridge from $z$ to $X^z$. Then, set $Z_t^z = H_{t-t_0}(X^z,U')$ for all $t\geqslant t_0$ and $z=x,y$, where $U'$ is uniformly distributed  over $[0,1]$ and independent from $(Z_t^x,Z_t^y)_{t\in[0,t_0]}$, and $H$ is a representation of $z \mapsto (V_t^z)_{t\geqslant 0}$. By the Markov property, $(Z_t^z)_{t\geqslant 0}$ is then a Markov process with initial condition $Z_0^z = z$ and such that $Z_{t_0}^z = X^z$.

Define a second coupling $(\tilde Z_t^x,\tilde Z_t^y)_{t\geqslant 0}$ as follows. Let $\tau = \inf\{t\geqslant 0\ :\ Z^x_t = Z^y_t\}$. For $z=x,y$,  set $\tilde Z_t^z = Z_t^z$ for all $t\leqslant \tau$ and $\tilde Z_t^z = Z_t^x$ for $t\geqslant \tau$. By the strong Markov property, $(\tilde Z_t^z)_{t\geqslant 0}$ is a Markov process associated to $(P_t)_{t\geqslant 0}$ with $\tilde Z^z_0=z$, and the two other conditions are satisfied.
\end{proof}

Examples of couplings such that \eqref{Eq-LemCoupl1} holds and \eqref{Eq-LemCoupl2} does not are easily constructed for discrete-time Markov chains. Here is an example with continuous time. Consider the process $(X_t)_{t\geqslant 0}$ on $[-1,1]$ with generator
\[L f(x) = sign(x) f'(x) +  \frac{f(1)+f(-1)}{2} - f(0)\,.\]
In other words, starting from $x\in [-1,1]\setminus\{0\}$, then $X_t=x-sign(x) t$ for $t\in [0,|x|]$, after which $X_t=0$ for all $t\in[|x|,|x|+T)$ where $T$ is an exponential r.v. with parameter  1 and $X_{|x|+T}=R$ where $R$ is a Rademacher r.v. with parameter $1/2$. It is clear that if $|x|=|y|$ then $\delta_x P_t=\delta_y P_t$ for all $t\geqslant |x|$, since $\delta_x P_{|x|}=\delta_0$. Now consider for some $t_0>2$ the processes $(X_t,Y_t)_{t\geqslant 0}$ that jump at the same time $|x|+T$, the first one to a Rademacher r.v. $R$, the second one to  $R'=R(1-2\1_{|x|+T<t_0-1})$, which is indeed a Rademacher r.v. independent from $T$. In other words, as long as the jump occurs early enough for the processes to have the time to come back to 0 before time $t_0$, then we send them to opposite points, else to the same one. Obviously \eqref{Eq-LemCoupl1} holds and \eqref{Eq-LemCoupl2} does not.


\subsection{Study of the non-linear process}\label{Sec-Preuves-NL}

We use in this section the notations of Section \ref{Section-def-non-lineaire}. In particular we consider a Markov semi-group $(P_t)_{t\geqslant0}$ and non-linear jump rate and kernel $\nu\mapsto \lambda_\nu$, $Q_\nu$.

\subsubsection{Coupling with  different inhomogeneous jump mechanisms}\label{Section-Couplage-NL}

{Before giving the formal construction of the coupling, let us briefly highlight the general strategy: we want to define simultaneously two processes, driven by different time-inhomogeneous jump mechanisms, in such a way that, as much as possible, they are equal at some times. To do so, we use the linear Markov part of the dynamics in order to merge the two processes with some positive probability under a Doeblin condition (provided that no non-linear jump occurs meanwhile). To define the non-linear jump, we force the two processes to jump as much as possible simultaneously and, as much as possible, to the same location. Thanks to the Lipschitz assumptions on the jump mechanismsn the probability to achieve this is controlled by the distance between the time-dependent measures associated with the processes.}

Consider $\mu^1,\mu^2\in\MRP$. For all $x,y\in E$, $t\geqslant 0$ and $i=1,2 $, denote 
\[\tilde Q_t^i(x) \ = \ \frac{\lambda_{\mu_t^i}(x)}{\lambda_*} Q_{\mu_t^i}(x) + \po 1 - \frac{\lambda_{\mu_t^i}(x)}{\lambda_*}\pf \delta_x\,,\]
and $p_t(x,y) = {(\tilde  Q_{t}^1(x)\wedge\tilde  Q_{t}^2(y))(E)}$. If $p_t(x,y)= 0$, set
\[Q^=_t(x,y) = \delta_x \,,\qquad  Q_t^{i,\neq} (x,y) = \tilde Q_t^i(x) \,.\]
If $p_t(x,y)\in(0,1)$, set
\[Q_t^=(x,y) = \frac{\tilde Q_t^1(x)\wedge\tilde Q_t^2(y)}{p_t(x,y)}\,,\qquad   Q_t^{i,\neq}(x,y) = \frac{\tilde Q_t^i(x) - \tilde Q_t^1(x)\wedge \tilde Q_t^2(y)}{1- p_t(x,y)}\,.\]
Finally, if $p_t(x,y)=1$, set
\[Q_i^=(x,y) =  {\tilde Q_t^i(x)}\,,\qquad  Q_t^{i,\neq}(x,y) = \delta_x \,.\]
Then $Q_t^{=}$ and $Q_t^{i,\neq}$ for $i=1,2$ are Markov kernels from $E^2$ to $\mathcal P(E)$ such that, for all $x,y\in E$ and $i=1,2$,
\begin{equation}\label{eq:tildeQti}
\tilde  Q_t^i(x)\ = \ p_t(x,y) Q_t^=(x,y) + \po 1-p_t(x,y)\pf  Q_t^{i,\neq} (x,y)\,,
\end{equation}
and, as shown in the proof of Lemma \ref{Lem-couplage}, $\|\tilde Q_t^1(x) - \tilde Q_t^2(y)\|_{TV} = 2(1-p_t(x,y))$, so that the condition \eqref{Eq-Hyp-NonLin-Contract} of  Assumption~\ref{Hyp-MeanField-Exist} reads
\begin{eqnarray}\label{Eq-Hyp-pt}
\forall x\in E,\ \forall t\geqslant 0\,,\qquad 1-p_t(x,x) & \leqslant & \frac\theta{2\lambda_*} \|\mu_t^1 - \mu_t^2\|_{\V}\,,
\end{eqnarray}
and similarly with $\V = 1$ for Assumption \ref{Hyp-taux_Lipschitz}. For all $t\geqslant 0$ and $i=1,2$ let $  G_t^=$ (resp. $  G_t^{i,\neq}$) be a representation of  $Q_t^{=}$ (resp. $Q_t^{i,\neq}$), in the sense that for all $x,y\in E$, if $U$ is a r.v. uniformly distributed over $[0,1]$, then $G_t^=(x,y,U)\sim Q_t^=(x,y)$.

{
The meaning of the notations introduced here, namely $p_t, Q_t^{=}$ and $Q_t^{i,\neq}$ lies in the decomposition \eqref{eq:tildeQti}  of $\tilde Q_t^i(x)$, the latter being the law of the position after a jump at time $t$ from the position $x$ (except that, with respect to the initial kernel $Q_{\mu_t^i}$, we have artificially increased the jump rate to the upper bound $\lambda_*$, at the cost of adding ``phantom jumps" which are jumps from $x$ to $x$). The decomposition \eqref{eq:tildeQti} gives the following way to sample a position according to $\tilde Q_t^1(x)$ (when the other process is at a position $y$): with probability $p_t(x,y)$, draw a position according to $Q_t^=(x,y)$, otherwise draw a position according to $Q_t^{1,\neq}(x,y)$. The same goes for $\tilde Q_t^2(x)$ and since  $Q_t^=$ is the same for $i=1,2$ we can take the same random variable for the two processes in that case. In other words, $Q_t^=$ is the law of the common position of the process after a synchronous jump, while $Q_t^{i,\neq}$ for $i=1,2$ are the law of the positions when the coupling fails. 
}

\medskip


We define a Markov process $(X_t,Y_t)_{t\geqslant 0}$ on $E^2$ as follows. Let $m_0,h_0\in\mathcal P(E)$ and let $(X_0,Y_0)$ be an optimal coupling of $m_0$ and $h_0$ such as constructed in the proof of Lemma~\ref{Lem-couplage}. For a given $t_0$,  let $(Z_{t},\tilde Z_{t})_{t\geqslant 0}$ be an optimal coupling of $(\delta_{X_0} P_t)_{t\geqslant 0}$ and $(\delta_{Y_0} P_t)_{t\geqslant 0}$ in the sense of Lemma \ref{Lem-CouplageMarkov}. Consider an i.i.d. sequence $(S_k,U_k,V_k,W_k)_{k\in\mathbb N}$, independent from $(Z_{t},\tilde Z_{t})_{t\geqslant 0}$ where, for all $k\in\N$, $S_k$, $U_k$, $V_k$ and $W_k$  are independent one from the other, $S_k$ follows a standard exponential law with parameter 1 and $U_k$, $V_k$ and $W_k$   a uniform distribution over $[0,1]$.  Set $T_0 = 0$, $T_{n+1} = T_n +   S_n / \lambda_*$. for all $n\in\N$  and suppose that $(X_{t},Y_t)_{t\in [0,T_n]}$ have been defined for some $n\in \mathbb N$ and are independent from $(S_k,{U_{k}},V_k,W_k)_{k\geqslant n}$.

\begin{itemize}
\item Set $(X_{t},Y_{t}) = (Z_{t},\tilde Z_{t}) $ for all $t\in[0,T_{1})$ and  $(\tilde X_{T_{1}},\tilde Y_{T_{1}}) = (Z_{T_{1}},\tilde Z_{T_{1}}) $. If $n\geqslant 1$, for all $t\in(T_n,T_{n+1})$, set 
\[ (X_{t},Y_{t}) =  \po H_{t-T_n}(X_{T_n},U_{n}),H_{t-T_n}(Y_{T_n},U_{n})\pf \]
and
\[ (\tilde X_{T_{n+1}},\tilde Y_{T_{n+1}}) = \po H_{T_{n+1}-T_n}(X_{T_n},U_{n}),H_{T_{n+1}-T_n}(Y_{T_n},U_{n})\pf \,.\]
\item If $V_n \leqslant p_{T_n}\po \tilde X_{T_{n+1}}, \tilde Y_{T_{n+1}}\pf  $  then set 
\[X_{T_{n+1}} \ = \ Y_{T_{n+1}} \ = \ G_{T_n}^=\po \tilde X_{T_{n+1}}, \tilde Y_{T_{n+1}},W_n\pf\,, \]
and else set 
\[X_{T_{n+1}} \ = \  G_{T_n}^{1,\neq}\po \tilde X_{T_{n+1}}, \tilde Y_{T_{n+1}},W_n\pf\qquad Y_{T_{n+1}} \ = \  G_{T_n}^{2,\neq}\po \tilde Y_{T_{n+1}}, \tilde X_{T_{n+1}},W_n\pf\,. \] 
\end{itemize}

\begin{prpstn}\label{Prop-CouplNonLin-Ok}
For all $t\geqslant 0$, $X_t\sim m_0 R_{0,t}^{\mu^1}$ and $Y_t\sim h_0 R_{0,t}^{\mu^2}$.
\end{prpstn}

\begin{proof}
By symmetry, we only consider the case of $(X_t)_{t\geqslant 0}$. Note that in  the definition  of $(X_t)_{t\geqslant 0}$, as in the definition of $(X_t^{\mu^1,x})_{t\geqslant 0}$ in Section \ref{Section-def-non-lineaire}, $(T_n)_{n\geqslant 0}$ is a Poisson process with intensity $\lambda_*$. Moreover, for all $n\geqslant 0$, conditionally to $X_{T_n}$, $(X_t)_{t\in[T_n,T_{n+1})}$ is a Markov process associated to $(P_t)_{t\geqslant 0}$, i.e.  has the same distribution as $(Z_t^{X_{T_n}})_{t\in [0,T_{n+1}-T_n)}$, and the same goes for $(X_t^{\mu^1,x})_{t\in[T_n,T_{n+1})}$. As a consequence, it only remains to check that the Markov chain $(X_{T_n})_{n\in\N}$ has the same distribution as the Markov chain $(X_{T_n}^{\mu^1,x})_{n\in\N}$, which is equivalent to say that they have the same transition kernel. Yet, for any bounded measurable $f$ on $E$,
\begin{eqnarray*}
\mathbb E \po f(X_{T_1}) \ |\ T_1=t,\ \tilde X_{T_1} = \tilde  x,\ \tilde Y_{T_1} =\tilde  y\pf  & = & \mathbb E \big[ \1_{V_0 \leqslant p_{t}\po \tilde  x, \tilde  y\pf }  f(G_{t}^=\po \tilde x, \tilde y,W_0\pf ) \\
& & \ +  \ \1_{V_0 > p_{t}\po \tilde  x, \tilde  y\pf }  f(G_{t}^{1,\neq}\po \tilde x, \tilde y,W_n\pf )  \big]\\
& = & p_{t}\po \tilde  x, \tilde  y\pf  Q_t^=(\tilde x,\tilde y) + (1 -p_{t}\po \tilde  x, \tilde  y\pf) Q_t^{1,\neq}(\tilde x,\tilde y)\\
& = & \tilde Q^1_t f (\tilde x)\,.
\end{eqnarray*}
Hence,
\begin{eqnarray*}
\mathbb E \po f(X_{T_1}) \pf & =& \int_0^{+\infty} P_t^1 \tilde Q^1_t f(x) \lambda_* e^{-\lambda_* t} \dd t\\
& =& \int_0^{+\infty} P_t^1 \po \lambda_{\mu_t^i}   Q_{\mu_t^i} f + \po \lambda_* - \lambda_{\mu_t^i}\pf f \pf e^{-\lambda_* t} \dd t \,.\end{eqnarray*}
Similarly,
\begin{eqnarray*}
\mathbb E \po f(X_{T_1}^{\mu^1,x}) \ |\ T_1=t,\ H_{E_0 / \lambda_*}(x ,U_0) = \tilde  x\pf  & = & \mathbb E \big[ \1_{V_0 \leqslant \lambda_{\mu^1_t}(\tilde x)/\lambda_* }  f(G_{\mu^1_t} \po \tilde x,W_0\pf ) \\
& & \ +  \ \1_{V_0 > \lambda_{\mu^1_t}(\tilde x)/\lambda_* }  f(\tilde x) \big]\\
& = & \frac{\lambda_{\mu^1_t}(\tilde x)}{\lambda_* }Q_{\mu^1_t} f(x) + \po 1 - \frac{\lambda_{\mu^1_t}(\tilde x)}{\lambda_* }\pf f(\tilde x)\\
& = & \tilde Q^1_t f (\tilde x)\,,
\end{eqnarray*}
and again 
\begin{eqnarray*}
\mathbb E \po f(X_{T_1}^{\mu^1,x}) \pf & =& \int_0^{+\infty} P_t^1 \po \lambda_{\mu_t^i}   Q_{\mu_t^i} f + \po \lambda_* - \lambda_{\mu_t^i}\pf f \pf e^{-\lambda_* t} \dd t \,,
\end{eqnarray*}
which concludes.
\end{proof}

\subsubsection{The total variation case}\label{Section-NLTV}

\begin{proof}[Proof of Proposition \ref{Prop-ExistUniqueNonLin}]
For $(\mu^1_t,\mu^2_t)_{t\geqslant 0}$, $m_0,h_0\in\mathcal P(E)$, consider the associated   process $(X_t,Y_t)_{t\geqslant 0}$  defined in section \ref{Section-Couplage-NL}.  From Proposition \ref{Prop-CouplNonLin-Ok}, for all $t\geqslant 0$,
\begin{eqnarray*}
\frac12 \| m_0 R_{0,t}^{\mu^1} - h_0R_{0,t}^{\mu^2}\|_{TV}  & \leqslant  & \mathbb P\po X_t \neq Y_t\pf \\
& \leqslant & \mathbb P\po X_0 \neq Y_0\pf  + \mathbb P\po \tau_{split} < t\ | \ X_0 = Y_0\pf 
\end{eqnarray*}
where $\tau_{split}  =   \inf \{ T_n\ :\ n\in \N,\ V_n \geqslant p_{T_n}( \tilde X_{T_{n+1}}, \tilde Y_{T_{n+1}})  \}$. Note that, by construction, $X_0=Y_0$ implies that $X_t=Y_t$ for all $t< \tau_{split}$ so that, conditionally to $X_0=Y_0$, 
\begin{eqnarray*}
\tau_{split} &  \geqslant & \inf \{ T_n\ :\ n\in \N,\ 1- V_n \leqslant \theta \|\mu_{T_n}^1 - \mu_{T_n}^2\|_{TV}/(2\lambda_*)  \}\ : =\ \tilde \tau_{split}\,,
\end{eqnarray*}
where we used \eqref{Eq-Hyp-pt}. Note that $\tilde \tau_{split}$ is independent from $(X_0,Y_0)$, so that
\begin{eqnarray*}
\mathbb P\po \tau_{split} < t\ | \ X_0 = Y_0\pf & \leqslant  & \mathbb P\po \tilde \tau_{split} < t\pf\\
& = & 1 - \exp\po-\frac\theta 2\int_0^t \|\mu_s^1-\mu^2_s\|_{TV}\dd s \pf \,.
\end{eqnarray*}
Since $(X_0,Y_0)$ is an optimal coupling of $m_0$ and $h_0$, using that $1-\exp(-a)\leqslant a$ for all $a\in\R$, we have thus obtained
\begin{eqnarray}\label{Eq-NLContractTV}
 \| m_0 R_{0,t}^{\mu^1} - h_0R_{0,t}^{\mu^2}\|_{TV}  & \leqslant  &  \  \| m_0  - h_0 \|_{TV}  + \theta\int_0^t \|\mu_s^1-\mu^2_s\|_{TV}\dd s\,.
\end{eqnarray}
Now, remark that, for fixed $m_0\in\mathcal P(E)$, $t_1\geqslant 0$, $(m_0 R_{0,t}^{\mu})_{t\in[0,t_1]}$ only depends on $({\mu_t})_{t\in[0,t_1]}$. Set $t_1 = 1/(2\theta)$ and consider $\Psi : ({\mu_t})_{t\in[0,t_1]}\rightarrow(m_0 R_{0,t}^{\mu})_{t\in[0,t_1]}$. Then
\[\sup_{s\in[0,t_1]} \|m_0 R_{0,s}^{\mu^1} - m_0 R_{0,s}^{\mu^2}\|_{TV} \ \leqslant \ \frac12 \sup_{s\in[0,t_1]} \|\mu^1_s-\mu^2_s\|_{TV}\,, \]
In other words $\Psi$ is a contraction of $L^\infty\po [0,t_1],\mathcal P(E)\pf$, which is complete, hence admits a unique fixed point, which is a solution of \eqref{Eq-EquaNonLin} in the sense of Definition \ref{Def-solution} on $[0,t_1]$. Then there exists a unique solution on $[nt_1,(n+1)t_1]$ for all $n\in \N$, thus on $\R_+$.

Considering two such solutions with respective initial distributions $m_0$ and $h_0$, \eqref{Eq-NLContractTV} reads 
\[ \| m_t - h_t\|_{TV}  \ \leqslant  \  \  \| m_0  - h_0 \|_{TV}  + \theta\int_0^t \|m_s-h_s\|_{TV}\dd s\,,
\]
and Gronwall's Lemma concludes.
\end{proof}

\begin{proof}[Proof of Theorem \ref{Thm-TVNonLin}]
Let $t\mapsto m_t,h_t$ be two solutions of \eqref{Eq-EquaNonLin}. Consider again the proces $(X_t,Y_t)$ defined in Section \ref{Section-Couplage-NL}, with $\mu^1 = m$ and $\mu^2 = h$. In particular,
\begin{eqnarray*}
\frac12 \| m_{t_0} - h_{t_0}\|_{TV}  & \leqslant  & \mathbb P\po X_{t_0} \neq Y_{t_0} \pf \\
& \leqslant & \mathbb P\po  X_{t_0} \neq Y_{t_0},\ X_0 \neq Y_0\pf  + \mathbb P\po X_{t_0} \neq Y_{t_0}\ | \ X_0 = Y_0\pf\,.
\end{eqnarray*}
In the proof of Proposition \ref{Prop-ExistUniqueNonLin}, we established that
\begin{eqnarray*}
2\mathbb P\po X_{t_0} \neq Y_{t_0}\ | \ X_0 = Y_0\pf & \leqslant &  \theta \int_0^t \| m_s - h_s\|_{TV} \dd s\\
& \leqslant &  \theta \int_0^t e^{\theta s} \| m_0 - h_0\|_{TV} \dd s\ = \ (e^{\theta t_0}-1) \| m_0 - h_0\|_{TV} \,.
\end{eqnarray*}
By construction, for all $t<T_1$, $(X_t,Y_t)=(Z_t,\tilde Z_t)$, so that
\begin{eqnarray*}
\mathbb P\po  X_{t_0} = Y_{t_0}\ |\  X_0 \neq Y_0\pf & \geqslant & \mathbb P\po  Z_{t_0} = \tilde Z_{t_0},\ {T_1>t_0}\ |\  X_0 \neq Y_0\pf\\
& = & \mathbb P\po  Z_{t_0} = \tilde Z_{t_0} \ |\ Z_{0} \neq \tilde Z_{0} \pf \mathbb P ({T_1>t_0}) \ \geqslant \alpha e^{-\lambda_* t_0}\,.
\end{eqnarray*}
Since $\mathbb P(X_0 \neq Y_0) = \|m_0-h_0\|_{TV}/2$, we have obtained 
\begin{eqnarray*}
\| m_{t_0} - h_{t_0}\|_{TV}  & \leqslant  & \po e^{\theta t_0} -\alpha e^{-\lambda_* t_0} \pf \|m_0-h_0\|_{TV} \,.
\end{eqnarray*}
Now, for any $t\geqslant 0$, using again Proposition \ref{Prop-ExistUniqueNonLin}, we get
\begin{eqnarray*}
\| m_{t} - h_{t}\|_{TV}  & \leqslant  & e^{\theta (t-\lfloor t/t_0 \rfloor t_0)} \| m_{\lfloor t/t_0 \rfloor t_0} - h_{\lfloor t/t_0 \rfloor t_0}\|_{TV} \\
& \leqslant & e^{\theta t_0} \po e^{\theta t_0} -\alpha e^{-\lambda_* t_0} \pf^{\lfloor t/t_0 \rfloor } \|m_0-h_0\|_{TV}\,.
\end{eqnarray*}
\end{proof}

\subsubsection{The general $\V$ norm case}\label{Section-NL-Gene}


\begin{lmm}\label{Lem_NL_x=y}
Under Assumptions \ref{Hyp-taux_borne} and \ref{Hyp-MeanField-Exist}, for all $\mu^1,\mu^2\in \MRP$, $x\in E$ and $t\geqslant 0$,
\begin{eqnarray*}
\|\delta_x R_{0,t}^{\mu^1}- \delta_x R_{0,t}^{\mu^2}\|_\V & \leqslant &  \theta   \gamma_*  \V(x) e^{\po \rho_* + \lambda_*(\gamma_* -1)\pf t}\int_0^{t} \|\mu^1_{u} - \mu^2_{u}\|_{\V}  \dd u   \,.
\end{eqnarray*}
\end{lmm}

\begin{proof}
Consider the process $(X_t,Y_t)_{t\geqslant 0}$ defined in Section \ref{Section-Couplage-NL} with $X_0=Y_0=x$, and $t_1\geqslant 0$. As in the proof of Proposition \ref{Prop-ExistUniqueNonLin}, let $\tau_{split}  =   \inf \{ T_n\ :\ n\in \N,\ V_n \geqslant p_{T_n}( \tilde X_{T_{n+1}}, \tilde Y_{T_{n+1}})  \}$. Since $X_t=Y_t$ for all $t< \tau_{split}$,
\begin{eqnarray*}
\tau_{split} &  \geqslant & \inf \{ T_n\ :\ n\in \N,\ 1- V_n \leqslant \theta \|\mu^1_{T_n} - \mu^2_{T_n}\|_{\V}/(2\lambda_*)  \}\ : =\ \tilde \tau_{split}\,,
\end{eqnarray*}
and we bound
\begin{eqnarray*}
\mathbb E \po   d_\V \po X_{t_1},Y_{t_1}\pf \pf & \leqslant &\mathbb E\po \1_{\tilde \tau_{split} \leqslant t_1} \po    \V(X_{t_1})+ \V(Y_{t_1})\pf \pf\,.
\end{eqnarray*}
Set $\Gamma = (S_k,V_k)_{k\in \N}$, and remark that $\tilde \tau_{split}$ and $K=\sharp\{n\in\N\ : \ T_n \leqslant t_1 \}$ are deterministic functions of $\Gamma$, while $(Z_t,\tilde Z_t)_{t\geqslant 0}$ and $(U_k,W_k)_{k\in\N}$ are independent from $\Gamma$. Using \eqref{Eq-Hyp-NonLin-Gene-1} and \eqref{Eq-Hyp-NonLin-Gene-4}, for all $n\in \N$,
\begin{eqnarray*}
\mathbb E \po   \V(X_{T_{n}})  \ |\ \Gamma \pf & \leqslant &  \gamma_*     \mathbb E \po   \V(\tilde X_{T_{n}})  \ |\ \Gamma \pf \\
& = &  \gamma_*  \mathbb E \po   P_{T_{n}-T_{n-1}} \V(X_{T_{n-1}})  \ |\ \Gamma \pf \\
& \leqslant &   \gamma_* e^{\rho_*( T_{n}-T_{n-1})}  \mathbb E \po    \V(X_{T_{n-1}})  \ |\ \Gamma \pf \,.
\end{eqnarray*}
Similarly, and then dy direct induction,
\begin{eqnarray*}
\mathbb E \po   \V(X_{t_1})  \ |\ \Gamma \pf & \leqslant &     e^{\rho_*( t_1-T_{K})}  \mathbb E \po    \V(X_{T_{K}})  \ |\ \Gamma \pf \ \leqslant \   \gamma_*^{K}e^{\rho_* t_1} \V(x)  \,.
\end{eqnarray*}
The case of $Y_{t_1}$ being identical, with $Y_0=X_0=x$, we have obtained 
\begin{eqnarray*}
\mathbb E \po   d_\V \po X_{t_1},Y_{t_1}\pf \pf & \leqslant & 2 e^{\rho_* t_1}\V(x) \mathbb E\po \1_{\tilde \tau_{split} \leqslant t_1}   \gamma_*^{K}   \pf\,.
\end{eqnarray*}
 Besides, $(V_k)_{k\in\N}$ is independent from $K$, and conditionnally to $K$, $\{T_i\}_{i\in \cco 1,K\ccf }$ are distributed like $K$ independent uniform r.v. over  $[0,t_1]$ and independent from $(V_k)_{k\in\N}$ , so that 
\begin{eqnarray*}
\mathbb P \po \tilde \tau_{split} > t_1 \ | \ K \pf & \geqslant & \mathbb P \po 1-V_n \geqslant \theta  \|\mu^1_{T_{n+1}} - \mu^2_{T_{n+1}}\|_{\V}/(2\lambda_*)\,, \quad \forall n\in \cco 0,K-1\ccf \ | \ K  \pf \\
& = & \po 1 -\frac1{t_1} \int_0^{t_1}   1 \wedge \po \frac\theta{2\lambda_*} \|\mu^1_{u} - \mu^2_{u}\|_{\V}\pf \dd u \pf^K\ := \ s^K\,.
\end{eqnarray*}
Since $K$ follows a Poisson law with parameter $\lambda_* t_1$, 
\begin{eqnarray*}
\mathbb E \po   d_\V\po X_{t_1},Y_{t_1}\pf \pf & \leqslant & 2 e^{\rho_* t_1}  \V(x) \mathbb E\po (1-s^K)      \gamma_*^{K}  \pf\\
&=& 2 e^{\rho_* t_1}  \V(x) \po e^{\lambda_*t_1(\gamma_* -1)} - e^{\lambda_*t_1(\gamma_* s-1)}\pf\\
& \leqslant & 2\lambda_*t_1\gamma_* (1-s)      \V(x) e^{\po \rho_* + \lambda_*(\gamma_* -1)\pf t_1}  \,,
\end{eqnarray*}
which concludes.
\end{proof}

\begin{proof}[Proof of Theorem \ref{Thm-NonLin-Exist}]
Let $t\mapsto \mu^1_t,\mu^2_t$ be two measurable functions from $\R_+$ to $\mathcal P_\V$, and $m_0\in\mathcal P_\V$. From Lemma~\ref{Lem_NL_x=y}, 
\begin{eqnarray*}
\| m_0 R_{0,t}^{\mu^1} - m_0 R_{0,t}^{\mu^2}\|_\V & = & \int  \| \delta_{x} R_{0,t}^{\mu^1} - \delta_{x} R_{0,t}^{\mu^2}\|_\V m_0(\dd x) \\
& \leqslant & \theta  \gamma_*  m_0(\V) e^{\po \rho_* + \lambda_*(\gamma_* -1)\pf t}\int_0^{t} \|\mu^1_{u} - \mu^2_{u}\|_{\V}  \dd u   \,.
\end{eqnarray*}
Moreover, if $\mu_t^1(\V) \leqslant m_0(\V) \vee (M/(1-\eta))$ for all $t\geqslant 0$, then \eqref{Eq-Hyp-NonLin-Lyap} implies that
\[m_0R_{0,t}^{\mu^1}(\V) \ \leqslant \ e^{-\rho t} m_0(\V) + (1-e^{-\rho t}) \po \eta \po m_0(\V)\vee \frac{M}{1-\eta}\pf + M\pf \ \leqslant \ m_0(\V)\vee \frac{M}{1-\eta}\,.\]
Let $t_1$ be small enough so that 
\begin{eqnarray*}
t_1 \theta   \gamma_*  \po m_0(\V)\vee \frac{M}{1-\eta}\pf  e^{\po \rho_* + \lambda_*(\gamma_* -1)\pf t_1} & \leqslant & \frac12\,.
\end{eqnarray*}
Let $\mathcal A$ be the set of measurable functions $t\mapsto \mu_t$ from $[0,t_1]$ to $\mathcal P_\V$ such that $\mu_t(\V) \leqslant m_0(\V)\vee (M/(1-\eta))$ for all $t\in[0,t_1]$. Remark that $\mathcal A$ is not empty since it contains the function with constant value $m_0$. As a closed subset of $L^\infty ( [0,t_1],\mathcal P_\V)$, $\mathcal A$ is complete with the distance
\[dist_{\mathcal A}(\mu^1,\mu^2) \ = \ \sup_{s\in[0,t_1]} \|\mu^1_s- \mu_s^2\|_\V\,.\]
The map $\Psi : (\mu_t)_{t\in[0,t_1]} \rightarrow (m_0R_{0,t}^\mu)_{t\in[0,t_1]}$ being a contraction of $\mathcal A$, it admits a unique fixed point in $\mathcal A$, which is then a solution of \eqref{Eq-EquaNonLin} over the time interval $[0,t_1]$. In particular, $m_{t_1} \leqslant m_0(\V) \vee (M/(1-\eta))$, and the same arguments (with the same time $t_1$) yield a solution on the time interval $[t_1,2t_1]$, and then by induction on $\R_+$.

If $m$ is such a solution, then \eqref{Eq-Hyp-NonLin-Lyap} integrated with respect to $m_0$ and applied with $\mu_s = m_s$ for all $s\geqslant 0$ yields
\[e^{\rho t} m_t(\V)  \ \leqslant \  m_0(\V) + \eta\rho \int_0^t e^{\rho s}  m_s(\V) \dd s + (e^{\rho t} - 1) M \ := \ j_t \]
with $j_t'\leqslant \eta \rho j_t + \rho M \exp(\rho t)$, hence
\[m_t(\V) \ \leqslant \ e^{-\rho t} j_t \ \leqslant \ e^{-\rho t}\po e^{\eta \rho t} j_0 + M \frac{e^{\rho t} - e^{\eta \rho t}}{1-\eta}\pf \,,\]
which is exactly \eqref{Eq-NL-Lyap}.
\end{proof}

A direct corollary of \eqref{Eq-Hyp-NonLin-Lyap}  and \eqref{Eq-NL-Lyap} is that, for all $t\mapsto m_t$ solution of \eqref{Eq-EquaNonLin} with $m_0(\V) \leqslant M/(1-\eta) + 1$, all $t\geqslant 0$ and all $x\in E$,
\begin{eqnarray}\label{Eq-LyapNL-Contract-V(x)}
R_{0,t}^m \V(x) & \leqslant & e^{-\rho t} \V(x) + \po 1-e^{-\rho t} \pf \po \frac{M}{1-\eta} +1\pf\,.
\end{eqnarray}

\begin{lmm}\label{Eq-NL-Expansion}
Under Assumptions \ref{Hyp-taux_borne} and \ref{Hyp-MeanField-Exist}, if $t\mapsto m_t,h_t$ are solutions of \eqref{Eq-EquaNonLin} with $h_0(\V)\vee m_0(\V) \leqslant M/(1-\eta) + 1$, then for all $t\geqslant 0$,
\begin{eqnarray*}
\| m_t- h_t\|_\V& \leqslant & e^{{\po\po\rho + \theta  \gamma_* \pf \frac{M}{1-\eta} + \theta  \gamma_*\pf } t} \| m_0- h_0\|_\V\,.
\end{eqnarray*}
\end{lmm}

\begin{proof}
Consider $(X_0,Y_0)$ an optimal coupling of $m_0$ and $h_0$ such as given by Lemma \ref{Lem-couplage} and let $(X_t,Y_t)$ be an optimal coupling of $\delta_{X_0} R_{0,t}^m$ and  $\delta_{Y_0} R_{0,t}^h$, so that
\begin{eqnarray*}
\| m_t - h_t\|_\V & \leqslant  & \mathbb E\po d_\V(X_t,Y_t)\pf \  = \ \mathbb E\po \po \1_{X_0=Y_0} + \1_{X_0\neq Y_0} \pf d_\V(X_t,Y_t)\pf
\end{eqnarray*}
Since $ \mathbb E(  \1_{X_0=Y_0} \V(X_0))  \leqslant m_0(\V) \wedge h_0(\V)$, from  Lemma~\ref{Lem_NL_x=y}, 
\begin{eqnarray*}
  \mathbb E\po    \1_{X_0 = Y_0}  {d_\V(X_t,Y_t)} \pf  & \leqslant & \theta  \gamma_* \po \frac{M}{1-\eta} +1\pf  e^{\po \rho_* + \lambda_*(\gamma_* -1)\pf t}\int_0^{t} \|m_{u} - h_{u}\|_{\V}  \dd u \\
  & := & a_1(t)\int_0^{t} \|m_{u} - h_{u}\|_{\V}  \dd u \,.
\end{eqnarray*}
From  \eqref{Eq-LyapNL-Contract-V(x)},
\begin{eqnarray*}
  \mathbb E\po    \1_{X_0 \neq Y_0}   \po \V(X_t)+\V(Y_t)\pf \pf  & \leqslant & e^{-\rho t} \mathbb E \1_{X_0 \neq Y_0}  \po \V(X_0)+ \V(Y_0)\pf \\
  & &  \ +\  2\po 1-e^{-\rho t} \pf \po \frac{M}{1-\eta} +1\pf \mathbb P\po X_0\neq Y_0\pf\\
  & \leqslant & \co e^{-\rho t} + \po 1-e^{-\rho t} \pf \po \frac{M}{1-\eta} +1\pf\cf \| m_0-h_0\|_\V \\
  & := &  a_2(t) \| m_0-h_0\|_\V\,.
\end{eqnarray*}
Thus
\begin{eqnarray*}
\| m_t - h_t\|_\V & \leqslant  &   a_2(t) \| m_0 - h_0\|_\V  + a_1(t) \int_0^t \| m_s - h_s\|_\V  \dd s\,.
\end{eqnarray*}
As in  the proof of Theorem \ref{Thm-NonLin-Exist}, the Gronwall's Lemma yields, for all $t\geqslant 0$,
\begin{eqnarray*}
\| m_t - h_t\|_\V & \leqslant  &    \po a_2(t) + a_1(t) \int_0^t a_2(s) e^{\int_s^t a_1(u)\dd u}\dd s \pf \| m_0 - h_0\|_\V\,.
\end{eqnarray*}
For $t\geqslant 0$ and $n\geqslant 0$, applying $n$ times this result with the time $t/n$ we get
\begin{eqnarray*}
\| m_t - h_t\|_\V & \leqslant  &   \po a_2(t/n) + a_1(t/n) \int_0^{t/n} a_2(s) e^{\int_s^{t/n} a_1(u)\dd u}\dd s \pf ^n \| m_0 - h_0\|_\V
\end{eqnarray*}
where we used that, from \eqref{Eq-NL-Lyap}, $m_{kt/n}(\V) \leqslant M/(1-\eta) + 1$ for all $k\in\N$. Letting $n$ go to infinity,
\begin{eqnarray*}
\| m_t - h_t\|_\V & \leqslant  &  e^{\po a_2'(0) + a_1(0) a_2(0)\pf t}\| m_0 - h_0\|_\V\,,
\end{eqnarray*}
which concludes.
\end{proof}

In the following, as in \cite{HairerMattingly2008}, we consider for a scale parameter $\beta>0$ the distance on $\mathcal P_\V(E)$ given by $\rho_\beta(\mu_1,\mu_2) = \| \mu_1 - \mu_2\|_{\beta + \V}$. In other words, denoting $d_\beta (x,y)  = \1_{x\neq y} \po 2\beta + \V(x) + \V(y)\pf$ and applying Lemma \ref{Lem-couplage} means that
\[\rho_\beta(\mu_1,\mu_2)\ =\ \inf\left\{\mathbb E\po d_\beta(X,Y)\pf\,:\,X\sim\mu_1,\ Y\sim \mu_2 \right \}\,.\]

\begin{lmm}\label{Lem_NL_HairerMatt}
Under Assumptions \ref{Hyp-taux_borne} and \ref{Hyp-MeanField-Exist}, if $\alpha>0$ and if $t\mapsto m_t,h_t$ are solutions of \eqref{Eq-EquaNonLin} with $h_0(\V)\vee m_0(\V) \leqslant M/(1-\eta) + 1$, then for all $x,y\in E$ such that $x\neq y$,
\begin{eqnarray*}
\rho_\beta \po \delta_x R_{0,t_0}^m, \delta_y R_{0,t_0}^h\pf & \leqslant  & \kappa d_\beta(x,y)
\end{eqnarray*}
with
\begin{eqnarray}\label{Eq-Def-Beta} 
\beta & =&    \frac{2 \po 1-e^{-\rho t_0} \pf e^{\lambda_* t_0} }{\alpha} \po \frac{M}{1-\eta} +1\pf  \,,\qquad  \kappa \ = \ \frac{1+e^{-\rho t_0}}{2} \vee \po  1-\frac12 \alpha e^{-\lambda_* t_0}\pf\,.
\end{eqnarray}
\end{lmm}

\begin{proof}
This is essentially the proof of \cite[Theorem 3.1]{HairerMattingly2008}. Consider the process $(X_t,Y_t)_{t\geqslant 0}$  defined in Section \ref{Section-Couplage-NL} with $X_0 = x$, $Y_0 = y$, and bound
\begin{eqnarray*}
\rho_\beta \po \delta_x R_{0,t_0}^m, \delta_y R_{0,t_0}^h\pf & \leqslant  &  \mathbb E \po   d_\beta \po X_{t_0},Y_{t_0}\pf \pf\\
 & \leqslant & 2\beta \mathbb P\po X_{t_0}\neq Y_{t_0}\pf +   R_{0,t_0}^m \V(x) +   R_{0,t_0}^h \V(y)\,.
\end{eqnarray*}
Then \eqref{Eq-LyapNL-Contract-V(x)} yields
\begin{eqnarray*}
R_{0,t_0}^m \V(x)  + R_{0,t_0}^h \V(y) & \leqslant & e^{-\rho t_0} \po \V(x) + \V(y)\pf + 2 \po 1-e^{-\rho t_0} \pf \po \frac{M}{1-\eta} +1\pf \,. 
\end{eqnarray*}
Two cases are distinguished. First, if $\V(x)+\V(y) \leqslant 8(1+M/(1-\eta))$, similarly to the proof of Theorem \ref{Thm-TVNonLin}, from \eqref{Eq-NonLin-Gene-Doeblin},
\[\mathbb P\po X_{t_0}=Y_{t_0}\pf\  \geqslant \ \mathbb P \po Z_{t_0} = \tilde Z_{t_0},\ T_1>t_0 \pf \ \geqslant \ \alpha e^{-\lambda_* t_0}\,,\] 
so that 
\begin{eqnarray*}
\mathbb E \po   d_\beta \po X_{t_0},Y_{t_0}\pf \pf & \leqslant & 2\beta \po 1-\alpha e^{-\lambda_* t_0}\pf +     e^{-\rho t_0} \po \V(x) + \V(y)\pf + 2 \po 1-e^{-\rho t_0} \pf \po \frac{M}{1-\eta} +1\pf\\
& = & 2\beta \po 1-\frac12 \alpha e^{-\lambda_* t_0}\pf +     e^{-\rho t_0} \po \V(x) + \V(y)\pf\,,
\end{eqnarray*}
where we used the definition \eqref{Eq-Def-Beta} of $\beta$. Second, if $\V(x)+\V(y) \geqslant 8(1+M/(1-\eta))$, then
\begin{eqnarray*}
\mathbb E \po d_{\beta}(X_{t_0},Y_{t_0})\pf & \leqslant & 2  \beta  +     e^{-\rho t_0} \po \V(x) + \V(y)\pf +2 \po 1-e^{-\rho t_0} \pf \po \frac{M}{1-\eta} +1\pf\\
& \leqslant & 2\beta  +     \frac{1+e^{-\rho t_0}}{2}\po \V(x) + \V(y)\pf  - 2 \po 1-e^{-\rho t_0} \pf \po \frac{M}{1-\eta} +1\pf \\
& = & 2\beta \po 1-\frac12 \alpha e^{-\lambda_* t_0}\pf +      \frac{1+e^{-\rho t_0}}{2} \po \V(x) + \V(y)\pf \,,
\end{eqnarray*}
where we used again the definitions \eqref{Eq-Def-Beta} of $\beta$.

\end{proof}

\begin{proof}[Proof of Theorem \ref{Thm_non-lin}]
Applying \eqref{Eq-Hyp-NonLin-Lyap} with $x=x_0$ and $\mu_s = \delta_{x_0}$ for all $s\geqslant 0$ where $x_0\in E$ is such that $\V(x_0)$ is arbitrarily close to $\inf \V$, we obtain that $\inf \V \leqslant M(1-\eta)$. In particular, $\PV:=\{\nu\in\mathcal P_\V,\ \nu(\V) \leqslant M(1-\eta) + 1\}$ is not empty since it contains $\delta_{x_0'}$ for all $x_0'\in E$ with $\V(x_0')$ sufficiently close to $\inf\V$. From \eqref{Eq-LyapNL-Contract-V(x)}, if $t\mapsto m_t$ is a solution of \eqref{Eq-EquaNonLin} with $m_0 \in \PV$ and then $m_t\in\PV$ for all $t\geqslant 0$. 

Let $m_0,h_0\in \PV$   and $(X_0,Y_0)$ be an optimal coupling of $m_0$ and $h_0$ given by Lemma \ref{Lem-couplage}, i.e. be such that, for any $\beta>0$,
\[\rho_\beta(m_0,h_0) \ = \ \mathbb E \po d_\beta(X_0,Y_0)\pf\,.\]
Conditionning on the initial value, let $(X_{t_0},Y_{t_0})$ be an optimal coupling of $\delta_{X_0} R_{0,t_0}^m$ and $\delta_{Y_0} R_{0,t_0}^h$, so that
\begin{eqnarray*}
\mathbb E \po d_\beta(X_{t_0},Y_{t_0})\pf & = & \mathbb E \po \rho_\beta\po \delta_{X_0} R_{0,t_0}^m,\delta_{Y_0} R_{0,t_0}^h\pf\pf\\
& = & \mathbb E \po \po \1_{X_0=Y_0} + \1_{X_0 \neq Y_0}\pf \rho_\beta\po \delta_{X_0} R_{0,t_0}^m,\delta_{Y_0} R_{0,t_0}^h\pf\pf\,.
\end{eqnarray*}
Considering $\kappa$ and $\beta$ such as defined in Lemma \ref{Lem_NL_HairerMatt} and using Lemma \ref{Lem_NL_x=y} and the fact that $d_\beta \leqslant (1+\beta) d_\V$,
\begin{eqnarray*}
\mathbb E \po d_\beta({X_{t_0},Y_{t_0}})\pf & \leqslant &    \kappa    \mathbb E \po d_\beta(X_{0},Y_{0})\pf\\
&  & \  +\ 2(1+\beta) \theta  \gamma_*  \po \frac{M}{1-\eta}+1\pf  e^{\po \rho_* + \lambda_*(\gamma_* -1)\pf t_0}\int_0^{t_0} \|m_{u}- h_{u}\|_\V  \dd u    \,.
\end{eqnarray*}
Together with Lemma \ref{Eq-NL-Expansion} and the fact that, $\|\cdot\|_\V \leqslant \rho_\beta$ (since $d_\V\leqslant d_\beta$), this means
\begin{eqnarray}
\notag \rho_{\beta }(m_{t_0}, h_{t_0}) & \leqslant & \rho_{\beta }(m_{0}, h_{0})   \co \kappa + \theta \frac{  2(1+\beta) \gamma_*   }{\po\rho + \theta  \gamma_* \pf  }  e^{\po \rho_* + \lambda_*(\gamma_* -1)\pf t_0}  \po e^{\po\rho + \theta  \gamma_* \pf \po\frac{M}{1-\eta} +1\pf t_0} - 1\pf \cf \\
& \leqslant & \tilde  \kappa \rho_{\beta }(m_{0}, h_{0}) \label{Eq-NL-Demo-Kappa} \,.
\end{eqnarray}
By assumption, $\tilde  \kappa<1$, so that
 $\Psi_{t_0}:\PV \rightarrow \PV$ that maps $m_0$ to $m_{t_0}$ where $t\mapsto m_t$ is a solution of \eqref{Eq-EquaNonLin} is a contraction. If a sequence $(\mu_n)_{n\in\N}$ in $\mathcal P_\V$ is such that $\rho_\beta(\mu_n,\nu) \rightarrow 0$ for some $\nu\in\mathcal P(E)$ as $n\rightarrow \infty$, then $\mu_n(\V)\rightarrow \nu(\V)$. Hence, $\PV$  is a closed subset of $\mathcal P_\V(E)$ endowed with the metric $\rho_\beta$, which is complete. As a consequence, $\Psi_{t_0}$ admits a unique fixed point $\mu_\infty$ in $\PV$, which for now may depend on $t_0$.

From \eqref{Eq-NL-Lyap},  for all $n\in\N$ and $m_0\in\mathcal P_\V(E)$,
\begin{eqnarray*}
m_{nt_0}(\V) & \leqslant & e^{-\rho(1-\eta) n t_0} m_0(\V) + \po 1 - e^{-\rho(1-\eta) n t_0}\pf \frac{M}{1-\eta}.
\end{eqnarray*}
which, applied to $m_0=\mu_\infty$ and letting $n$ go to infinity implies that  $\V(\mu_\infty)\leqslant M/(1-\eta)$. More generally,    \eqref{Eq-NL-Lyap} means that for all $m_0\in \mathcal P_\V(E)$ and all $t\geqslant s_0:= \ln( m_0(\V))/(\rho(1-\eta))$, $m_t\in\PV$.
 For all $\nu \in \PV$,
\[\rho_\beta(\mu_\infty,\nu) \ \leqslant \ 2\beta + 2 \frac{M}{1-\eta} + 2\,.\]
Combining this with \eqref{Eq-NL-Demo-Kappa}, we obtain that, for all $m_0\in\mathcal P_\V(E)$,
\begin{eqnarray*}
 \rho_\beta \po \mu_\infty ,m_{t}\pf & \leqslant &   \tilde  \kappa^{\lfloor (t-s_0)/t_0\rfloor}  \rho_\beta \po \mu_\infty ,m_{t -  t_0\lfloor (t-s_0)/t_0\rfloor}\pf \\
 & \leqslant & 2 \tilde \kappa^{ (t-s_0)/t_0 - 1}\po \beta + 1 + \frac{M}{1-\eta }\pf \,.
\end{eqnarray*}
Since $\tilde  \kappa \geqslant e^{-\rho (1-\eta) t_0}$, $\tilde  \kappa^{ -s_0/t_0 } \leqslant m_0(\V)$.
We have then obtained, for all $t\geqslant 0$ and all $m_0\in\mathcal P_\V$,
\begin{eqnarray*}
 \rho_\beta \po \mu_\infty ,m_{t}\pf  & \leqslant & 2 \tilde  \kappa^{ t/t_0 -1}\po \beta + 1 + \frac{M}{1-\eta}\pf   m_0(\V)\,.
\end{eqnarray*}
In particular, if $m_0 = \mu_\infty$, since in that case $m_s = h_{s+nt_0}$ for all $s\geqslant 0$ and $n\in \N$, letting $n$ go to infinity, we get that $m_s = \mu_\infty$ for all $s\geqslant 0$, in other words $\mu_\infty$ is an equilibrium of \eqref{Eq-EquaNonLin}. Finally, $\rho_\beta \geqslant \|\cdot\|_\V$, which concludes.
\end{proof}

\subsection{Study of the particle system}\label{Sec-Preuves-Particules}

We use in this section the notations of Section \ref{Section-Def-Particule}. In particular for $i\in\cco 1,N\ccf$ we consider a Markov semi-group $(P_t^i)_{t\geqslant0}$ and jump rate and kernel $\lambda_i$, $Q_i$.

\subsubsection{Coupling for interacting particles}\label{Section-Particules-Couplage}

%
%
%

For all $x,y\in E$ and $i\in\cco 1,N\ccf$, denote 
\[\tilde Q_i(x) \ = \ \frac{\lambda_i(x)}{\lambda_*} Q_i(x) + \po 1 - \frac{\lambda_i(x)}{\lambda_*}\pf \delta_{x_i}\,,\]
and $p_i(x,y) = (\tilde  Q_i(x)\wedge\tilde  Q_i(y))(E_i)$. If $p_i(x,y)= 0$, set
\[Q^=_i(x,y) = \delta_{x_i} \,,\qquad  Q_i^{\neq} (x,y) = \tilde Q_i(x) \,.\]
If $p_i(x,y)\in(0,1)$, set
\[Q_i^=(x,y) = \frac{\tilde Q_i(x)\wedge\tilde Q_i(y)}{p_i(x,y)}\,,\qquad   Q_i^{\neq}(x,y) = \frac{\tilde Q_i(x) - \tilde Q_i(x)\wedge\tilde  Q_i(y)}{1- p_i(x,y)}\,.\]
Finally, if $p_i(x,y)=1$, set
\[Q_i^=(x,y) =  \tilde Q_i(x)\,,\qquad  Q_i^{\neq}(x,y) = \delta_{x_i} \,.\]
Then $Q_i^{=}$ and $Q_i^{\neq}$ are Markov kernels from $E^2$ to $\mathcal P(E_i)$ such that, for all $x,y\in E$,
\[\tilde  Q_i(x)\ = \ p_i(x,y) Q_i^=(x,y) + \po 1-p_i(x,y)\pf  Q_i^{\neq} (x,y)\,,\]
and, as shown in the proof of Lemma \ref{Lem-couplage}, $\|\tilde Q_i(x) - \tilde Q_i(y)\|_{TV} = 2(1-p_i(x,y))$, so that \eqref{Eq-TauxLipschitz-Particules} reads
\[\forall x,y\in E\text{ such that } x_i=y_i\,,\qquad 1-p_i(x,y) \ \leqslant \ \frac{\theta \dist(x,y)}{2N\lambda_*}\,,\]
and similarly for \eqref{Eq-TauxLipschitz-Particules-V}. For all $i\in\cco 1,N\ccf$ let $  G_i^=:E^2\times[0,1]\rightarrow E_i$ (resp. $  G_i^{\neq}$) be a representation of  $Q_i^{=}$ (resp. $Q_i^{\neq}$), in the sense that for all $x,y\in E$, if $U$ is a r.v. uniformly distributed over $[0,1]$, then $G_i^=(x,y,U)\sim Q_i^=(x,y)$. 

Let $x,y\in E$. We define a Markov process $(X_t,Y_t)_{t\geqslant 0} = (X_{i,t},Y_{i,t})_{i\in\cco 1,N\ccf, t\geqslant 0}$ on $E^2$ and an auxiliary Markov process $(J_t)_{t\geqslant0}$ on $\cco 1,N\ccf$ as follows. For a given $t_0$ and all $i\in \cco 1,N\ccf$  let $(Z_{i,t},\tilde Z_{i,t})_{t\geqslant 0}$ be an optimal coupling of $(\delta_x P_t^i)_{t\geqslant 0}$ and $(\delta_y P_t^i)_{t\geqslant 0}$ in the sense of Lemma \ref{Lem-CouplageMarkov}, independent one from the other for $j\neq i$. Consider an i.i.d. sequence $(S_k,(U_{i,k})_{i\in\cco1,N\ccf},V_k,W_k,I_k)_{k\in\mathbb N}$, independent from  $(X_0,Y_0)$ and $(Z_{i,t},\tilde Z_{i,t})_{t\geqslant 0,i\in\cco 1,N\ccf}$ and where, for all $k\in\N$, $S_k$, $V_k$, $W_k$, $I_k$ and $(U_{i,k})_{i\in\cco1,N\ccf}$ are independent one from the other, $S_k$ follows a standard exponential law with parameter 1, $V_k$ and $W_k$ (resp. $I_k$) a uniform distribution over $[0,1]$ (resp. over $\cco 1,N\ccf$), and $(U_{i,k})_{i\in\cco1,N\ccf}$ are i.i.d. r.v. uniformly distributed over $[0,1]$.  Set $T_0 = 0$, $X_0=x$, $Y_0=y$ and $J_0 = \dist(x,y)/2$ and suppose that $T_n\geqslant 0$ and $(X_{t},Y_t,J_t)_{ t\in [0,T_n]}$ have been defined for some $n\in \mathbb N$ and are independent from $(S_k,(U_{i,k})_{i\in\cco1,N\ccf},V_k,W_k,I_k)_{k\geqslant n}$. Set $T_{n+1} = T_n + (N\lambda_*)^{-1} S_n$. 
\begin{itemize}
\item For all $i\in  \cco 1,N\ccf$ and all $t\in(T_n,T_{n+1})$, set $J_t = J_{T_n}$,
\[ (X_{i,t},Y_{i,t}) = \left\{\begin{array}{l}
(Z_{i,t},\tilde Z_{i,t}) \qquad \text{if } (X_{i,T_n},Y_{i,T_n})=(Z_{i,T_n},\tilde Z_{i,T_n})\\
\po H^i_{t-T_n}(X_{i,T_n},U_{i,n}),H^i_{t-T_n}(Y_{i,T_n},U_{i,n})\pf \qquad \text{else, } 
\end{array}\right.\]
and
\[ (\tilde X_{i,T_{n+1}},\tilde Y_{i,T_{n+1}}) = \left\{\begin{array}{l}
(Z_{i,T_{n+1}},\tilde Z_{i,T_{n+1}}) \qquad \text{if } (X_{i,T_n},Y_{i,T_n})=(Z_{i,T_n},\tilde Z_{i,T_n})\\
\po H^i_{T_{n+1}-T_n}(X_{i,T_n},U_{i,n}),H^i_{T_{n+1}-T_n}(Y_{i,T_n},U_{i,n})\pf \qquad \text{else. } 
\end{array}\right.\]
\item If $V_n \leqslant p_{I_n}\po \tilde X_{T_{n+1}}, \tilde Y_{T_{n+1}}\pf  $  then set 
\[X_{T_{n+1}} \ = \ Y_{T_{n+1}} \ = \ G_{I_n}^=\po \tilde X_{T_{n+1}}, \tilde Y_{T_{n+1}},W_n\pf\,, \]
and else set 
\[X_{T_{n+1}} \ = \  G_{I_n}^{\neq}\po \tilde X_{T_{n+1}}, \tilde Y_{T_{n+1}},W_n\pf\qquad Y_{T_{n+1}} \ = \  G_{I_n}^{\neq}\po \tilde Y_{T_{n+1}}, \tilde X_{T_{n+1}},W_n\pf\,. \]

\item If $V_n \geqslant 1- \theta J_{T_n}/(N\lambda_*)$ and $x_{I_n}=y_{I_n}$ then set $J_{T_{n+1}} = J_{T_n} + 1$, else set $J_{T_{n+1}} = J_{T_n}$. 
\end{itemize}
Then $(X_{t},Y_t,J_t)$ is defined for all $t\geqslant 0$.

\begin{prpstn}\label{Prop-Coupl-Particule-Ok}
For all $t\geqslant 0$, $X_t\sim \delta_x R_t$, $Y_t\sim\delta_y R_t $ and $2 J_t \geqslant \dist(X_t,Y_t)$.
\end{prpstn}
\begin{proof}
The proof of the first statements is similar  to the proof of Proposition \ref{Prop-CouplNonLin-Ok}, hence omitted. 
 For the last one, we note that
\[J_t \geqslant \tilde J_t := \frac12 \dist(x,y) + \sharp  \{ i\in \cco 1,N\ccf\, : \ x_i =y_i,\ X_{i,t} \neq Y_{i,t}\}\,.\]
Indeed, it is clear by definition that $2\tilde J_t \geqslant \dist(X_t,Y_t)$ for all $t\geqslant 0$. Since $J_0 = \tilde J_0$, suppose that $J_{T_n} \geqslant \tilde J_{T_n}$ for some $n\in \N$. Then $J_t = J_{T_n} \geqslant \tilde J_{T_n} =  \tilde J_t$ for all $t\in[T_n,T_{n+1})$. From 
\[\dist \po \tilde X_{T_{n+1}}, \tilde Y_{T_{n+1}}\pf \ \leqslant \dist\po X_{T_{n}},  Y_{T_{n}}\pf\ \leqslant \ 2  \tilde J_{T_{n}}\,,\]
we get that
\begin{eqnarray*}
\tilde J_{T_{n+1}} & =   & \tilde J_{T_n} + \1_{x_{I_n}= y_{I_n}} \1_{V_n \geqslant p_{I_n}\po \tilde X_{T_{n+1}}, \tilde Y_{T_{n+1}}\pf }\\
& \leqslant & \tilde J_{T_n} + \1_{x_{I_n}= y_{I_n}} \1_{V_n \geqslant 1- \theta \tilde J_{T_n}/N} \\
& \leqslant &  J_{T_n} + \1_{x_{I_n}= y_{I_n}} \1_{V_n \geqslant 1- \theta  J_{T_n}/N} \ =\  J_{T_{n+1}}\,,
\end{eqnarray*}
which concludes. 
\end{proof}
\subsubsection{The total variation case}

\begin{proof}[Proof of Theorem \ref{Thm_particules_TV}]
Keep all the notations of Section \ref{Section-Particules-Couplage} and, for $t\leqslant 0$, let 
\[\mathcal{A}_t = \{i\in\N\ :\, Z_{i,t}= \tilde Z_{i,t},\ x_i\neq y_i\text{ and }I_n\neq i\ \forall n\in\N\text{ such that }T_n<t \}\,.\] 
Remark that, if $i\in\mathcal A_t$, then $X_{i,t} =Y_{i,t}$, so that
\begin{multline*}
\frac12 \dist \po   X_t, Y_t \pf \ = \ \frac12 \dist(x,y)+ \sharp  \{ i\in \cco 1,N\ccf\, : \ x_i =y_i,\ X_{i,t} \neq Y_{i,t}\} \\
 \ - \sharp  \{ i\in \cco 1,N\ccf\, : \ x_i \neq y_i,\ X_{i,t} = Y_{i,t}\} \ \leqslant \  J_t - \sharp \mathcal A_t\,.
\end{multline*}

For all $i\in \N$, $\mathcal T_i = \{T_n>0\ : \ I_n = i,\ n\in\N\}$ are the jump times of a Poisson process of intensity $\lambda_*$, independent from $\mathcal T_j$ for $j\neq i$ and from $(Z_{i,t},\tilde Z_{i,t})_{t\geqslant 0}$. As a consequence, for all $i\in \cco 1,N\ccf$ with $x_i\neq y_i$, $\1_{i\in\mathcal A_t}$ is a Bernoulli r.v. with parameter $\exp(-\lambda_*t) \mathbb P ( Z_{i,t} = \tilde Z_{i,t})$, independent from $\1_{j\in\mathcal A_t}$ if $j\neq i$ and from $J_t$. In particular, 
\[\mathbb E(\sharp \mathcal A_{t_0}) \ \geqslant\ \frac12 \dist(x,y) \exp(-\lambda_*t_0) \alpha\,.\]
%
On the other hand, the generator of the Markov process $(J_t)_{t\geqslant 0}$ is
\[K f(s) \ = \   \theta s \po 1-\frac{\dist(x,y)}{2N}\pf  \po f(s+1) - f(s) \pf\,. \]
Applied to $f(s) =s$, this yields
\begin{eqnarray}\label{Eq-EspJ0}
\partial_t \mathbb E\po J_t\pf & = &   \theta   \po 1-\frac{\dist(x,y)}{2N}\pf \mathbb E\po J_t\pf \ \leqslant \   \theta     \mathbb E\po J_t\pf\,,
\end{eqnarray}
and thus 
\[\mathbb E\po \dist \po   X_{t_0}, Y_{t_0}\pf \pf \ \leqslant \ \dist(x,y) e^{ \theta t_0   } - \alpha e^{-\lambda_* t_0} \dist(x,y)\,. \]
On the other hand, for all $t\in [0,t_0]$,
\[\mathbb E\po \dist \po   X_{t}, Y_{t}\pf \pf \ \leqslant \ 2\mathbb E\po J_t\pf \ \leqslant \ \dist(x,y) e^{  \theta t_0  } \,. \]
Considering on $\mathcal P(E)$ the distance $\rho$ defined by
\[\rho(\mu_1,\mu_2) \ = \ \inf\left\{\mathbb E\po \dist(X,Y)\pf\,:\, X\sim \mu_1,\ Y\sim\mu_2 \right \}\,,\] 
we have thus obtained that, for all $t\in [0,t_0]$ and $x,y\in E$,
\[\rho \po\delta_x R_t , \delta_y R_t\pf \ \leqslant \ \po  e^{ \theta t_0} - \1_{t=t_0}\alpha e^{-\lambda_* t_0} \pf \dist(x,y)\,.\] 
Finally, for all $t\geqslant 0$ and $x,y\in E$,
\begin{eqnarray*}
 \|\delta_x R_t - \delta_y R_t\|_{TV} & \leqslant & \rho \po\delta_x R_t , \delta_y R_t\pf \\
& \leqslant & \dist(x,y) e^{  \theta t_0  } \po e^{ \theta t_0 } - \alpha e^{-\lambda_* t_0}\pf^{\lfloor t/t_0\rfloor}\\
& \leqslant & 2N\1_{x\neq y } e^{ \theta t_0  } \po e^{ \theta t_0 } - \alpha e^{-\lambda_* t_0}\pf^{\lfloor t/t_0\rfloor} \,,
\end{eqnarray*}
which, integrated with respect to any initial distribution, concludes the proof of the first statement of Theorem \ref{Thm_particules_TV}.

Now, let $m_t'$ and $h_t'$ be the respective distributions of $X_{I,t}$ and $Y_{I,t}$, where $I$ is independent from $(X_t,Y_t)_{t\geqslant 0}$. Then
\[ \|m_t' - h_t'\|_{TV} \ \leqslant \ 2\mathbb P\po   X_{I,t} \neq Y_{I,t}\pf \ = \ \frac1N \mathbb E\po \dist \po   X_{t}, Y_{t}\pf  \pf\,.\]
Considering an optimal coupling $(X_0,Y_0)$ of $m_0$ and $h_0$ and, conditionally to $(X_0,Y_0)$, an optimal coupling $(X_t,Y_t)$ of $\delta_{X_0}R_t$ and $\delta_{Y_0} R_t$,
\begin{eqnarray*}
\mathbb E\po \dist \po   X_{t}, Y_{t}\pf  \pf & \leqslant & \mathbb E\po \dist \po   X_{0}, Y_{0}\pf  \pf e^{ \theta t_0  }  \po e^{ \theta t_0 } - \alpha e^{-\lambda_* t_0}\pf^{\lfloor t/t_0\rfloor}\\
& =& 2 N \mathbb P\po   X_{I,0}\neq  Y_{I,0}\pf e^{ \theta t_0  }  \po e^{ \theta t_0 } - \alpha e^{-\lambda_* t_0}\pf^{\lfloor t/t_0\rfloor}\\
&= & N\|m_0' - h_0'\|_{TV} e^{ \theta t_0  }  \po e^{ \theta t_0 } - \alpha e^{-\lambda_* t_0}\pf^{\lfloor t/t_0\rfloor}\,,
\end{eqnarray*}
which concludes.
\end{proof}

\subsubsection{The general $\V$ case}

In all this section,  Assumptions \ref{Hyp-TauxBorne-Particules} and \ref{Hyp-Particules-Gene} are enforced. 
 For $i\in\cco 1,N\ccf$ and $x\in E$, set 
\begin{eqnarray*}
\W_i(x) & = & \V_i(x_i) + \frac{2\eta}{N(1-\eta)}\V(x)\,.
\end{eqnarray*}

\begin{lmm}
For all $t\geqslant 0$ and $x\in E$, and $i\in \cco 1,N\ccf$,
\begin{eqnarray}\label{Eq-Lyap-Particules-V}
R_{t} \V(x) & \leqslant & e^{- \rho(1-\eta) t }\V(x) + \po 1 - e^{- \rho(1-\eta) t } \pf \frac{N M}{1-\eta}\\
R_{t} \W_i(x) & \leqslant & e^{- \frac12\rho(1-\eta) t }\W_i(x) +      \frac{2(1+\eta) M }{(1-\eta)^2} \po 1 - e^{- \frac12\rho(1-\eta) t } \pf\label{Eq-Lyap-Particules-Vi} \,.
\end{eqnarray}
\end{lmm}
\begin{proof}
Summing \eqref{Eq-Hyp-Particule-Lyap} over $i\in \cco 1,N\ccf$ reads
\begin{eqnarray*}
R_{t} \V(x) & \leqslant & e^{- \rho t }\V(x) + \rho \int_0^t e^{ \rho(s-t) } \po  \eta  R_s\V(x) + N M\pf \dd s 
\end{eqnarray*}
for all $t\geqslant 0$ and $x\in E$. The Gronwall's Lemma then yields \eqref{Eq-Lyap-Particules-V} which, reintegrated in \eqref{Eq-Hyp-Particule-Lyap}, gives 
\begin{eqnarray*}
R_{t} \W_i(x) & \leqslant & e^{- \rho t }\V_i(x_i) + e^{- \rho(1-\eta) t }\frac{2\eta \V(x)}{N(1-\eta)} +   \po 1 - e^{- \rho(1-\eta) t } \pf \frac{2\eta M}{(1-\eta)^2}\\
& & {+  \rho \int_0^t e^{ \rho(s-t) }  \po   \frac{\eta}{N} e^{- \rho(1-\eta) s }\V(x) + \po 1 - e^{- \rho(1-\eta) s } \pf \frac{\eta M}{1-\eta} +M\pf  \dd s }\\
& = &  e^{- \rho t }\V_i(x_i) + \po e^{- \rho(1-\eta) t } + \frac{1-\eta}{2\eta} \po e^{- \rho(1-\eta) t }  - e^{- \rho t }  \pf \pf \frac{2\eta \V(x)}{N(1-\eta)} \\
& & {+   \po 1 - e^{- \rho(1-\eta) t } \pf \frac{2\eta M}{(1-\eta)^2} +  \frac{\rho M e^{-\rho t}  }{1-\eta}\int_0^t e^{\rho s} \po     1  - e^{- \rho(1-\eta) s }   \eta  \pf \dd s }\\
&\leqslant & a(t) \W_i(x) +b(t)
\end{eqnarray*}
with
\begin{eqnarray*}
a(t) &=&  e^{- \rho(1-\eta) t } + \frac{1-\eta}{2\eta} \po e^{- \rho(1-\eta) t }  - e^{- \rho t }  \pf\\
b(t)&=& {\po \frac{2\eta M}{(1-\eta)^2} + \frac{M}{1-\eta}\pf \po 1 - e^{-\rho(1-\eta)t}\pf =  \frac{(1+\eta) M}{(1-\eta)^2}  \po 1 - e^{-\rho(1-\eta)t}\pf\,.}
\end{eqnarray*}
Applying $R_s$ for some $s\geqslant0$ to this inequality and using the semi-group property, we have thus obtained for all $s,t\geqslant 0$
\[R_{t+s} \W_i(x) \ \leqslant \ a(t)R_{s} \W_i(x) +b(t)\,.\]
Hence, for all $n\in\N$ and $t\geqslant 0$,
\begin{eqnarray*}
R_t \W_i(x) & \leqslant & a(t/n) R_{t-t/n}\W_i(x) + b(t/n) \\ 
& \leqslant & \dots \ \leqslant \ \po a(t/n)\pf^n \W_i(x) +b(t/n) \frac{1-\po a(t/n)\pf^n}{1-a(t/n)}\,.
\end{eqnarray*}
Letting $n$ go to infinity yields
\[R_t \W_i(x)\leqslant \ e^{a'(0)t} \W_i(x) - \frac{b'(0)}{a'(0)} \po 1 - e^{a'(0)t}\pf  \,,\]
which is \eqref{Eq-Lyap-Particules-Vi}.
\end{proof}

Similarly to \cite{HairerMattingly2008} or Section \ref{Section-NL-Gene}, we now consider  some parameter $\beta>0$ and, for all $x,y\in E$ and $i\in\cco 1,N\ccf$,
\[d_{i,\beta}(x,y) \ = \ \1_{x_i\neq y_i} \po 2\beta +  \W_i(x) + \W_i(y)\pf\,.\]

\begin{lmm}\label{Lem-Particules-XneqY}
For  all $i\in\cco 1,N\ccf$ and all $x,y\in E$ with $x_i\neq y_i$, let   $(X_t,Y_t)_{t\geqslant 0}$ be the coupling of $(\delta_x R_t)_{t\geqslant0}$ and $(\delta_y R_t)_{t\geqslant0}$ defined in Section \ref{Section-Particules-Couplage}. Then
\begin{eqnarray*}
\mathbb E \po d_{i,\beta} \po X_{t_0},Y_{t_0}\pf\pf & \leqslant & \kappa  d_{i,\beta}(x,y) 
\end{eqnarray*}
with 
\begin{eqnarray}\label{Eq-Def-Beta-Particules} 
\beta & =&    \frac{4(1+\eta) M e^{\lambda_* t_0} }{\alpha (1-\eta)^2} \po 1 - e^{- \frac12\rho(1-\eta) t_0 } \pf \\
 \kappa & = & \frac{1+e^{- \frac12\rho(1-\eta) t_0 }}2  \vee \po  1-\frac12 \alpha e^{-\lambda_* t_0}\pf\,.
\end{eqnarray}
\end{lmm}
\begin{proof}
This is again essentially the proof of \cite[Theorem 3.1]{HairerMattingly2008}. We bound
\begin{eqnarray*}
\mathbb E \po   d_{i,\beta} \po X_{t_0},Y_{t_0}\pf \pf & \leqslant & 2\beta \mathbb P\po X_{i,t_0}\neq Y_{i,t_0}\pf +   R_{t_0} \W_i(x) +   R_{t_0} \W_i(y)\,.
\end{eqnarray*}
From \eqref{Eq-Lyap-Particules-Vi},
\begin{eqnarray*}
R_{t_0} \W_i(x)  + R_{t_0} \W_i(y) & \leqslant & e^{- \frac12\rho(1-\eta) t_0 } \po \W_i(x) + \W_i(y)\pf  +      \frac{4(1+\eta) M }{(1-\eta)^2} \po 1 - e^{- \frac12\rho(1-\eta) t_0 } \pf\,.
\end{eqnarray*}
Two cases are distinguished. First, if $\V_i(x_i)+\V_i(y_i) \leqslant 16(1+\eta )M/(1-\eta)^2$,  as in the previous section, from \eqref{Eq-Particule-Gene-Doeblin},
\begin{eqnarray*}
\mathbb P\po X_{i,t_0}=Y_{i,t_0}\pf&  \geqslant & \mathbb P \po Z_{i,t_0} = \tilde Z_{i,t_0} \text{ and }I_n\neq i\ \forall n\text{ such that }T_n<t \pf\\
& \geqslant & \alpha  e^{-\lambda_* t_0}\,,
\end{eqnarray*}
so that 
\begin{eqnarray*}
\mathbb E \po   d_{i,\beta } \po X_{t_0},Y_{t_0}\pf \pf & \leqslant & 2\beta \po 1-\alpha  e^{-\lambda_* t_0}\pf +  e^{- \frac12\rho(1-\eta) t_0 } \po \W_i(x) + \W_i(y)\pf \\
& & \ +      \frac{4(1+\eta) M }{(1-\eta)^2} \po 1 - e^{- \frac12\rho(1-\eta) t_0 } \pf\\
& \leqslant & 2\beta \po 1-\frac12 \alpha  e^{-\lambda_* t_0}\pf + e^{- \frac12\rho(1-\eta) t_0 } \po \W_i(x) + \W_i(y)\pf 
\end{eqnarray*}
where we used the definition \eqref{Eq-Def-Beta-Particules}  of $\beta$. Second, if $\V_i(x_i)+\V_i(y_i) \geqslant 16(1+\eta )M/(1-\eta)^2 $, since $\W_i \geqslant \V_i$, then
\begin{eqnarray*}
\mathbb E \po d_{i,\beta } (X_{i,t_0},Y_{i,t_0})\pf & \leqslant & 2\beta   + e^{- \frac12\rho(1-\eta) t_0 } \po \W_i(x) + \W_i(y)\pf  +      \frac{4(1+\eta) M }{(1-\eta)^2} \po 1 - e^{- \frac12\rho(1-\eta) t_0 } \pf\\ 
& \leqslant & 2\beta   + \frac{1+e^{- \frac12\rho(1-\eta) t_0 }}2 \po \W_i(x) + \W_i(y)\pf \\
& & \  -\      \frac{4(1+\eta) M }{(1-\eta)^2} \po 1 - e^{- \frac12\rho(1-\eta) t_0 } \pf\\
& = & 2\beta \po 1-\frac12 \alpha  e^{-\lambda_* t_0}\pf + \frac{1+e^{- \frac12\rho(1-\eta) t_0 }}2 \po \W_i(x) + \W_i(y)\pf
\end{eqnarray*}
where we used again the definition \eqref{Eq-Def-Beta-Particules}  of $\beta$.
\end{proof}

\begin{lmm}\label{Lem-Particules-X=Y}
For  all $i\in\cco 1,N\ccf$, $\beta>0$, $t\geqslant 0$ and  $x,y\in E$ with $x_i = y_i$, if  $(X_t,Y_t)_{t\geqslant 0}$ is the coupling of $(\delta_x R_t)_{t\geqslant0}$ and $(\delta_y R_t)_{t\geqslant0}$ defined in Section \ref{Section-Particules-Couplage}, then
\begin{eqnarray*}
\mathbb E \po d_{i,\beta} \po X_{t},Y_{t}\pf\pf 
& \leqslant  &  \frac{\dist(x,y)}{2N}\po e^{\theta t} - 1\pf \po 2\beta + c_*(t) \po \W_i(x) + \W_i(y)\pf  \pf \,,
\end{eqnarray*}
with $c_*(t) = (\gamma_*+1) \exp(\po \rho_* +\lambda_*(\gamma_*-1)\pf t)$.
\end{lmm}
\begin{proof}
Let $t_1\geqslant 0$ and $i\in\cco 1,N\ccf$ be such that $x_i = y_i$. For all $j\in \cco 1,N\ccf$, let $\tau_j= \inf\{T_{n}\geqslant 0\ :\ n\in\N,\ I_n=j,\ V_n \geqslant 1- \theta J_{T_n}/N\}$. By construction and Lemma \ref{Prop-Coupl-Particule-Ok}, $X_{i,t} = Y_{i,t}$ for all $t\leqslant \tau_i$, so that 
\begin{eqnarray}\label{Eq-Preuve-X=Y-particule_1}
\mathbb E \po d_{i,\beta } \po X_{t_1},Y_{t_1}\pf\pf & \leqslant & \mathbb E \po \1_{\tau_i \leqslant t_1} \po 2\beta +  \W_i(X_{t_1}) + \W_i(Y_{t_1})\pf \pf\,.
\end{eqnarray}
Let $\Gamma = \po  S_k,V_k,I_k\pf _{k\in\mathbb N}$, $K=\sharp \{n\in\N_*\ : \ T_{n} \leqslant t_1\}$ and for all $j\in\cco 1,N\ccf$,  $\mathcal A_j= \{n\in\N_*\ : \ T_{n} \leqslant t_1,\ I_n=j\}$ and $K_j=\sharp \mathcal A_j$. 
Note that  $K$ and $\{\tau_j,K_{j}\}_{j\in\cco 1,N\ccf}$ are deterministic functions of $\Gamma$, and in particular are independent from $\po  U_{i,k}\pf _{i\in\cco 1,N\ccf,k\in\mathbb N}$. Hence, for all $n \in \N$ and $j\in\cco 1,N\ccf$, using \eqref{Eq-Hyp-Particule-Gene-1} and \eqref{Eq-Hyp-Particule-Gene-4},
\begin{eqnarray*}
\mathbb E \po   \V_j(X_{j,T_{n}})  \ |\ \Gamma \pf & \leqslant &  \po \1_{I_n \neq j} + \gamma_*\1_{I_n=j}\pf    \mathbb E \po   \V_j(\tilde X_{j,T_{n}})  \ |\ \Gamma \pf \\
& = & \po \1_{I_n \neq j} + \gamma_*\1_{I_n=j}\pf   \mathbb E \po   P_{T_{n}-T_{n-1}}^j\V_j(X_{j,T_{n-1}})  \ |\ \Gamma \pf \\
& \leqslant &   \po \1_{I_n \neq j} + \gamma_*\1_{I_n=j}\pf  e^{\rho_*(T_{n}-T_{n-1})}  \mathbb E \po    \V_j(X_{j,T_{n-1}})  \ |\ \Gamma \pf\,.
\end{eqnarray*}
Similarly, and then by direct induction,
\begin{eqnarray*}
\mathbb E \po   \V_j(X_{j,t_1})  \ |\ \Gamma \pf & \leqslant &   e^{\rho_*(t_1-T_{K})} \mathbb E \po    \V_i(X_{j,T_{K}})  \ |\ \Gamma \pf   \\
& \leqslant & e^{\rho_* t_1} \gamma_*^{K_j}\V_j(x_j) \,.
\end{eqnarray*}
We have thus obtained for all $j\in\cco 1,N\ccf$
\begin{eqnarray}\label{Eq-Preuve-X=Y-particule_2}
\mathbb E \po \1_{X_{i,t_1}\neq Y_{i,t_1}} \V_j(X_{j,t_1})  \pf & \leqslant & e^{\rho_* t_1} \V_j(x_j)  \mathbb E \po \1_{\tau_i \leqslant t_1} \gamma_*^{K_j} \pf\,.
\end{eqnarray}
For $j\in \cco 1,N\ccf$, set $\Gamma_j = \{  (S_k,V_k,I_k)\ : \ I_k =  j,\ k\in\mathbb N\}$. Remark that $(J_t)_{t\geqslant 0}$ is a deterministic function of $\{\Gamma_j\ :\ j\in\cco 1,N\ccf,\ x_j= y_j\}$ and that, by Poisson thinning, $\Gamma_j$ (hence $K_j$) is independent from $\Gamma_k$ for all $k\neq j$. In particular, if $j\in\cco 1,N\ccf$ is such that $x_j \neq  y_j$,
\begin{eqnarray}\label{Eq-Preuve-X=Y-particule_3}
\mathbb E \po \1_{\tau_i \leqslant t_1} \gamma_*^{K_j} \pf & = & \mathbb P \po \tau_i \leqslant t_1\pf \mathbb E \po   \gamma_*^{K_j} \pf \ = \ \mathbb P \po \tau_i \leqslant t_1\pf  e^{\lambda_*(\gamma_*-1) t_1} \,.
\end{eqnarray}
For $j\in \cco 1,N\ccf$ with $x_j\neq y_j$ we define the process $(J_{j,t})_{t\geqslant 0}$ as follows. Set $J_{j,0}=0$ and, for all $n\in\N$ and   $t\in(T_n,T_{n+1})$, set $J_{j,t} = J_{j,T_n}$ and
\[J_{j,T_{n+1}} = J_{j,T_n} + \1_{I_n\neq j } \1_{V_n \geqslant 1- \theta J_{j,T_n}/(N\lambda_*)}\,.\] 
In other words, $(J_{j,t})_{t\geqslant 0}$ is like $(J_t)_{t\geqslant 0}$ but ignoring the jumps of the $j^{th}$ particle. By thinning of Poisson processes, $(J_{j,t})_{t\geqslant 0}$  is independent from $\Gamma_j$. Moreover $J_{j,t} \leqslant J_t$ for all $t$ and, conditionnally to $\tau_j > t$, $J_{j,t} = J_t$. As a consequence, if $i\neq j$, we bound
\begin{eqnarray*}
\mathbb P \po \tau_i \leqslant  t_1,\ \tau_j>t_1\ |\ K_j\pf &= & \mathbb P \po \tau_j>t_1,\ \exists n\in \mathcal A_i\text{ s.t. } V_n \geqslant 1- \theta J_{j,T_n}/(N\lambda_*)\ |\ K_j \pf  \\
 & \leqslant & \mathbb P \po \exists n\in \mathcal A_i\text{ s.t. } V_n \geqslant 1- \theta J_{j,T_n}/(N\lambda_*)\ |\ K_j \pf  \\
 & = & \mathbb P \po \exists n\in \mathcal A_i\text{ s.t. } V_n \geqslant 1- \theta J_{j,T_n}/(N\lambda_*) \pf \\
 & \leqslant & \mathbb P \po \exists n\in \mathcal A_i\text{ s.t. } V_n \geqslant 1- \theta J_{T_n}/(N\lambda_*) \pf \\
 & = &  \mathbb P \po \tau_i \leqslant t_1\pf\,,
\end{eqnarray*}
and thus, for all $j\neq i$ with $x_j=y_j$,
\begin{eqnarray}
\notag \mathbb E \po \1_{\tau_i \leqslant t_1} \gamma_*^{K_j} \pf & \leqslant & \mathbb E \po \1_{\tau_j \leqslant t_1} \gamma_*^{K_j} \pf  + \mathbb E \po \1_{\tau_i \leqslant t_1} \1_{\tau_j > t_1}  \gamma_*^{K_j} \pf\\
\notag & \leqslant & \mathbb E \po \1_{\tau_j \leqslant t_1} \gamma_*^{K_j} \pf + \mathbb P \po \tau_i \leqslant t_1\pf \mathbb E \po \gamma_*^{K_j} \pf   \\
& = & \mathbb E \po \1_{\tau_j \leqslant t_1} \gamma_*^{K_j} \pf + \mathbb P \po \tau_i \leqslant t_1\pf  e^{\lambda_*(\gamma_*-1) t_1} \label{Eq-Preuve-X=Y-particule_4} \,. 
\end{eqnarray}
Integrating \eqref{Eq-Preuve-X=Y-particule_3} and \eqref{Eq-Preuve-X=Y-particule_4} in \eqref{Eq-Preuve-X=Y-particule_2} and then in \eqref{Eq-Preuve-X=Y-particule_1} by definition of $\W_i$, we have obtained so far that
\begin{multline}\label{Eq-Preuve-X=Y-particule_5}
  \mathbb E \po d_{i,\beta } \po X_{t_1},Y_{t_1}\pf\pf \ \leqslant  \  \mathbb P \po \tau_i \leqslant t_1\pf\po 2\beta + \frac{2\eta e^{\po \rho_* +\lambda_*(\gamma_*-1)\pf t_1} }{N(1-\eta)}  \po \V(x)+\V(y)\pf\pf \\
 \ +  2 e^{  \rho_*   t_1} \V_i(x_i) \mathbb E \po \1_{\tau_i \leqslant t_1} \gamma_*^{K_i} \pf +  \frac{4\eta e^{  \rho_*   t_1}}{N(1-\eta)}   \sum_{j=1}^N \1_{x_j=y_j} \V_j(x_j) \mathbb E \po \1_{\tau_j \leqslant t_1} \gamma_*^{K_j} \pf \,.
\end{multline}
%
Using that $(J_{i,t})_{t\geqslant 0}$, $K_i$ and $(V_n)_{n\in\mathcal A_i}$ are independent r.v. and that, conditionally to $K_i$, $\{T_n\ : \ n\in\mathcal A_i\}$ are distributed like $K_i$ independent r.v. uniformly distributed over $[0,t_1]$,
\begin{eqnarray*}
\mathbb P\po \tau_i > t_1\ |\ K_i\pf & = & \mathbb P \po  V_n \leqslant 1- \theta J_{i,T_n}/(N\lambda_*),\, \forall n\in\mathcal A_i \ |\ K_i\pf \\
& = &   \po 1 - \frac1{t_1}\int_0^{t_1} 1 \wedge \mathbb E\po  \theta  J_{i,w} /(N\lambda_*)\pf \dd w \pf^{K_i} \\
& \geqslant &  \po 1 - \frac1{t_1}\int_0^{t_1} 1 \wedge \mathbb E\po  \theta  J_{w} /(N\lambda_*)\pf \dd w \pf^{K_i} \ :=\ s^{K_i}\,.
\end{eqnarray*}
With the same computations as in the proof of Lemma \ref{Lem_NL_x=y}, this leads to 
\begin{eqnarray*}
\mathbb E \po \1_{\tau_i \leqslant t_1} \gamma_*^{K_i} \pf
& \leqslant &  \lambda_*t_1\gamma_* (1-s)        e^{  \lambda_*(\gamma_* -1) t_1} \\
& \leqslant &  \gamma_*        e^{  \lambda_*(\gamma_* -1) t_1} \frac{\dist(x,y)}{2N} \int_0^{t_1} \theta e^{\theta u}\dd u \,, 
\end{eqnarray*}
where we used that, from \eqref{Eq-EspJ0}, $\mathbb E(J_{t} )\leqslant \exp(\theta  t) \dist(x,y)/2$ for all $t\geqslant 0$. The case of $j\neq i $ with $x_j=y_j$ is identical and, similarly (i.e. applying this bound with $\gamma_*=1$), 
\begin{eqnarray*}
\mathbb P\po \tau_i \leqslant t_1 \pf & \leqslant & \frac{\dist(x,y)}{2N}\po e^{\theta t_1} - 1\pf\,.
\end{eqnarray*}
Together with \eqref{Eq-Preuve-X=Y-particule_5}, this yields
\begin{eqnarray*}
  \mathbb E \po d_{i,\beta } \po X_{t_1},Y_{t_1}\pf\pf & \leqslant  &  \frac{\dist(x,y)}{2N}\po e^{\theta t_1} - 1\pf \po 2\beta + \po \gamma_*+1\pf e^{\po \rho_* +\lambda_*(\gamma_*-1)\pf t_1} \po \W_i(x) + \W_i(y)\pf  \pf \,.
\end{eqnarray*}
\end{proof}


\begin{proof}[Proof of Theorem \ref{Thm_particules_gene}]

Let $\beta$ and $\kappa$ be given by Lemma \ref{Lem-Particules-XneqY}, and consider on $\mathcal P(E)$ the distance $\bar \rho_\beta$ given for all $\mu,\nu\in\mathcal P(E)$ by
\begin{eqnarray*}
 \bar \rho_\beta(\mu,\nu) & = & \inf\left\{\mathbb E\po \sum_{i=1}^N d_{i,\beta}(X_i,Y_i)\pf\,:\,X\sim\mu,\ Y\sim \nu \right \}\,.
\end{eqnarray*}
Lemmas \ref{Lem-Particules-XneqY} and \ref{Lem-Particules-X=Y} means that for all $x,y\in E$,
\begin{eqnarray*}
 \bar \rho_\beta(\delta_x R_{t_0}, \delta_y R_{t_0})& \leqslant & \kappa \bar \rho_\beta(\delta_x  , \delta_y  ) + \frac{\dist(x,y)}{2}\po e^{\theta t} - 1\pf \po 2\beta + \frac{c_*(t)}N \sum_{i=1}^N \po \W_i(x) + \W_i(y)\pf  \pf\\
 &\leqslant & \kappa \bar \rho_\beta(\delta_x  , \delta_y  ) + \frac{\dist(x,y)}{2}\po e^{\theta t} - 1\pf c_*(t) \po 2\beta + \frac{1+\eta}{N(1-\eta) }\po \V(x) + \V(y)\pf\pf\,.
\end{eqnarray*}
Since $\eta\geqslant 1/2$ by assumption,
\begin{eqnarray*}
\frac{\dist(x,y)}{2 N(1-\eta) }\po \V(x) + \V(y)\pf   &  = & \sum_{i=1}^N \1_{x_i \neq y_i} \frac{\V(x) + \V(y)}{2 N(1-\eta) }\\
 & \leqslant & 2 \sum_{i=1}^N \1_{x_i \neq y_i} \po \W_i(x) + \W_i(y)\pf\,,
\end{eqnarray*}
and we get that 
\begin{eqnarray*}
 \bar \rho_\beta(\delta_x R_{t_0}, \delta_y R_{t_0})& \leqslant &  \po \kappa + \po e^{\theta t} - 1\pf c_*(t_0) (1+\eta) \pf  \bar \rho_\beta(\delta_x  , \delta_y  ) \ = \ \tilde  \kappa  \bar \rho_\beta(\delta_x  , \delta_y  )\,.
\end{eqnarray*}
Then for all initial distributions $\mu,\nu\in\mathcal P(E)$, conditioning with respect to the initial values
\begin{eqnarray*}
\bar\rho_\beta(\mu R_{t_0},\nu R_{t_0}) & = & \inf\left\{\mathbb E\po \bar\rho_\beta(\delta_{X_0} R_{t_0},\delta_{Y_0} R_{t_0})\pf\,:\,X_0\sim\mu,\ Y_0\sim \nu \right \}\ \leqslant \ \tilde \kappa \bar\rho_\beta(\mu ,\nu  )\,.
\end{eqnarray*}
As in the proof of Theorem \ref{Thm_non-lin}, the assumption that $\tilde \kappa <1$ implies that $R_{t_0}$ admits a unique fixed point $\mu_\infty$ in $\mathcal P_\V(E)$. Moreover, by the semi-group property, for all $s\geqslant 0$,
\[\mu_\infty R_s R_{t_0} = \mu_\infty R_{t_0} R_s =  \mu_\infty R_s \,,\]
so that $\mu_\infty R_s$ is a fixed point of $R_{t_0}$ and thus by uniqueness it is equal to $\mu_\infty$. 
 Integrating \eqref{Eq-Lyap-Particules-V} with respect to $\mu_\infty$ and 
letting $t$ go to infinity yields $\mu_\infty(\V)  \leqslant NM/(1-\eta)$, and thus for all $\nu \in\mathcal P_\V(E)$,
\begin{eqnarray*}
\bar \rho_\beta(\nu,\mu_\infty) & \leqslant &  2\beta N + \frac{1+\eta}{1-\eta }\po \nu(\V) + \mu_\infty(\V)\pf \\
& \leqslant &  \po 2\beta  + M \frac{1+\eta}{(1-\eta)^2 }\pf N+ \frac{1+\eta}{1-\eta }  \nu(\V) \,.
\end{eqnarray*}
Finally, for all $\nu\in\mathcal P_\V(E)$ and $t\geqslant 0$
\begin{eqnarray*}
\bar \rho_\beta(\nu R_t,\mu_\infty) & \leqslant & \tilde \kappa^{\lfloor t/t_0\rfloor} \bar \rho_\beta(\nu R_{t-\lfloor t/t_0\rfloor t_0},\mu_\infty)\\
& \leqslant & \tilde \kappa^{\lfloor t/t_0\rfloor} \po \po 2\beta  + M \frac{1+\eta}{(1-\eta)^2 }\pf N+ \frac{1+\eta}{1-\eta }  \nu R_{t-\lfloor t/t_0\rfloor t_0}(\V) \pf\\
& \leqslant & \tilde \kappa^{\lfloor t/t_0\rfloor} \po 2N \po \beta  + M \frac{1+\eta}{(1-\eta)^2 }\pf  + \frac{1+\eta}{1-\eta }  \nu (\V) \pf\,,
\end{eqnarray*}
where we used \eqref{Eq-Lyap-Particules-V}. The conclusion then follows from the fact that $\bar \rho_\beta \geqslant \bar d_\V$ for all $\beta\geqslant 0$.

\end{proof}

\section{Examples}\label{Sec-Exemples}

This section provides illustrations  of our main results. For the sake of clarity and since the approach would be the same in more general or sophisticated applications, the models are chosen to be simple in order to highlight the core arguments.

As has already been mentionned, the neuron network model of \cite{CanizoYoldas} provides a first example where Theorem \ref{Thm-TVNonLin} applies.

\subsection{Mean-field run-\&-tumble process}\label{SectionRTP}

Run-\&-tumble processes model, among other things, the motion of some bacteria \cite{Patlak,ErbanOthmer} (see also \cite{Fetique,Calvez,FontbonaGuerinMalrieu2016} for more recent details and references). The rate at which a bacterium tumbles depends on the concentration of given chemo-attractants in the medium. For a large population of bacterium, mean-field interaction is a natural extension of these dynamics.

We consider the non-linear integro-differential equation on $E=\R\times\{-1,1\}$ given by
\begin{eqnarray}\label{Eq-RTP}
\partial_t m_t(x,y) + y \partial_x m_t(x,y) & = & \lambda_{m_t}(x,-y)m_t(x,-y) -  \lambda_{m_t}(x,y)m_t(x,y)\,,
\end{eqnarray}
with ${\lambda_\nu(x,y)} = r(y(x - \theta x_\nu))$ where $x_\nu = \int_{E} x \nu(\dd x,\dd y)$ is  the barycentre of $\nu$, $\theta\in(0,1)$ and $r:\R\rightarrow\R_+$ is a Lipschitz function such that 
\begin{eqnarray}\label{Eq-CondRTP}
  \lim_{s\rightarrow -\infty} r(s) \  :=a & <&  b\  :=\  \lim_{s\rightarrow {+\infty}} r(s).
\end{eqnarray}
 The linear case $\theta = 0$ corresponds to the run-\&-tumble process attracted to a neighborhood of origin studied in \cite{FontbonaGuerinMalrieu2016} (except that we don't assume that $r$ is non-decreasing). Indeed, when  $|x|$ is large, the jump rate $r(xy)$ is larger when $xy>0$, namely when the process is drifting away from a given compact, than when $xy<0$, i.e. when the process is going toward the compact. Similarly, the case $\theta = 1$ would correspond to a process attracted toward the barycenter of its law. For $\theta\in(0,1)$, the process is attracted to an average of the origin and of the barycenter of its law. Note that, if $s\mapsto r(s) -r(0)$ were anti-symmetric, then the linear process would admit a symmetric equilibrium, which would then be  an equilibrium for the non-linear equation since its mean is zero. In the following, this symmetry is not assumed.

A similar multi-dimensional process could be considered with the same assumptions as in \cite{Fetique} to  ensure the existence of a Lyapunov function. The arguments below could then be straightforwardly adapted.

\begin{prpstn}\label{Prop-RTP}
Assume  \eqref{Eq-CondRTP} holds for some $0<a<b$ and $c:= \inf_{s\in \R} r(s)>0$. Then there exist $C,\theta_*>0$ and $\kappa\in(0,1)$ such that, if $\theta\leqslant \theta_*$, then \eqref{Eq-RTP} admits a unique equilibrium $\mu\in\mathcal P(E)$ such that for all solution $t\mapsto m_t$ of \eqref{Eq-RTP} and all $t\geqslant 0$,
\[\|m_t - \mu\|_\W \ \leqslant \ C \kappa^t m_0(\W)\,,\]
with $\W(x,y) = \exp((b-a)|x|/4 )$.
\end{prpstn}
\begin{proof}

Let us establish that Assumption \ref{Hyp-MeanField-Exist} holds for \eqref{Eq-RTP}. To recover the notations of Section \ref{Section-def-non-lineaire}, consider the generator $L$ on $E$ given by
\begin{eqnarray*}
Lf(x,y) & = & y \partial_x f(x,y) + c \po f(x,-y) - f(x,y)\pf\,,
\end{eqnarray*}
the jump rate ${\tilde{\lambda}_\nu}(x,y) = r(y(x - \theta  x_\nu)) - c$ and the kernel $Q_{\nu}(x,y) = Q(x,y) = \delta_{(x,-y)}$. In particular, for $z=(x,y)\in E$,
\begin{eqnarray}\label{Eq-RTP-Lip}
\notag \| {\tilde{\lambda}_\nu}(z) (Q(z)-\delta_z) - {\tilde{\lambda}_\mu}(z) (Q(z)-\delta_z) \|_{TV} & = & |\lambda_{\nu}(z)-\lambda_{\mu}(z)|\\
&\leqslant & \theta \| r'\|_\infty|x_\nu - x_\mu|\,, 
\end{eqnarray}
and 
\[|x_\nu - x_\mu| \ \leqslant\ \int |x| |\mu-\nu| \ \leqslant\ k^{-1}\|\mu-\nu\|_\V \]
for all $\V$ and $k>0$ such that $\V(x,y) \geqslant k|x|$ for all $(x,y)\in E$. 

The Markov process with generator $L$ is the so-called integrated telegraph process and it is clear, either by simple controllability argument as in \cite{MonmarcheRTP}, more precise coupling estimates as \cite{FontbonaGuerinMalrieu2016} or just explicit computations of the density transition, that for any compact set $\mathcal K$ of $\R$ there exist $t_0,\alpha>0$ such that for all $x,x'\in \mathcal K$ and $y,y'\in\{-1,1 \}$, the Doeblin condition 
\begin{eqnarray}\label{Eq-RTP-Doeblin}
\| \delta_{(x,y)} P_{t_0} - \delta_{(x',y')} P_{t_0}\|_{TV} & \leqslant &  2(1-\alpha)
\end{eqnarray}
holds. Note that $t_0$ and $\alpha$ only depends on $\mathcal K$ and $c$.

 It only remains to construct a Lyapunov function $\V$  that satisfies \eqref{Eq-Hyp-NonLin-Lyap}, \eqref{Eq-Hyp-NonLin-Gene-1} and \eqref{Eq-Hyp-NonLin-Gene-4} and such that $\V(x,y) \geqslant k|x|$ for all $(x,y)\in E$ for some $k>0$. Given a measurable function $t\mapsto \mu_t$ from $\R_+$ to $\mathcal P(E)$, recall that we denote $(R_{s,t}^\mu)_{t\geqslant s\geqslant 0}$ the semi-group associated to the inhomogeneous process defined in Section \ref{Section-def-non-lineaire} with the generator
\[L_t f(x,y) \ = \ y \partial_x f(x,y) + \lambda_{\mu_t}(x,y)\po f(x,-y) - f(x,y)\pf\,.\]

%
Let $\varphi:\R\rightarrow [0,1]$ be defined by $\varphi(s) = 1$ for $s\geqslant 1$, $\varphi(s)= 0$ for $s\leqslant -1$ and $\varphi(s) = (1+\sin(\pi s/2))/2$ for $s\in[-1,1]$, and let $h:\R\rightarrow \R_+$ be defined by $h(s)=|s|$ for ${|s|\geqslant 1}$ and $h(s) = (s^2 + 1 )/2$ for $s\in[-1,1]$. Set
\begin{eqnarray}\label{Eq-VRTP}
\V(x,y) & = & e^{\frac{b-a}{4} h(x)}\po  \frac{5a+3b}{2(b-a)} + \varphi(yx) \pf\,.
\end{eqnarray}
As a smooth function in $x$, it belongs to the domain of the generalized generator of $L_{t}$ for all $t\geqslant 0$, see \cite{Davis,DurmusGuillinMonmarche2018}.
For $(x,y)\in E$ with $|x|\geqslant 1$,
\[L_t \V(x,y) = e^{\frac{b-a}{4} |x|}\co sign(xy)\frac{b-a}{4} \po  \frac{5a+3b}{2(b-a)}  + \1_{yx>0} \pf + r\po y(x - \theta  x_{\mu_t}\pf \po \1_{xy<0} - \1_{xy>0}\pf \cf\,.\]
Take $R_0\geqslant 1$  large enough so that  $r( z) >(7b+a)/8$ and $r(-z) < (b+7a)/8$ for all $z\geqslant R_0$. Then, for all $t\geqslant 0$ and $(x,y)\in E$ with $|x| \geqslant R_0 + \theta |x_{\mu_t}|$, if $xy>0$ then
\[r\po y(x - \theta  x_{\mu_t})\pf   >(7b+a)/8\,,\]
so that
\begin{eqnarray*}
L_t \V(x,y) & \leqslant & e^{\frac{b-a}{4} |x|}\co   \frac{5a+3b}{8}  + \frac{b-a}{4}   - \frac{7b+a}{8} \cf \ \leqslant -\frac{b-a}{4}e^{\frac{b-a}{4} |x|}\, 
\end{eqnarray*}
while if $xy<0$ then
\begin{eqnarray*}
L_t \V(x,y) & \leqslant & e^{\frac{b-a}{4} |x|}\co   -\frac{5a+3b}{8}    + \frac{b+7a}{8} \cf \ \leqslant -\frac{b-a}{4} e^{\frac{b-a}{4} |x|}\,. 
\end{eqnarray*}
Hence, for all $x,y\in E$ with $|x| \geqslant R_0 + \theta |x_{\mu_t}|$ and all $t\geqslant 0$,
\begin{eqnarray*}
L_t \V(x,y) & \leqslant &  -\rho\V(x,y)
\end{eqnarray*}
with $\rho = (b-a)^2/[10(b+a)]$. Besides, 
for all $x,y\in E$ with $ |x| \leqslant R_0 + \theta |x_{\mu_t}|$, a rough bound is
\begin{eqnarray*}
L_t \V(x,y)  + \rho \V(x,y)& \leqslant &  e^{\frac{b-a}{4}h(x)} (a+2b+\pi/4)\\
& \leqslant &  e^{\frac{b-a}{4}(1+R_0+\theta |x_{\mu_t}|)} (a+2b+\pi/4)\\
& \leqslant &  e^{\frac{b-a}{4}(1+R_0)} \po 1 + \theta e^{\frac{b-a}{4}|x_{\mu_t}|}\pf (a+2b+\pi/4)\,,
\end{eqnarray*}
where we used that $\theta \in[0, 1]$. Moreover, by the Jensen Inequality,
\[e^{\frac{b-a}{4}|x_{\mu_t}|} \ \leqslant \ \int_E e^{\frac{b-a}{4}|x|} \mu_t(\dd x,\dd y) \ \leqslant \ \frac23 \mu_t(\V)\,. \]
We have obtained that,  
 for all $(x,y)\in E$ and all $t\geqslant 0$,
\begin{eqnarray*}
L_t \V(x,y)  & \leqslant &  - \rho \po \V(x,y) - \eta \mu_t(\V) - M  \pf  \,.
\end{eqnarray*}
with 
\[M \ =\ \rho^{-1}(a+2b+\pi/4) e^{\frac{b-a}{4}(1+R_0)}\,,\qquad \eta \ = \ 2\theta M/3\,.\]


For all $(x,y)\in E$, let  $(X_t,Y_t)_{t\geqslant0}$ be a process associated to the semi-group $(R_{s,t}^{\mu})_{t\geqslant s\geqslant 0}$  and initial conditions $(x,y)$ and let $N_t$ be its number of jumps in $[0,t]$, which is stochastically less than a r.v. with Poisson law with parameter $\|r\|_\infty t$. For $n\in \N$, let $\tau_n = \inf\{t,\ N_t \leqslant n,\ |X_t|\geqslant n\}$, which almost surely goes to infinity as $n$ goes to infinity (note that $|X_t| \leqslant |x|+ t$ for all $t\geqslant 0$). For all $n\in\N$, the Dynkin formula yields
\begin{eqnarray*}
\mathbb E\po e^{\rho(t\wedge \tau_n)} \V\po X_{t\wedge\tau_n}\pf\pf & = & \V(x) + \mathbb E\po \int_0^{t\wedge \tau_n} e^{\rho s} \po L_s\V(X_s) + \rho \V(X_s) \pf \dd s\pf\\
& \leqslant & \V(x) + \rho \int_0^t e^{\rho s} \po \eta \mu_s(\V) + M\pf \dd s\,.
\end{eqnarray*}
Letting $n$ go to infinity we obtain \eqref{Eq-Hyp-NonLin-Lyap}. If $\theta \leqslant \theta_* \leqslant 1/M$, $\eta \leqslant 2/3$. 
%

The two other Lyapunov conditions of Assumption \ref{Hyp-MeanField-Exist} are readily checked. Indeed, similar computations  yield  $L\V(x,y) \leqslant (b+\pi/2) \V(x,y)$ and then
\[P_t \V(x,y)\ \leqslant \ e^{(b+\pi/2) t} \V(x,y)\]
and, since $G_{\nu}((x,y),u)=(x,-y)$ for all $\nu\in\mathcal P(E)$, $(x,y)\in E$ and almost all $u\in[0,1]$,
\[ \frac{\V\po G_{\nu}((x,y),u)\pf }{\V\po x,y\pf } \leqslant \frac{\frac{5a+3b}{2(b-a)} + \varphi(|x|)}{\frac{5a+3b}{2(b-a)} + \varphi(-|x|)} \leqslant 2\,.\]
{Condition \eqref{Eq-NonLin-Gene-Doeblin}}  is deduced from \eqref{Eq-RTP-Doeblin}. {Finally, \eqref{Eq-Hyp-NonLin-Contract}}  ensues from \eqref{Eq-RTP-Lip} since
\[|x| \leqslant \frac{4}{b-a}  e^{\frac{b-a}{4} |x |} \leqslant \frac{8}{3(b-a)}\V(x,y)\]
for all $(x,y)\in E$. 
As a conclusion, if $\theta \leqslant \theta_* \leqslant 1/M$, Assumption \ref{Hyp-MeanField-Exist} holds, and Theorems \ref{Thm-NonLin-Exist} and \ref{Thm_non-lin} state that \eqref{Eq-RTP} admits a solution for all initial condition $m_0\in \mathcal P_\V$ and that, provided $\theta_*$ is sufficiently small (depending on $r$ or more precisely on $a$, $b$, $c$ and $\| r'\|_\infty$), then \eqref{Eq-RTP} admits a unique equilibrium toward which all solutions $t\mapsto m_t$ with $m_0\in\mathcal P_\V(E)$ converges geometrically. The equivalence of $\V$ and $\W$, hence of $\|\cdot\|_\V$ and $\|\cdot\|_\W$, concludes.

\end{proof}
\subsection{MCMC for granular media equilibrium}\label{Sec-MH}

Let $E_1= (\R/\Z)^d$ for some $d\in\N_*$ and let $U:E_1\rightarrow \R$ and $W:E_1^2\rightarrow \R$ be $\mathcal C^1$ functions respectively called the exterior and interaction potential. We want to sample according to the probability law $\mu_{V}$ on $E_1$ with density propotional to $\exp(-\beta V)$, where $\beta>0$ is the so-called inverse temperature and $V$ solves
\begin{eqnarray}\label{Eq-VPoisson-MCMC}
V(x) & = & U(x) + \frac{\int_{E_1} W(x,y) e^{-\beta V(z)}\dd  z }{\int_{E_1} e^{-\beta V(z)}\dd  z } \,.
\end{eqnarray}  
In fact, this problem doesn't necessarily  have a unique solution.  More precisely, $V$ solves \eqref{Eq-VPoisson-MCMC} iff $\mu_{V}$ is an equilibrium of the McKean-Vlasov (or granular media) equation \cite{CattiauxGuillinMalrieu,McKean,Malrieu2001}, which admits a unique such equilibrium at large temperature \cite{TugautSmallNoise,Tugaut2014}, i.e. for $\beta$ smaller than some threshold $\beta_0>0$.  We will only consider this large temperature regime.

Sampling according to $\mu_V$ can be achieved through interacting particles MCMC. To that purpose, we consider three classical Markov chains: the Metropolis-Hastings (MH) chain with Gaussian proposal, the Unadjusted Langevin Algorithm (ULA) and the Metropolis Adjusted Langevin Algorithm (MALA), which we now define. Let $N\in\N_*$ and, for $x\in E= E_1^N$, write
\[U_N(x)\ = \ \sum_{i=1}^N \po U(x_i) + \frac1{2N} \sum_{j=1}^N \1_{j\neq i} W(x_i,x_j)\pf \,.\]
For $i\in\cco 1,N\ccf$, we denote by $e_i$ the $(Nd)\times d$ matrix whose $d\times d$ blocks are all zero except the $i^{th}$   one which is the $d$-dimensional identity matrix. In other words, $e_i$ is such that, if $x\in E$ and $y\in \R^d$, then $x+e_i y = (x_1,\dots,x_{i-1},x_i+y,x_{i+1},\dots,x_N)$, where $x_i+y$ is understood as a sum in $(\R/\Z)^d$. For $i\in\cco 1,N\ccf$, $x\in E^N$ and $y\in\R^d$, set
\begin{eqnarray*}
p_{i}(x,y) & = &1 \wedge  \exp \po \beta \po U_N(x) - U_N(x+e_i (y-x_i))\pf \pf \\
\tilde p_{i}(x,y) & = &  1\wedge \frac{\exp \po \beta  U_N(x) + \frac{1}{4 \tau^2}|y  - \tau  \beta  \nabla_{x_i} U_N(x_i)|^2  \pf}{\exp  \po \beta  U_N(x+e_i y) +  \frac{1}{4 \tau^2} |-y - \tau  \beta  \nabla_{x_i} U_N(x+e_i y)|^2  \pf}
\end{eqnarray*}
for some $\tau>0$  and let $q:E_1^2 \mapsto \R_+$ be a symmetric Markov density kernel, i.e. be such that for all $x$, $q(x,\cdot)$ is the density of a probability measure on $E_1$ and $q(x,y) = q(y,x)$ for all $x,y\in E_1$. For instance, the image on the periodic torus of $q(x,y) \propto \exp(-|x-y|^2/(2\sigma^2))$   for some $\sigma>0$, or $q(x,y) \propto 1$. Define on $E$ the Markov kernels $Q_{MH}$, $Q_{ULA}$ and $Q_{MALA}$ by
\begin{eqnarray*}
Q_{MH} f(x) &=& \frac1N \sum_{i=1}^N \int_{E_1}  \po p_{i}(x,y)  f\po x+e_i (y-x_i)\pf+(1-p_{i}(x,y) ) f(x)\pf q(x_i,y)\dd y\\
Q_{ULA} f(x) &=& \frac1N \sum_{i=1}^N \int_{\R^d}   f(x+e_i y) \frac{e^{-\frac{1}{4\tau}|y-\tau  \beta  \nabla_{x_i} U_N(x_i)|^2}}{\sqrt{8\pi\tau}}\dd y\\
Q_{MALA} f(x) &=& \frac1N \sum_{i=1}^N \int_{\R^d}   \po \tilde  p_{i}(x,y)  f(x+e_i y)+(1-\tilde p_{i}(x,y) ) f(x)\pf \frac{e^{-\frac{1}{4\tau}|y-\tau  \beta  \nabla_{x_i} U_N(x_i)|^2}}{\sqrt{8\pi\tau}}\dd y
\end{eqnarray*} 
Let $x\in E$, and consider independent r.v. $Y$, $Z$, $U$ and $I$ where $Y$ follows the standard (mean zero, variance one) Gaussian distribution, and $Z$, $U$ and $I$ the uniform one respectively on $E_1$, $[0,1]$ and $\cco 1,N\ccf$. Then
\begin{eqnarray*}
x +  e_I (-\tau \beta \nabla_{x_{I}} U_N(x) + \sqrt{2\tau} Y & \sim & Q_{ULA}(x)\\
x + \1_{\{\tilde p_{I}(x,-\tau \beta \nabla_{x_{I}} U_N(x) + \sqrt{2\tau} Y)<U\}} e_I (-\tau \beta \nabla_{x_{I}} U_N(x) + \sqrt{2\tau} Y) & \sim &  Q_{MALA}(x)\\
x + \1_{\{p_{I}(x,\sigma Y)<U\}} e_I \sigma Y & \sim & Q_{MH}(x)  
\end{eqnarray*} 
if $q(x,y) \propto \exp(-|x-y|^2/(2\sigma^2))$ while, if $q(x,y) \propto 1$,
\begin{eqnarray*}
x + \1_{\{p_{I}(x,Z-x_I)<U\}} e_I (Z-x_I) & \sim & Q_{MH}(x) \,. 
\end{eqnarray*} 

More discussions, motivations and comparison of these processes can be found in \cite{BrosseDurmusMoulinesSabanis,BourabeeHairer,Durmus2018} and references within. As far as the present work is concerned, these three cases are similar, so that we focus in the following on the MH case alone, with $q(x,y)\propto 1$. For a given jump rate $\bar \lambda >0$, consider the continuous-time Markov chain on $E$  with generator
\[L_{MH}f(x) = \bar \lambda N \po Q_{MH} f(x) - f(x)\pf\,.\]
Let $(R_t)_{t\geqslant 0}$ be the associated semi-group. The associated process is constructed as follows. Let $(S_k,I_k,U_k,Z_k)_{k\geqslant 0}$ be an i.i.d. sequence where, for all $k\in\N$, $S_k$ is a standard exponential r.v., and $Z_k$, $U_k$ and $I_k$ are uniformly distributed respectively on $E_1$, $[0,1]$ and $\cco 1,N\ccf$. Let $x\in E$, set $X_{i,0} = x_i$ for all $i\in\cco 1,N\ccf$, $T_0=0$ and $T_{n+1} = T_n + (N\bar \lambda)^{-1} S_n$ for all $n\in\N$. Suppose that $(X_t)_{t\in [0,T_n]}$ has been defined for some $n\in\N$ and is independent from  $(S_k,I_k,U_k,Y_k)_{k\geqslant n}$. Set $X_t = X_{T_n}$ for all $t\in(T_n,T_{n+1})$. For all $j\neq I_n$, set $X_{j,T_{n+1}} = X_{j,T_n}$. If $U_n < p_{I,1}(X_{T_n},Z_n-X_{I,n})$, set $X_{I_n,T_{n+1}} = Z_n$ and else, set $X_{I_n,T_{n+1}} = X_{I_n,T_n}$.

\begin{prpstn}
Suppose that $\beta = \tilde \beta / \po \po osc(U)+osc(W)\pf\pf$ with $\tilde \beta$ such that 
\begin{eqnarray*}
\rho & := &    e^{-\tilde  \beta  } -   \frac{4 osc(W)}{osc(U)+osc(W)} \tilde  \beta\po e^{\tilde \beta }  - 1   \pf    \ > \ 0\,.
\end{eqnarray*}
 Then \eqref{Eq-VPoisson-MCMC} admits a unique solution $V$ and, moreover, for all $m_0,h_0 \in \mathcal P(E)$ and all $t\geqslant 0$,
 \begin{eqnarray*}
 \|m_0 R_{t} - h_0 R_{t}\|_{TV} & \leqslant &  N  e^{- \bar \lambda  \rho t}\|m_0   - h_0  \|_{TV}\,,
\end{eqnarray*}

\end{prpstn}

\begin{proof}

To recover the notations of Section \ref{Section-Def-Particule}, denote $p_*  =    \exp(-\beta(osc(U) + osc(W)) $ and for all $i\in\cco 1,N\ccf$ consider the generator $L_i$ on $E_i=E_1$ defined by 
\begin{eqnarray*}
L _if(z) & = & \bar \lambda p_* \po \int_{E_1}  \   f(y) \dd y  - f(z) \pf
\end{eqnarray*} 
for all $z\in E_1$ and bounded measurable function $f$ on $E_1$. Let $(P_t^i)_{t\geqslant 0}$ be the associated Markov semi-group, which is simply
\begin{eqnarray}\label{Eq-PtMH}
P_t^i f(z) & = & e^{-\bar \lambda p_* t} f(z) + \po 1 - e^{-\bar \lambda p_* t}\pf \int_{E_1}  \   f(y) \dd y  \,.
\end{eqnarray}
 For all $i\in\cco 1,N\ccf $ and $x\in E$, set $\lambda_i(x) = \bar \lambda (1-p_*)$ and, for a measurable bounded $f$ on $E_i$,
\begin{eqnarray*}
Q_i f(x) & = &  \int_{E_1}  \po \frac{p_{i}(x,y)-p_* }{1-p_*}  f\po y\pf+\frac{1  -p_{i}(x,y)}{1-p_*} f(x_i)\pf  \dd y\,.
\end{eqnarray*}
Remark that $p_* \leqslant p_i(x,y)$ for all $i\in\cco 1,N\ccf$, $x\in E$ and $y\in E_1$, so that $Q_i$ is indeed a Markov kernel.  Then the MH process with generator $L_{MH}$ is the system of particles such as defined in Section \ref{Section-Def-Particule} associated to  $((P_t^i)_{t\geqslant 0},\lambda_i,Q_i)_{i\in\cco 1,N\ccf }$. Note that  for all $i\in\cco 1,N\ccf$ and $x,z\in E$, 
\begin{eqnarray*}
  \|\lambda_i(x) \po Q_i (x)-\delta_x\pf  - \lambda_i(z) \po Q_i(z) -\delta_z\pf \|_{TV} &  =  &  \bar \lambda (1-p_*)  \|  Q_i (x) -Q_i(z)  \|_{TV}\,.
\end{eqnarray*}
If $x,z\in E$ are such that $x_i = z_i$, then an optimal coupling of $Q_i (x)$ and $Q_i(z) $ is constructed as follows. Let $Y$ and $U$ be independent r.v. uniformly distributed over, respectively, $E_1$ and $[0,1]$. Set
\[X_i = x_i + \1_{(1-p_*) U < p_{i}(x,Y)-p_*} (Y-x_i) \,,\qquad Z_i := z_i + \1_{(1-p_*)U < p_{i}(z,Y)-p_*} (Y-z_i)\,. \]
Then $X_i\sim Q_i(x)$, $Z_i\sim Q_i(z)$ and
\begin{eqnarray*}
\mathbb P(X_i \neq Z_i) & = & \mathbb P \co p_{i}(x,Y) \wedge  p_{i}(z,Y) -p_* < (1-p_* )U < p_{i}(x,Y) \vee  p_{i}(z,Y) -p_*\cf\\
&  =  & (1-p_*)^{-1} \int_{E_1}   | p_{i}(x,y) -p_{i}(z,y) |  \dd y \,.
\end{eqnarray*}
Now, for all $a,b\in\R$, $|1\wedge e^a-1\wedge e^b| \leqslant | e^a- e^b|\leqslant e^{b\vee a}|b-a|$, so that if $x_i=z_i$ then 
\begin{eqnarray*}
| p_{i}(x,y) -p_{i}(z,y) |  & \leqslant & \beta e^{\beta \po osc(U)+osc(W)\pf}|  U_N(x) - U_N(x+e_i y) -  U_N(z) + U_N(z+e_i y)|\\
& =  & \frac{\beta e^{\tilde \beta}}N\left |  \sum_{j\neq i} \po W(x_i,x_j) - W(x_i+y,x_j) - W(z_i,z_j) + W(z_i+y,z_j) \pf \right |\\
& \leqslant & 2 osc(W) \beta e^{\tilde  \beta  } \frac{\dist(x,y)}{N}\,.
\end{eqnarray*}
 We have thus obtained, for all $i\in\cco 1,N\ccf$ and $x,z\in E$ with $x_i=z_i$,
\begin{eqnarray*}
  \|\lambda_i(x) \po Q_i (x)-\delta_x\pf  - \lambda_i(z) \po Q_i(z) -\delta_z\pf \|_{TV} &  \leqslant &  \theta  \frac{{\dist(x,z)}}{N}\,.
\end{eqnarray*}
which is \eqref{Eq-TauxLipschitz-Particules}, with $\theta = 4\bar \lambda (1-p_*)  osc(W) \beta \exp(\tilde  \beta )$. On the other hand, the Doeblin condition \eqref{Eq_Prop-Particule_HypDoeblin_Unif} is clear, since \eqref{Eq-PtMH} immediately yields that for all $t\geqslant 0$ and $x,y\in E_1$ and all $i\in\cco 1,N\ccf$,
\[\| \delta_x P_t^i  - \delta_y P_t^i \|_{TV} \ =  \ 2 \1_{x\neq y} e^{-\bar \lambda p_* t}\,.\]
Hence, Assumption \ref{Hyp-TauxLipschitz-Particules}  holds and denoting $(R_t)_{t\geqslant 0}$ the semi-group on $E$ associated to the particle system $(X{i,t})_{i\in\cco 1,N\ccf,t\geqslant 0}$, Theorem \ref{Thm_particules_TV} means that for all $t_0>0$,  $n\in \N$ and $h_0,m_0\in \mathcal P(E)$,
\begin{eqnarray*}
 \|m_0 R_{nt_0} - h_0 R_{nt_0}\|_{TV} & \leqslant &  N e^{  \theta t_0   } \po e^{  \theta t_0  } - \po 1 - e^{-\bar \lambda p_* t_0} \pf  e^{-\bar \lambda(1-p_*) t_0 }\pf^n \|m_0   - h_0  \|_{TV}\,,
\end{eqnarray*}
which, applied for $t_0 = t/n$ for a fixed $t\geqslant 0$ and letting $n$ go to infinity, reads
\begin{eqnarray*}
 \|m_0 R_{t} - h_0 R_{t}\|_{TV} & \leqslant &  N  e^{(\theta - \bar \lambda p_*) t}\|m_0   - h_0  \|_{TV}\,.
\end{eqnarray*}
Similarly, provided $m_0 = (m_0')^{\otimes N}$ for some $m_0'\in\mathcal P(E_1)$ and similarly for $h_0$, following the remark made after Theorem \ref{Thm_particules_TV},
\begin{eqnarray*}
 \|m_t'- h_t'\|_{TV} & \leqslant &   e^{(\theta - \bar \lambda p_*) t}\|m_0'   - h_0'  \|_{TV} 
\end{eqnarray*}
holds for all $t\geqslant 0$ if $m_t'$ is  either  the law of $X_{1,t}$ or the solution of the non-linear mean-field limit of the particle system (applying Theorem \ref{Thm-TVNonLin}). In particular, remark that $V$ solves  \eqref{Eq-VPoisson-MCMC} iff the probability density proportional to $\exp(-\beta V)$ is an equilibrium of this non-linear equation. Yet, if $\rho =  p_* -\theta / \bar \lambda >0 $ then the contraction of the total variation norm implies that the latter admits a unique equilibrium, which concludes.


\end{proof}

\subsection{The Zig-Zag process with a close to tensor target}\label{Sec-ZZ}

Let $\pi\in \mathcal P(\R^N)$ be a probability law with a density proportional to $\exp(-U)$, where $U\in\mathcal C^1(\R^N)$. Consider the so-called Zig-Zag process \cite{BierkensFearnheadRoberts}  on $E=\R^N\times\{-1,1\}^N$ with generator
\[L_{ZZ} f(x,y) \ = \ y \cdot \nabla_x f(x,y) + \sum_{i=1}^N  \po y_i \partial_{x_i} U(x)\pf_+   \po f(x,y_{-i}) - f(x,y)\pf\]
where, for $y\in\{-1,1\}^N$ and $i\in \cco 1,N\ccf$, $y_{-i} = (y_1,\dots,y_{i-1},-y_i,y_{i+1},\dots,y_N)$. The Zig-Zag process admits $\pi(\dd x) \otimes \po (\delta_{-1}+\delta_1)/2\pf^{\otimes N}$ as an invariant measure and it is ergodic under general conditions on $U$ (see \cite{BierkensRobertsZitt}), so that for all reasonable (e.g. bounded) functions $f$ on $\R^N$,
\[\frac1t \int_0^t f(X_s) \dd s \ \underset{t\rightarrow \infty}\longrightarrow \ \pi(f)\,.\]
This makes it suitable for MCMC purpose. If the target measure is of tensor form, i.e. if $U(x) = U_1(x_1) + \dots + U_N(x_N)$ for some one-dimensional functions $U_i$, then the $N$ coordinates of the process are independent one-dimensional Zig-Zag processes, so that the convergence to equilibrium of each of these coordinates is independent of $N$. Let us prove that this property is stable under the addition of correlations between the coordinates, provided they are small.

\begin{prpstn}
Suppose that there exist $U_1,\dots,U_N\in \mathcal C^1(\R)$, $W\in \mathcal C^1(\R^N)$, $\rho>0$ and $R\geqslant 1$ such that for all $i\in\cco 1,N\ccf$ and $x\in \R$, if $|x|\geqslant R$ then $x U_i'(x) \geqslant \rho|x|$. Let $(R_t)_{t\geqslant 0}$ be the semi-group associated to the generator $L_{ZZ}$ with $U(x) = W(x) + \sum_{i=1}^N U_i(x_i)$. Then there exist $\theta_*,C'>0$ and $\bar\kappa\in(0,1)$ that depend only on $R$, $\rho$ and $C:= \sup\{ |U'_i(x)|,\ i\in \cco 1,N\ccf,\ |x|<R\}$ such that, if $\sup\{\|\partial_{x_i} W\|_\infty,\ i\in\cco 1,N\ccf\} < \theta_*$, then for all $\mu\in\mathcal P(\R^N)$,
\begin{eqnarray*}
\| \mu R_t - \pi\|_{TV} & \leqslant & C' \bar \kappa^{\lfloor t/t_0\rfloor} \mu(\V)\,.
\end{eqnarray*}
\end{prpstn}

\begin{proof}
Let $h,\varphi\in \mathcal C^1(\R)$ be as defined in the proof of Proposition \ref{Prop-RTP}.  For $i\in\cco 1,N\ccf$ and $(z,w)\in E =\R\times\{-1,1\}$, set $\V_i(z,w) = \exp(\rho h(z)/4 ) (1+\varphi(zw))$. We extend $\V_i$ as a function on $E^N$ by $\V_i(x,y) = \V_i(x_i,y_i)$ for $(x,y)\in E^N$ and write $\V = \sum_{i=1}^N \V_i $. Denote $\theta = \sup_{i\in\cco 1,N\ccf}\|\partial_{x_i} W\|_\infty$, and suppose that $\theta\leqslant \rho/8$.  If $(x,y)\in E^N$ is such that $|x_i|\geqslant R$,
\begin{eqnarray*}
L_{ZZ}\V_i(x,y) & = & \co \frac14 \rho   (1+\1_{x_i z_i>0})sign(x_i z_i)  - \po y_i \partial_{x_i} U(x)\pf_+ sign(x_i z_i) \cf e^{\frac\rho 2 h(x)}\\
& \leqslant & \co - \frac14  \rho   + \|\partial_{x_i} W\|_\infty\   \cf e^{\frac\rho 2 h(x)} \ \leqslant \ -\frac\rho{16}\V_i(x,y)\,.
\end{eqnarray*}
On the other hand if  $|x_i|\leqslant R$,
\begin{eqnarray*}
L_{ZZ}\V_i(x,y) + \frac\rho {16} \V_i(x,y) & \leqslant &  2 e^{\rho R/2} \po \frac{\rho}{16} + \frac{\rho}{2} + \frac{\pi}{4} + C + \|\partial_{x_i} W\|_\infty \pf \ \leqslant \ \frac\rho {16} C'
\end{eqnarray*}
with $C'=32\exp(\rho R/2)(\rho + 1 + C)/\rho $. Integrating in time similarly as in the proof of Proposition \ref{Prop-RTP}, we get that
\begin{eqnarray*}
R_t \V_i(x,y) & \leqslant & e^{-\frac\rho {16}  t}\V_i(x_i,y_i) + \po 1 - e^{-\frac\rho {16}  t}\pf C'\,.
\end{eqnarray*}
for all $t\geqslant 0$, $(x,y)\in E^N$ and $i\in\cco 1,N\ccf$. In particular, summing over $i\in\cco 1,N\ccf$,
\begin{eqnarray*}
R_t \V(x,y) & \leqslant & e^{-\frac\rho {16}  t}\V(x,y) + N \po 1 - e^{-\frac\rho {16}  t}\pf C'\,.
\end{eqnarray*}
To check the other conditions of Assumption \ref{Hyp-TauxBorne-Particules}  and \ref{Hyp-Particules-Gene}, we consider the semi-group $(P_t^i)_{t\geqslant 0}$ on $R$ with generator $L_i$ given by
\begin{eqnarray*}
L_i f(z,w) & = & w \partial_z f(z,w) + \po w U'_i(z) - \theta \pf_+ \po f(z,-w) - f(z,w)\pf\,,
\end{eqnarray*}
and the jump rates and kernels
\begin{eqnarray*}
 \lambda_i(x,y) & = & \po y_i \partial_{x_i} U(x)\pf_+  - \po y_i U'_i(x_i) - \theta \pf_+\\
 Q_i(x,y) & = & \delta_{(x,y_{-i})}\,.
\end{eqnarray*}
Then the Zig-Zag process corresponds to the process defined in Section \ref{Section-Def-Particule} associated to the semi-groups $(P_t^i)_{t\geqslant 0,i\in\cco 1,N\ccf}$ and the jump mechanims $(\lambda_i,Q_i)_{i\in\cco 1,N\ccf}$. Note that, if $G_i$ is a representation of $Q_i$, then for all $x\in E^N$ and almost all $u\in [0,1]$, $G_i((x,y),u) = (x,y_{-i})$ so that $\V_i(G_i((x,y),u)) \leqslant 2 \V_i(x_i,y_i)$. To check \eqref{Eq-Hyp-Particule-Gene-1}, through computations similar to the previous ones, it is clear that
\[L_i \V_i(z,w) \ \leqslant \ -\frac{\rho }{16}\V_i(z,w) + \tilde C\]
for some $\tilde C$ that does not depend on $\theta$. Finally, the Doeblin condition \eqref{Eq-Particule-Gene-Doeblin} is a consequence of the ergodicity of the one-dimensional Zig-Zag process as established in \cite{FontbonaGuerinMalrieu2016,BierkensRobertsZitt}. Hence, Theorem \ref{Thm_particules_gene} holds, which concludes.
\end{proof}

\subsection{Hybrid drift/bounce kinetic samplers}

Let $E_1=\R^d$ for some $d\in\N$,  $U\in\mathcal C^1(E_1)$ be an exterior potential, $W\in\mathcal C^1(E_1^2)$  an interaction potential and $\beta>0$ be an inverse temperature. Similarly as in Section \ref{Sec-MH}, we want to compute expectations with respect to the granular media equilibrium, i.e. the probability measure with density proportional to $\exp(-\beta V)$ where $V$ solves \eqref{Eq-VPoisson-MCMC}. Again, this can be approximated with a system of $N$ interacting particles. Denoting, for $x\in E:=E_1^N$,
\[U_N(x) \ = \ \sum_{i=1}^N  U(x_i) \,,\qquad W_N(x) \ =\   \frac{1}{2N} \sum_{j=1}^N W(x_i,x_j)\,,\]
consider the Markov process on $E^2$ with generator $L_N$ defined by
\begin{multline*}
L_N f(x,v) \ = \ v\cdot \na_x f(x,v) - \beta  \na U(x) \cdot \na_v f(x,y)\\
 +  \beta \sum_{i=1}^N \po v_i\cdot \na_{x_i} W_N(x)\pf_+ \po f(x,R_{N,i}(x,v)) - f(x,v)\pf + D f(x,v)
\end{multline*}
with 
\[R_{N,i}(x,v) \ = \ \po v_1,v_2,\dots, v_{i-1},v_i - 2 \frac{v_i\cdot \na_{x_i} W_N(x)}{|\na_{x_i} W_N(x)|^2}\na_{x_i} W_N(x),v_{i+1},\dots,v_N\pf\,,\]
and $D$ is some dissipativity operator on the velocities, ergodic with respect to the standard Gaussian measure $\gamma_{dN}$ on $E=\R^{dN}$, for instance $D=h D_i$, $i=1..3$ with $h>0$ and
\begin{eqnarray*}
D_1 f(x,v) & = & -v\cdot \na_v f(x,v) + \Delta_v f(x,v)\\
D_2 f(x,v) & =& \int \po f\po x,\sqrt{p} v + \sqrt{1-p} w\pf  - f(x,v)\pf \gamma_{dN}(\dd w)\\
D_3 f(x,v) & =& \sum_{i=1}^N \int \po f\po x,\po v_1,\dots,v_{i-1},\sqrt{p} v_i + \sqrt{1-p} w,v_{i+1},\dots,v_N\pf\pf - f(x,v)\pf \gamma_d(\dd w )
\end{eqnarray*}
for some $p\in[0,1)$. Then the probability measure $\mu_N \propto \exp(-\beta (U_N+W_N))$  on $E$ is invariant for $L_N$.

The motivation to use such an hybrid process to sample $\mu_N$ is the following. The computations of $\na U_N$ and $\na_x W_N$   have a respective numerical cost of $\mathcal O(N)$ and $\mathcal O(N^2)$. Suppose that $W$ is Lipschitz with a known bound $\|\na W\|_\infty \leqslant \eta$. Then, by thinning method \cite{LewisShedler,Thieullen2016}, it is rather simple to sample jump times with rate $(v_i \cdot \na_{x_i} W_N(x))_+$, by proposing jumps at rate $|v_i|\eta$ and then accepting them with some probability. In that case, we only have to compute $\na_{x_i} W$ at each of these proposed jump times, instead of computing it at all times $n\delta t$, $n\in \N$, where $\delta t$ is the timestep used for the discretization of the trajectory. If, on the other hand, no efficient bounds on $\na U$ are available, then it makes more sense to deal with $U$ with a drift operator (with a discretization scheme with time-step $\delta t$) rather than with a jump one. This argument is in fact not restricted to mean-field processes: each time the forces can be decomposed as a bounded but expensive part and a cheap but singular one (for instance, long-range interactions versus short-range interactions in molecular dynamics) it is reasonable to treat the bounded part with jumps and the singular one with drift. Contrary to multiple time-step methods like RESPA \cite{Tuckerman}, there is no contribution of the bounded forces to the systematic bias on the invariant measure.

Now, to study the long-time behavior of the particle system with generator $L_N$, we can decompose $L_N = L_N'+L_N''$ with
\begin{eqnarray*}
L_N' f(x,v) & = & v\cdot \na_x f(x,v) - \beta  \na U(x) \cdot \na_v f(x,y )+ D f(x,v) \\
L_N'' f(x,v) &=&  \beta \sum_{i=1}^N \po v_i\cdot \na_{x_i} W_N(x)\pf_+ \po f(x,R_{N,i}(x,v)) - f(x,v)\pf \,.
\end{eqnarray*}
If $D=D_1$ (resp. $D_3$) for instance, then $L_N'$ is the generator of $N$ independent kinetic Langevin  (resp. BGK-like)  processes. Remark that, in the case of $D=D_2$, then $L_N'$ is not the generator of $N$ independent processes because the partial refreshment of the velocities occurs at the same time for all the particles. As a consequence, this example doesn't exactly enter the framework of Section \ref{Section-Def-Particule}. Nevertheless, it is rather clear that, in this case, a coupling in the spirit of \cite{DurmusGuillinMonmarche2018} will still give for each coordinate, independently from the others, a probability to merge in some time $t_0$ independent from $N$. This is what is really required in the proof of Theorem \ref{Thm_particules_gene} (more precisely, of Lemma \ref{Lem-Particules-XneqY}). In fact, the partial refreshment of the velocities could also be done, for all particles at once, at a given determinist time period, as in Hamiltonian Monte-Carlo.

In any cases, as in the previous section, establishing the existence of a Lyapunov function $\mathcal \V(x) = \sum_{i=1}^N \mathcal V_i(x_i)$ and a local Doeblin condition for $L_N'$ is possible under some assumptions on $U$ and, assuming that $\eta$ is small enough, one can then prove through Theorem \ref{Thm_particules_gene} that the particle system converges towards its equilibrium at a rate independent from $N$.

\subsection{Selection/Mutation algorithms}

Let $(P_t)_{t\geqslant 0}$ be a Markov semi-group on a Polish space $E$ that satisfies Assumption~\ref{Hyp-DoeblinUnif}. Consider $\lambda_*>0$, $N\in\N_*$ and a function $p:E^2 \rightarrow [0,1]$. For all $i\in\cco 1,N\ccf$ and $x\in E^N$, set $\lambda_i(x) = \lambda_*$ and
\begin{eqnarray*}
Q_i f(x) & = & \frac1N\sum_{j=1}^N p(x_i,x_j) f(x_{j\rightarrow i}) + \po 1 - p(x_i,x_j)\pf f(x)
\end{eqnarray*}
where  $z =x_{j\rightarrow i}$ is defined by $z_k=x_k$ for all $k\in\cco 1,N\ccf\setminus \{i\}$ and $z_i = x_j$. In other words, if $x=(x_1,\dots,x_N)$ is the position of $N$ particles, a r.v. with law $Q_i(x)$ is drawn as follows: draw $J$ uniformly over $\cco 1,N\ccf$ and, with probability $p(x_i,x_J)$, kill the $i^{th}$ particle and replace it by a copy of the $J^{th}$ one. This kind of dynamics is used in a variety of algorithms, see \cite{DelMoral2013,DelMoralMiclo,CerouGuyader,FlemingViot79,CloezThai}
 and references within. Then Assumption \ref{Hyp-TauxLipschitz-Particules} clearly holds since, if $x,Z\in E^N$ are such that $x_i = z_i$,
\begin{eqnarray*}
\|Q_i(x) - Q_i(z)\|_{TV} & \leqslant & \frac1N \sum_{j=1}^N |p(x_i,x_j) - p(x_i,z_j)|\\
& \leqslant & \frac{\dist(x,z)}{N}\,.
\end{eqnarray*}
Hence, Theorem \ref{Thm_particules_TV} holds.

\subsection{The mean-field TCP process}\label{Sec-TCP}

We consider the non-linear integro-differential equation on $\R_+$ given by
\begin{eqnarray}\label{Eq-TCP}
\partial_t m_t(x) + \partial_x m_t(x,y) & = & \lambda_{m_t}(2x)m_t(2x) -  \lambda_{m_t}(x)m_t(x)\,,
\end{eqnarray}
with, for $\nu\in\mathcal P(\R_+)$ and $x\in\R_+$,
\[\lambda_\nu(x) \ = 1 + g_1(x) + \int_0^\infty g_2(x+y)\nu(\dd y)\,,\]
where $g_1$ and $g_2$ are both positive, non-decreasing functions on $\R_+$ that goes to infinity at infinity. For references on this model, called the TCP process, see \cite{ChafaiMalrieuParoux,BCGMZ} and references within. The  choice for this particular expression of $\lambda_\nu$ is not motivated by any modeling consideration; it is only meant to provide a very basic, yet interesting, example where Assumption \ref{Hyp-taux_borne} is not satisfied, since the non-linear jump rate is unbounded. In particular, it should be checked that, given a measurable function $t\mapsto \mu_t$ from $\R_+$ to $\mathcal P(\R_+)$, then the associated  process $(X_t^{\mu,x})_{t\geqslant 0}$, such as defined in Section \ref{Section-def-non-lineaire}, is well-defined for all time, namely that the probability that an infinite number of jumps occurs in a finite time is zero. 

This process is defined as follows. Suppose that $t\mapsto \mu_t $ is such that for all $x\geqslant 0$,
\begin{eqnarray}\label{Eq-TCP-gs}
g_*(x)& :=& \sup_{t\geqslant 0} \int_0^\infty g_2(x+y)\mu_t(\dd y) \ < \ \infty\,. 
\end{eqnarray}
Let $(S_k)_{k\geqslant 0}$  be an i.i.d. sequence of standard exponential r.v., $X_0 = x$ and 
$T_0=0$. Suppose that $T_n \geqslant 0$ and $(X_t)_{t\in [0,T_n]}$ have been defined for some $n\in\N$. Set
\[T_{n+1}\ = \ \inf\left  \{ t>T_n,\ S_n < \int_{T_n}^t \lambda_{\mu_s}(X_{T_n}+(s-T_n)) \dd s\right \}\,,\qquad X_{T_{n+1}} = \frac{X_{T_n}+T_{n+1}-T_n}{2}\]
and, for $t\in (T_n,T_{n+1})$, $X_t = X_{T_n}+t-T_n$. The process is then constructed by induction up to time $T_n$ for all $n\in\N$. Moreover, by construction, $X_{t\wedge T_n } \leqslant X_{0} + t$ for all $n\in\N$ and all $t\geqslant 0$, and since $g_1$ and $g_2$ are non-decreasing, $\lambda_{\mu_{s }}(X_{s}) \leqslant 1+g_1(X_{0}+t) + g_*(X_0+t)$ for all $s\leqslant t\wedge T_n$. Hence, for all $n$ such that $T_n\leqslant t$,
\[T_n \ \geqslant \ \frac{1}{1+g_1(X_{0}+t) + g_*(X_{0}+t)}\sum_{k=0}^{n-1} S_k\,.\]
In particular, for all $t\geqslant 0$, there are almost surely a finite number of jumps before time $t$. In other words, $T_n$ almost surely goes to infinity as $n$ goes to infinity, so that $X_t$ is almost surely defined for all $t\geqslant 0$. Let $(R_{s,t}^\mu)_{t\geqslant s \geqslant 0}$ be the associated inhomogeneous Markov semi-group, and $L_t$ be its generator, given by
\begin{eqnarray*}
L_t f(x) & =& f'(x) + \lambda_{\mu_t}(x) \po f\po x/2\pf - f(x)\pf\,.
\end{eqnarray*}


In the following we suppose that $g_2(x) \leqslant K\exp(\rho x)$ for some $K,\rho>0$, and that $\sup_{t\geqslant 0}\mu_t(\V)<\infty$, where $ \V(x) = \exp(\rho x)$ . In particular,  \eqref{Eq-TCP-gs} holds, so that the associated process and semi-group are well defined.

 Then, for $x\geqslant R:= \inf\{x\geqslant 0,\ g_1(x) \geqslant 2\rho\}$,
\begin{eqnarray*}
L_t \V(x) & = & \rho \V(x) + \po 1 + g_1(x) + \int_0^t g_2(x+y)\mu_t(\dd y)\pf \po \V \po \frac x2\pf - \V(x)\pf\\
& \leqslant & (1+2\rho) \sqrt{ \V (x)} - (1+\rho) \V(x) \ \leqslant \ - \frac12 \V(x) + \frac12 + \rho \,.
\end{eqnarray*}
For $x\in[0,R]$, $L_t \V(x) \leqslant \rho \V(x) \leqslant \rho \exp(\rho R)$ and thus, for all $x\in\R_+$,
\[L_t \V(x) \ \leqslant \ - \frac12 \po  \V(x) - \tilde C  \pf \,,\]
with $\tilde C  =\po 1 + 2\rho\pf \po 1 + e^{\rho R}\pf$. As in the proof of Proposition \ref{Prop-RTP}, this yields
\begin{eqnarray}\label{Eq-TCP-Lyap}
R_{0,t}^{\mu}\V(x) & \leqslant & e^{-t/2} \V(x) + \po 1 - e^{-t/2}\pf \tilde C \,.
\end{eqnarray}
In particular, for any $m_0\in\mathcal P_\V(\R_+)$,  $t\geqslant 0$ and $x\geqslant 0$,
\[m_0 R_{0,t}^\mu \po g_2(x+\cdot)\pf  \  \leqslant\ K e^{\rho x} m_0 R_{0,t}^\mu (\V)\  \leqslant \  K e^{\rho x}\po  \tilde C \vee m_0(\V)\pf\,.\]
Note that, simply by changing $\rho$ to $2\rho$, the same computations  shows that there also exist $\hat C$ such that
\begin{eqnarray}\label{Eq-TCP-Lyap2}
R_{0,t}^{\mu}(\V^2)(x) & \leqslant & e^{-t/2} \V^2(x) + \po 1 - e^{-t/2}\pf \hat C \,.
\end{eqnarray}

 Besides, for all $\nu_1,\nu_2\in\mathcal P(\R_+)$, denoting $(Y,\tilde Y)$ an optimal coupling of $\nu_1$ and $\nu_2$ such as given by Lemma \ref{Lem-couplage},
 \begin{eqnarray*}
  |\lambda_{\nu_1}(x) - \lambda_{\nu_2}(x)| & = & |\mathbb E \po g_2(x+Y) - g_2(x+\tilde Y)\pf| \\
  & \leqslant & \mathbb E \po \1_{Y\neq \tilde Y} \po g_2(x+Y) +g_2(x+ \tilde Y)\pf\pf \ \leqslant \ \frac{e^{\rho x}}K\|\nu_1-\nu_2\|_\V\,.  
 \end{eqnarray*}
 Let $m_0\in\mathcal P_\V(\R_+)$ and $t\mapsto \mu_t^1,\mu_t^2$ be such that $\mu_t^i(\V) \leqslant m_0(\V) \vee \tilde C$ for $i=1,2$ and $t\geqslant 0$. Consider the synchronous coupling $(X_t,\tilde X_t)$ of $m_0 R_{0,t}^{\mu^1}$ and $m_0 R_{0,t}^{\mu^2}$ with $X_0=\tilde X_0\sim m_0$, such as defined in Section \ref{Sec-Preuves-NL}. Then
 \begin{eqnarray*}
 \| m_0 R_{0,t}^{\mu^1} - m_0 R_{0,t}^{\mu^2}\|_\V & \leqslant & \mathbb E\po \1_{X_t \neq \tilde X_t}\po \V(X_t)+ \V(\tilde X_t)\pf \pf \\
 & \leqslant & 2 \mathbb E\po e^{\rho(X_0+t)} \mathbb P \po X_t\neq \tilde X_t\ |\ X_0\pf\pf\,.
 \end{eqnarray*}
 Denoting $\tau_{split} = \inf\{t\geqslant 0,\ X_t\neq \tilde X_t\}$, remark that, from $X_t \leqslant X_0+t$, almost surely, for all $s\in[0, \tau_{split})$,
 \begin{eqnarray}
 |\lambda_{\mu_s^1}(X_s) - \lambda_{\mu_s^2}(X_s)| & \leqslant &  \frac{e^{\rho (X_0+s)}}K\|\mu_s^1-\mu_s^2\|_\V\,,
 \end{eqnarray}
 and then
 \begin{eqnarray*}
 \mathbb P \po X_t\neq \tilde X_t\ |\ X_0\pf &=&  \mathbb P \po \tau_{split} \leqslant t\ |\ X_0\pf\\
 & \leqslant & 1 - \exp \po - \frac{e^{\rho (X_0+t)}}K\int_0^t \|\mu_s^1-\mu_s^2\|_\V\dd s\pf  \\
 & \leqslant & \frac{e^{\rho (X_0+t)}}K\int_0^t \|\mu_s^1-\mu_s^2\|_\V\dd s\,,
 \end{eqnarray*}
 so that, at the end,
  \begin{eqnarray}\label{Eq-TCP-3}
\notag  \| m_0 R_{0,t}^{\mu^1} - m_0 R_{0,t}^{\mu^2}\|_\V & \leqslant & \mathbb E\po \1_{X_t \neq \tilde X_t}\po \V(X_t)+ \V(\tilde X_t)\pf \pf \\
 & \leqslant & \frac{2}K e^{\rho t} m_0(\V^2)\int_0^t \|\mu_s^1-\mu_s^2\|_\V\dd s\,.
 \end{eqnarray}
 Similarly as in the proof of Lemma \ref{Thm-NonLin-Exist}, taking $t_1>0$ small enough so that $ 2 t_1  \exp(\rho t_1) [m_0(\V^2)\vee \hat C] / K \leqslant 1/2$, the map $(\mu_t)_{t\in[0,t_1]} \rightarrow (m_0 R_{0,t}^\mu)_{t\in[0,t_1]}$ admits a unique fixed point, which we call by definition a solution of \eqref{Eq-TCP}. Moreover, from \eqref{Eq-TCP-Lyap2}, $m_{t_1}(\V^2) \leqslant m_0(\V^2)\vee \hat C$, so that the same argument (with the same $t_1$) yields a solution on $[t_1,2t_1]$, and then on all $\R_+$.  From 
\eqref{Eq-TCP-Lyap}, such a solution $t\mapsto m_t$ satisfies 
\begin{eqnarray}\label{eq-TCP-4}
 \notag m_t(\V)  & \leqslant & e^{-t/2} m_0(\V) + \po 1 - e^{-t/2}\pf \tilde C\\
\Rightarrow\qquad \lambda_{m_t}(x) & \leqslant & 1 + g_1(x) + K e^{\rho x} \po m_0(\V) + \tilde C\pf\,,
\end{eqnarray}
for all $t,x\geqslant 0$. This bound is now uniform in time, although not uniform in $x$.

For $t$ large enough, from \eqref{Eq-TCP-Lyap2}, $m_t(\V^2) \leqslant \hat C +1$, so we can now restrict the study to initial distributions that satisfies $m_0(\V^2) \leqslant \hat C +1$. Then, with the same argument used to establish \eqref{Eq-TCP-3}, there exists $C''>0$ such that for all $m_0,h_0$ with $m_0(\V^2),h(\V^2) \leqslant \hat C +1$ and all $t_0\geqslant 0$,
  \begin{eqnarray*}
 \| m_t - h_t\|_\V & \leqslant & C'' \int_0^t \|m_s-h_s\|_\V\dd s\,,
 \end{eqnarray*}
 which is the ``splitting estimates'' part of the proof of Theorem \ref{Thm_non-lin}. Since the Lyapunov contraction \eqref{Eq-TCP-Lyap} has already been established, what remains to obtain the ``merging estimates'' part, i.e. the equivalent of Lemma \ref{Lem_NL_HairerMatt}, is just a time $t_0$ and a probability $\alpha>0$ to merge two processes in a time $t_0$, given that they started in some compact set. For the linear processes (i.e. with $g_2 = 0$), this has been done in \cite{ChafaiMalrieuParoux,BCGMZ}. As in Section \ref{Section-Couplage-NL}, to couple two non-linear processes, we simply couple two linear processes (i.e. with jump rate $1+g_1$) and we hope that no non-linear jump occurs in the time interval $[0,t_0]$. In the present case, however, we should be cautions, since the non-linear jump rate is not uniformly bounded in $x$. Nevertheless, if the starting point are taken in $[0,R_0]$ for some $R_0>0$, then during the time interval $[0,t_0]$ the processes remain in $[0,R_0+t_0]$, on which the non-linear jump rate is bounded by \eqref{eq-TCP-4}. Thus, conditionnally to the fact that they started in $[0,R_0]$ the probability that two coupled non-linear processes with initial law $m_0,h_0$ that satisfy $m_0(\V^2),h(\V^2) \leqslant \hat C +1$ have merged at time $t_0$ is uniformly bounded away from 0. Then the strategy of the proof of Theorem \ref{Thm_non-lin} can be adapted to get, provided that $K$ is small enough, that \eqref{Eq-TCP} admits a unique equilibrium, toward which all solutions with $m_0\in\mathcal P_\V$ converges exponentially fast.

\subsection{Markov processes with delay}\label{Sec-Delay}

We claimed in the introduction that the coupling arguments  used in \ref{Sec-Preuves} to deal with non-linear equation can also be used to deal with self-interacting processes. A general study of processes $(X_t)_{t\geqslant 0}$ interacting with a weighted empirical distribution
\[\nu_t \ = \ \po \int_0^t  w(t-s) \dd s \pf^{-1} \int_0^t \delta_{X_s} w(t-s)\dd s\]
for all $t\geqslant 0$ such that $\int_0^t  w(t-s) \dd s \neq 0$ exceeds the scope of the present paper, and is thus postponed. Nevertheless, as a proof of principle, consider the simple case where $w$ is a Dirac mass at some given time $t_0$. In other words, $X_t$ follows some Markovian dynamics and jumps at a rate and to a position that  depends on $X_{t-t_0}$. In particular, $Y_t=(X_s)_{s\in[t-t_0,t]}$ is a Markov process. Nevertheless it is not necessary to consider the somewhat complicated process $(Y_t)_{t\geqslant 0}$ to obtain a speed of convergence toward equilibrium for the law of $X_t$.  If there is some probability to couple two non-delayed processes starting at different positions and if the delayed jump rate is bounded, then there is some probability to merge two  delayed processes at some positive time $t_1$ and then there is some probability that no delayed jump occurs in a time interval of length $t_0$. If that happens, then after time $t_0+t_1$ the two processes share the same position and the same memory, and they will stay equal forever. An exponential convergence toward some equilibrium is obtained for the law of $X_t$ (which is not the solution of an integro-differential equation).

\section*{Acknowledgements}
Pierre Monmarch\'e acknowledges support from the French ANR grants EFI (ANR-17-CE40-0030) and SWIDIMS (ANR-20-CE40-0022).


\bibliographystyle{plain}
\bibliography{biblio}

\end{document}